\documentclass[11pt]{amsart}
\usepackage{amsmath,amsxtra,amssymb}
\numberwithin{equation}{section}
\newtheorem{lemma}{Lemma}[section]
\newtheorem{prop}[lemma]{Proposition}
\newtheorem{theorem}[lemma]{Theorem}
\newtheorem{cor}[lemma]{Corollary}

\newtheorem{prob}[lemma]{Problem}
\newtheorem{rem}[lemma]{Remark}

\newcommand{\re}{\begin{rem}\rm}
  \newcommand{\mar}{\end{rem}}
\newtheorem{exam}[lemma]{Example}

\newcommand{\kla}{\left ( }
\newcommand{\mer}{\right ) }

\renewcommand{\for}{\begin{eqnarray*}}

\newcommand{\mel}{\end{eqnarray*}}

\newcommand{\mitt}{\left | { \atop } \right.}
\newcommand{\kl}{\pl \le \pl}
\newcommand{\gl}{\pl \ge \pl}

\newcommand{\lel}{\pl = \pl}

\newcommand{\ez}{{\mathbb E}}

\newcommand{\nz}{{\rm  I\! N}}
\newcommand{\nen}{n \in \nz}

\newcommand{\rz}{{\mathbb R}}

\newcommand{\cz}{{\mathbb C}}
\newcommand{\kz}{{\rm  I\! K}}

\newcommand{\ten}{\otimes}

\newcommand{\wet}{\stackrel{\wedge}{\otimes}}

\newcommand{\p}{\hspace{.05cm}}
\newcommand{\pl}{\hspace{.1cm}}
\newcommand{\pll}{\hspace{.3cm}}

\newcommand{\hz}{\vspace{0.5cm}}
\newcommand{\lz}{\vspace{0.3cm}}

\newcommand{\qd}{\end{proof}\vspace{0.5ex}}

\newcommand{\Om}{\Omega}
\newcommand{\om}{\omega}

\newcommand{\al}{\alpha}
\newcommand{\si}{\sigma}
\newcommand{\Si}{\Sigma}

\newcommand{\la}{\lambda}
\newcommand{\eps}{\varepsilon}

\renewcommand{\L}{{\mathcal L}}
\newcommand{\F}{{\mathcal F}}
\newcommand{\E}{{\mathcal E}}

\newcommand{\K}{{\mathcal K}}

\newcommand{\M}{{\mathcal M}}

\newcommand{\N}{{\mathcal N}}

\newcommand{\U}{{\mathcal U}}
\newcommand{\B}{{\mathcal B}}
\newcommand{\Hh}{{\mathcal H}}
\newcommand{\ff}{{\mathbb  F}}

\newcommand{\noo}{\left \|}
\newcommand{\rrm}{\right \|}

\newcommand{\bet}{\left |}
\newcommand{\rag}{\right |}

\newcommand{\intt}{\int\limits}
\newcommand{\summ}{\sum\limits}

\newcommand{\prodd}{\prod\nolimits}

\newcommand{\lb}{\langle \langle }
\newcommand{\rb}{\rangle \rangle }

\newcommand{\8}{\infty}

\newcommand{\card}{\operatorname{card}}

\allowdisplaybreaks[3]
\begin{document}

\title[Embedding of OH]
{\bf \large Embedding of the operator space OH and the logarithmic
`little Grothendieck inequality'}

%\begin{center}
%\bf \large
% The projection constant of \boldmath $OH_n$ \unboldmath
% and the
%\end{center}
%\begin{center}
%\bf \large  logarithmic `little Grothendieck
%inequality'
%\end{center}

\author[Marius Junge]{Marius Junge$^\dag$}
\address{Department of Mathematics\\
University of Illinois at  Urbana-Champaign, Urbana, IL 61801}
\email{junge@math.uiuc.edu}

\thanks{$^\dag$Junge is  partially supported by the
National Science Foundation DMS-0301116}
%\begin{center}
%Marius Junge$^\S$
%\end{center}

\begin{abstract} Using free random varaibles we find an
embedding of the operator space $OH$ in the predual of a von
Neumann algebra. The properties of this  embedding allow us to
determined the projection constant of $OH_n$, i.e. there exists a
projection $P:B(\ell_2)\to OH_n$ such that
 \[  \frac{1}{96} \pl \sqrt{\frac{n}{1+\ln n}} \kl  \noo P\rrm_{cb} \kl 144 \pi \sqrt{\frac{n}{1+\ln n}} \pl .\]
According to recent results  of Pisier/Shlyahtenko, the lower
bound holds for every projection.  Improving a previous estimate
of order $(1+\ln n)$  of the author, Pisier/Shlyahtenko obtained
the following `logarithmic little Grothendieck inequality'
 \[ \kla \summ_{k=1}^n \noo u(x_k)\rrm^2 \mer^{\frac12} \kl 96
 \sqrt{1+\ln n}\pl \noo u\rrm_{cb} \pll \noo \summ_{k=1}^n x_k \ten \bar{x}_k
 \rrm_{B(H)\ten_{min}\bar{B}(H)}^{\frac12} \pl .\]
which holds for all linear maps $u:B(\ell_2)\to OH$ and vectors
$x_1,...,x_n$. We find a second proof of this inequality which
explains why the factor $\sqrt{1+\ln n}$ is indeed necessary. In
particular the operator space version of the `little Grothendieck
inequality' fails to hold. This `logarithmic little Grothendieck'
inequality characterizes $C^*$-algebras with the weak expectation
property of Lance.
 \end{abstract}
\maketitle

\newcommand{\llz}{\vspace{0.15cm}}

{\bf Plan:}
\begin{enumerate}
 \item[0.] Introduction and Notation.\llz
 \item[1.] Preliminaries.\llz
 \item[2.] A logarithmic characterization of C$^*$-algebras with WEP.\llz
 \item[3.] Pusz/Woronowicz' formula and the operator space $OH$. \llz
 \item[4.] The projection constant of the operator space $OH_n$. \llz
 \item[5.] Norm calculations in a quotient space.\llz
 \item[6.] $K$-functionals associated to a states and conditional expectations.\llz
 \item[7.] Sums of free mean zero variables. \llz
 \item[8.] Appendix: Concrete realization of $OH$ as functionals.\llz
\end{enumerate}

\newpage

\section{Introduction and Notation}
\setcounter{section}{0} \setcounter{equation}{0}

Probabilistic techniques  and concepts play an important role in
the theory of Banach spaces and operator algebras. For example
Khintchine's inequality and its application  to absolutely summing
sequences in $L_p$ spaces, the Grothendieck inequality and   the
`little Grothendieck inequality' are classical, fundamental tools
in the theory of Banach spaces. An analogue of Grothendieck's
work in the context of operator algebras and operator spaces, the
so-called `Grothendieck program', inspired important results in
operator algebras and motivates new research in operator spaces.
Let us mention Connes' \cite{Co} characterization of injective
von Neumann algebras, Pisier/Haagerup's noncommutative version of
Grothendiek's inequality (see e.g. \cite{Psl}), a noncommutative
Grothendieck inequality for exact operator spaces in \cite{JP}
and very recently Pisier and Shlyakhtenko's  \cite{PS}
Grothendieck theorem for operator spaces. In our  context
Effros/Ruan's \cite{ER0} notion of $1$-summing maps in operator
spaces is a further important example of this concept.

In this paper, we show that the operator space $OH$ introduced by
Pisier also admits a probabilistic realization in the predual of
a von Neumann algebra. The classical gaussian random variables
and semicircular random variables in a tracial probability space
are too commutative and have to be replaced by linear combinations
of free semicircular random variables in a non-tracial setting.
Indeed, after the first draft of this paper circulated Pisier
\cite{Ps5} shows the surprising fact that $OH$ does embed in the
predual of tracial probability spaces. Using the non-tracial
embedding, we follow the Banach space approach  and are then able
to determine the projection constant of the operator space $OH_n$
by calculating the norm of a canonical tensor in the operator
space projective tensor product. This is intimately connected to
the `little Grothendieck inequality' for operator spaces and we
can determine the right order for both inequalities.

The `little Grothendieck inequality' can be derived from
Grothendieck's inequality but also admits a very simple
probabilistic proof, whose roots might already been known to
Orlicz. Indeed, by duality it is sufficient to consider a linear
map $u:\ell_2^n\to L_1[0,1]$. Using Rademacher variables
$(\eps_i)$, we deduce for vectors $x_1,...,x_m\in \ell_2^n$ that
 \begin{eqnarray}
 \begin{minipage}{14cm}
 \for
   \noo \kla \summ_{i=1}^m |u(x_i)|^2 \mer^{\frac12}\rrm_1
   &\le& \sqrt{\frac{\pi}{2}} \pl \ez \noo \summ_{i=1}^m \eps_i
   u(x_i) \rrm_1  \nonumber \kl
  \sqrt{\frac{\pi}{2}} \pl  \noo u\rrm
   \ez \noo \summ_{i=1}^m \eps_i
   x_i \rrm_2  \nonumber \\
   &\le& \sqrt{\frac{\pi}{2}} \pl  \noo u\rrm \pl \kla
   \summ_{i=1}^m \noo x_i\rrm_2^2 \mer^{\frac12} \pl .
  \mel \end{minipage}
   \end{eqnarray}
Then duality yields the `little Grothendieck inequality' (for
Banach spaces). This proof inspired  our approach to the  operator
space version of the `little Grothendieck inequality'. The
strategy suggested by this approach  is to find a good embedding
of $OH_n \subset L_1(N)$ (using noncommutative gaussian random
variables) and then compare the norm of matrices in $OH_n\ten
OH_n\subset L_1(N\ten N)$ with the norm in $\ell_2^{n_2}$.
Following the commutative philosophy a good embedding (with
`gaussian' variables) will automatically provide the `little
Grothendieck inequality'.

Let us recall basic operator space notations in order to make
these statements more precise. An operator space $F$ comes either
with a concrete isometric embedding $\iota:F\to B(H)$ or with a
sequence $(\noo \pl \rrm_m)$ of matrix norms on $(M_m(F))$
derived from this embedding such that
 \[ \noo [x_{ij}] \rrm_{M_m(F)} \lel \noo
 [\iota(x_{ij})]\rrm_{B(H^m)} \pl .\]
Ruan's axioms tell us that whenever we have a sequence $(\noo \pl
\rrm_m)$ on $(M_m(F))$
 \[ \noo \left[ \begin{array}{cc} x& 0 \\ 0 & y \end{array}\right]
 \rrm_{n+m} \lel \max\{ \noo x\rrm_{n},\noo y\rrm_{m}\} \quad
 \mbox{and} \quad
 \noo [\summ_{k,l} a_{il}x_{kl}b_{lj}]\rrm_{m} \kl \noo
 a\rrm \noo x \rrm_{m} \noo b\rrm \]
then this can be obtained from some concrete embedding
$\iota:F\to B(H)$. Since a $C^*$-algebra $A$ admits a
representation  as operators on a Hilbert space and the norms of
the $C^*$-algebra  $M_m(A)$ are uniquely determined by this
embedding, it carries a natural operator space structure. In the
theory of operator spaces, hilbertian operator spaces are of
particular interest. For example  the column and row spaces of
matrices
 \[ H^c\lel B(\cz,H) \pl \subset B(H) \quad \mbox{and}
 \quad  H^r \lel B(H,\cz) \subset B(H) \pl \]
play a fundamental r\^{o}le.  Pisier discovered that the sequence
of norms  on $M_m\ten \ell_2$ obtained by the complex
interpolation method
 \[ M_m(OH) \lel [M_m(\ell_2^c),M_m(\ell_2^r)]_{\frac12} \]
defines a sequence of matrix norms on $M_m\ten \ell_2$ satisfying
Ruan's axioms. By Ruan's theorem these interpolation norms are
induced by a  sequence of operators $(T_k)\subset B(\ell_2)$
satisfying
 \[ \noo \summ_{k=1}^\8 a_k \ten T_k \rrm_{M_m(B(\ell_2))}
  \lel \noo \summ_{k=1}^\8 a_k\ten \bar{a}_k \rrm_{M_m\ten
  M_m}^{\frac12}
  \lel \sup_{\noo a\rrm_2,\noo b\rrm_2\le 1}
  \kla \summ_{k=1}^\8 tr(ax_kbx_k^*)\mer^{\frac12} \pl
  \]
for all sequences $(a_k)\in M_m$. The {\it operator Hilbert
space} $OH$ is defined to be the span of the $(T_k)$'s. It is
still unclear how to construct `concrete' operators $(T_k)$
satisfying this equality, but it is our aim to shed some light on
this question. In view of the Khintchine's inequality, it maybe
more natural to expect $OH$ as a subspace of a noncommutative
$L_1$ spaces, rather than finding directly the appropriate
embedding in a noncommutative $L_\8$ space. Since operator spaces
are closed under taking dual spaces, noncommutative $L_1$ spaces
carry a natural operator space structure which we now recall. The
morphisms in the category of operator space are {\it completely
bounded maps} $u:E\to F$ where
 \[ \noo u\rrm_{cb} \lel \sup_m \noo id_{M_m}\ten u: M_m(E) \to
 M_m(F)\rrm \pl .\]
The space of completely bounded maps is denoted by $CB(E,F)$. Then
the {\it standard dual} $E^*$ of an operator space $E$ is given by
the matrix norms on $M_m(E^*)$ by
 \[ M_m(E^*)\lel CB(E,M_m) \quad \quad \quad  \mbox{isometrically.} \]
Indeed, a matrix of functionals $[x^*_{ij}]$ induces a natural
operator $u:E\to M_m$ defined by $u(e)=[x^*(e)]_{ij}$. Note that
$OH$ is the only selfdual operator space, this means  the mapping
$id:OH\to OH^*$, $id(e_k)=e_k^*$ has cb-norm $1$. Here, $(e_k)$
$(e_k^*)$ denotes the sequence of the standard unit vector basis
and their biorthogonal functionals, i.e.   $e_k^*(e_j)=
\delta_{kj}$ (see \cite{Poh}). Motivated by the commutative
approach we could show in  preliminary version of this paper that
there is a linear map $u:OH_n\to A^*$ such that
 \for
 \frac{1}{C\p(1+\log n)} \noo x\rrm_{M_m(OH_n)} &\le&
 \noo (id_{M_m} \ten u)(x)\rrm_{M_m(A^*)}
    \kl C (1+\log n) \p \noo x\rrm_{M_m(OH_n)}
   \pl .
  \mel
Moreover,  for every $u:B(H)\to OH$
  \begin{eqnarray}
  \kla \summ_{k=1}^n \noo u(x_k)\rrm^2 \mer^{\frac12}
  &\le& C (1+ \ln n) \pl \noo u\rrm_{cb} \noo \summ_{k=1}^n x_k\ten
  \bar{x}_k\rrm^{\frac12} \pl .
  \end{eqnarray}
for all $x_1,...,x_n\in B(H)$. Pisier/ Shlyakhtenko's \cite{PS}
Grothendieck  theorem  for operator space implies that a map
$u:B(H)\to OH$ is completely bounded if and only if there  exists
a state $\phi$ and a constant $C$ such that
 \begin{eqnarray}
   \noo u(x)\rrm &\le&  C \pl [\phi(x^*x)\phi(xx^*)]^{\frac14}
 \end{eqnarray}
for all $x\in B(H)$. This can be understood as an the
interpolation inequality with respect to real interpolation
method and index $(\frac12,1)$. By a well-known application of
the Grothendieck-Pietsch factorization technique \cite{Pip}, the
operator space `little Grothendieck inequality' provides an
interpolation inequality with respect to the complex method with
index $\frac12$. Pisier/ Shlyakhtenko's Grothendieck \cite{PS}
derived from $(0.3)$ a logarithmic version
  \begin{eqnarray}
 \kla \summ_{k=1}^n \noo u(x_k)\rrm^2 \mer^{\frac12}
  \kl C \pl \sqrt{1+ \ln n} \pl \noo u\rrm_{cb} \noo \summ_{k=1}^n x_k\ten
  \bar{x}_k\rrm^{\frac12} \pl .
  \end{eqnarray}
In this paper, we will provide a second proof of $(0.4)$ and an
embedding of $OH$.  A key new ingredient (compared to our
preliminary version)  is a formula of Pusz/Woronowicz (and its
new dual version)  for the root of sesquilinear forms, see
\cite{Wo}. We refer to section 3 for more details on this
formula, but we will  indicate here  its consequences for $OH$.
Following Pusz/Woronowicz we use the the probability measure
$d\mu(t)=\frac{dt}{\pi \sqrt{t(1-t)}}$ on $[0,1]$ and the two
measures $d\nu_1(t)=t^{-1}d\mu(t)$,
$d\nu_2(t)=(1-t)^{-1}d\mu(t)$. Then the direct sum
 \[ E \lel L_2^c(\nu_1;\ell_2) \oplus L_2^r(\nu_2;\ell_2) \]
is an operator space. On $E$, we define the map $Q:E\to
L_0(\mu;\ell_2)$ by
 \[ Q(x_1,x_2)(t) \lel x_1(t)+x_2(t) \in \ell_2 \pl .\]
Using Ruan's theorem, we know that operator spaces  are closed by
taking quotients. Indeed, the matrix norm on the {\it operator
space quotient} $E/S$ is given by
 \[ M_m(E/S) \lel M_m(E)/M_m(S) \]
which is easily seen to  satisfy Ruan's axioms. Therefore, we
obtain an operator space
 \[ G \lel E/kerQ \pl .\]
Then, we may consider the subspace $F\subset G$ of equivalence
classes $(x_1,x_2)=kerQ$ such that $x_1+x_2$ is a $\mu$-a.e. a
constant element in $\ell_2$.

\hz \noindent {\bf Theorem 1.} {\it  $F$ is $2$ completely
isomorphic to $OH$.} \hz

Using Voiculescu's concept of free probability and an operator
valued  extension  of  Voiculescu's inequality for sums of
independent random variables \cite{Voir}, we obtain free
noncommutative gaussian variables in a non-tracial situation and
an embedding of $G$ in the predual of a von Neumann algebra. In
the appendix we provide a concrete realization of the underlying
von Neumann algebra $N$, a type III$_1$ factor with a free
quasi-free state in the sense of Shlyakhtenko \cite{S1} obtained
from $\mu$, $\nu_1$ and $\nu_2$. This concrete realization
provides linear combinations of semi-circular random variables
representing $OH$ when restricted to $N$. They are obtained from
a central limit procedure which justifies the term free
non-tracial free gaussian variables. For our applications to the
projection constant, we also have to know that the underlying von
Neumann algebra is QWEP. Let us recall that a $C^*$-algebra has
the weak expectation property (WEP) if there exists  a complete
contraction $P:B(H)\to A^{**}$ such that $P|_{A}=id_A$. A
$C^*$-algebra is QWEP if it is the quotient $A=B/I$ of  a
$C^*$-algebra with WEP by a two-sided ideal. It is an open
problem whether every $C^*$-algebra is QWEP, see \cite{Ki}.

\hz

\noindent  {\bf Theorem 2.} {\it $G$ completely isomorphic to a
complemented complemented subspace of the predual of  von Neumann
algebra with QWEP.}

\hz

Since $OH$ is completely isomorphic to a subspace of $G$,  we
immediately obtain the non-commutative analogue of the gaussian
embedding of $\ell_2$.

\hz  \noindent  {\bf Corollary 3.} {\it $OH$ completely embeds
into the predual of a von Neumann algebra $N$ with QWEP.}

\hz

In view of the discussion after $(0.1)$, we should expect that a
`good' embedding of $OH$ provides the `little Grothendieck
inequality'. At this point the analogy breaks down, or more
precisely produces different norms than in the classical case.
Indeed, the norm of matrices on $OH_n\ten OH_n$ viewed as a
subspace of $L_1(N\ten N)$ turns out to be too big.  This can be
expressed using  basic operator space tensor norms. The {\it
injective tensor norm} between tow operator spaces $E\subset
B(H)$, $F\subset B(K)$ is defined by
 \[  E\ten_{min} F \subset B(H\ten K) \pl .\]
Note that for  finite dimensional operator spaces $E$ or $F$ we
have $E^*\ten_{min}F=CB(E,F)$. The {\it projective tensor norm}
$E\wet F$ is defined such that
 \[ (E\wet F)^* \lel CB(E,F^*) \cong CB(F,E^*) \]
holds isometrically. The {\it operator $1$-summing norm} of a
linear map $u:E\to F$ is defined as
 \[ \pi_1^o(u) \lel \noo id \ten u: \K^*\ten_{min} E\to \K^*\wet F \rrm \pl .\]
In this definition Effros and Ruan followed Grothendieck's
approach replacing the space $c_0$ by the quantized analogue
$\K$, the compact operators on $\ell_2$.

\hz  \noindent  {\bf Theorem 4.} {\it Let $(f_i)$ the canonical
unit vector basis in $F$. Let $u:OH\to OH$ be a linear map
represented by the matrix $[a_{ij}]$, then}
 \[ \frac{1}{C} \pl  \pi_1^o(u) \kl \noo \summ_{ij}
 a_{ij} f_i\ten f_j \rrm_{G\wet G}
 \kl C \pl \pi_1^o(u) \pl .\]

The proof of Theorem 4 uses the complementation from Theorem 2,
an embedding  of $G\wet G$ in $M_*\wet \M_*$ and the QWEP
property of $\M_*$. The final punch line are norm calculations in
$G\wet G$, the quotient of (four) classical Banach spaces, see
section 5.

\hz \noindent {\bf Proposition 5.} {\it Let $\nen$, then }
 \[ \frac{1}{8\pi}  \pl \sqrt{n(1+\ln n)} \kl \noo \summ_{i=1}^n f_i\ten f_i\rrm_{G\wet G}
 \kl 16  \pl \sqrt{n (1+\ln n)} \pl .\]

A classical application of trace duality relates the `little
Grothendieck inequality' to the projection constant of $OH_n$.

\hz \noindent {\bf Theorem 6.} {\it There exists a constant $C\ge
1$ such that for all $\nen$, there is a projection
$P:B(\ell_2)\to OH_n$ such that}
 \[ \frac{1}{96}\pl \sqrt{ \frac{n}{1+\ln n}} \kl \noo P\rrm_{cb} \kl 144 \pi  \pl \sqrt{ \frac{2n}{1+\ln n}}
 \pl . \]\lz

A lower estimate for the $cb$-norm of this order for any
projection $P$ follows from the logarithmic 'little Grothendieck
inequality' (0.4).

\hz  \noindent {\bf Corollary 7.} {\it The logarithmic order
$\sqrt{1+\ln n}$ in $(0.4)$ is best possible.} \hz

Finally, let us discuss the class of $C^*$-algebras which admits
the logarithmic `little Grothendieck inequality'. C. le Merdy
observed that it can not hold for the $C^*$-algebra generated by
the left regular representation of the free groups in countably
many generators. Using a (unfortunately unpublished) result of
Haagerup, see \cite{Haa}, we obtain the following
characterization.

\hz \noindent {\bf Theorem 8.} {\it A $C^*$ algebra $A$ satisfies
the logarithmic `little Grothendieck inequality $(0.2)$ or $(0.4)$
if and only if $A$ is WEP.}\hz

The paper is organized as follows. Preliminary facts, including
the operator space structure on $L_1(N)$, are presented in
section 1. In section 2, we prove Theorem 8. modulo $(0.4)$. We
investigate the Pusz/Woronowicz's formula and its consequences
assuming Theorem 2 in section 3. In section 4, we use Theorem 2
and Theorem 4 in showing that the logarithmic `little
Grothendieck inequality' is optimal. The norm estimates for the
upper and lower bound of Proposition 5 can be found in section 5.
K-functionals and quotients similar as  $G$ are a central object
in our paper and are discussed in detail in section 6. We provide
the probabilistic tools from free probability in section 7.  In
the appendix, section 8,  we show how the central limit procedure
leads to a concrete representation. In the forthcoming paper
\cite{JO}, we will use different but similar probabilistic tools
in showing that $OH$ embeds into the predual of a hyperfinite von
Neumann algebra.

We would like to thank U. Haagerup on collaboration on
Proposition \ref{Voi} and  R. Speicher for stressing  the
consequences of his central limit theorem in this context.

\section{Preliminaries}
We use standard notation in operator algebras as in
\cite{TAK,KR}. In particular, for $C^*$-algebras $A\subset B(H)$,
$B\subset B(K)$, we will denote by $A\ten_{min}B\subset B(H\ten_2
K)$ the induced minimal $C^*$-norm. Here (and in the following)
$H\ten_2 K\cong S_2(H,K)=HS(H,K)$ stands for the unique tensor
product making $H\ten K$ a Hilbert space (often called the
Hilbert-Schmidt norm). If $N,M$ are von Neumann algebras, we
denote by $N\bar{\ten}M$ the closure of $N\ten M$ in the weak
operator topology in $B(H\ten_2 K)$ (assuming the inclusions are
normal). For a $C^*$-algebra $A$, we denote by $A^{op}$ the
C$^*$-algebra defined on the same underlying Banach space but
with the reversed multiplication $x\circ y=yx$. By $\bar{A}$, we
denote the $C^*$-algebra obtained by changing the complex
multiplication $\la.x\lel \bar{\la}x$ on $A$. Thus $A$ and
$\bar{A}$ coincide as real Banach algebras. Then, we see that the
map $j:A^{op}\to \bar{A}$ given by $j(x)=x^*$ is
$C^*$-isomorphism.

We will use the notation $\ten_{\eps}$, $\ten_{\pi}$ for
smallest, biggest tensor norm on \underline{Banach} spaces,
respectively. Thus for a Hilbert space $H\ten_{\pi}H=S_1(H)$
corresponds to the trace class operators and
$S_\8(H)=\K(H)=H\ten_{\eps}H$ is used for the compact operators
on $H$. We note the trivial inclusions
 \begin{eqnarray}
  H\ten_{\pi} K &\subset& H\ten_2 K \pl \subset \pl
  H\ten_{\eps} K \pl .
  \end{eqnarray}
We also assume the reader to be familiar with standard operator
space terminology and refer to \cite{ER,P} for more details. We
will need some basic facts about the column Hilbert space
$H^c=B(\cz,H)$ and the row Hilbert space $H^r=B(H,\cz)$ of  a
given Hilbert space $H$. Given an element $x=[x_{ij}]\in
M_m(H^c)$, we observe that
 \begin{eqnarray}
  \noo x\rrm_{M_m(H^c)} &=& \noo x\rrm_{B(\ell_2^m,\ell_2^m(H))}
  \lel  \noo x^*x\rrm_{B(\ell_2^m)}^{\frac12} \lel
  \noo [\summ_k (x_{ki}, x_{kj})]_{ij} \rrm_{M_m}^{\frac12} \pl .
  \end{eqnarray}
Here $(x,y)$ denotes the scalar product of $x$ and $y$ in $H$. In
this paper, we will assume the scalar product to be antilinear in
the first component.  Similarly, we have
 \begin{eqnarray}
 \noo x\rrm_{M_m(H^r)} &=&
 \noo [\summ_k (x_{ik}, x_{jk})]_{ij} \rrm_{M_m}^{\frac12} \pl .
  \end{eqnarray}
We refer to \cite{Poh} for more details on the operator Hilbert
space $H^{oh}=[H^c,H^r]_{\frac12}$ and interpolation norms. In
particular, let us consider $H=\ell_2$ and denote by $(e_k)$ the
natural unit vector basis. Then  for all sequences $(x_k)\subset
B(H)$ with  associated linear map $u:OH\to B(H)$ defined by
$u(e_k)=x_k$ we have
 \begin{eqnarray} \label{ohnorm}
 \noo u:OH\to Im(u)\rrm_{cb} &=& \noo u\rrm_{cb}
 \lel \noo \summ_k x_k\ten
 \bar{x}_k\rrm_{B(H)\ten_{min} \bar{B}(H)}^{\frac12}
 \pl .
 \end{eqnarray}
$H^c$, $H^r$ and $H^{oh}$ are homogeneous hilbertian operator
space, i.e. for $s\in \{ c,r,oh\}$ and every bounded linear map
$u:H^s\to H^s$, we have
 \[ \noo u\rrm_{cb} \lel \noo u\rrm \pl .\]
In terms of general operator space notation, let us recall that a
{\emph complete contraction} $u:E\to F$ is given by a completely
bounded map with $\noo u\rrm_{cb}\le 1$. As usual if $G=E/F$ is a
quotient operator space (i.e. $M_m(G)=M_M(E)/M_m(F)$
isometrically) and $T:E\to X$ is a complete contraction which
vanishes on $F$, then $T$ induces a unique map $\hat{T}:E/F\to X$
which is still a complete contraction and defined by
$\hat{T}(x+F)=T(x)$. We refer to the introduction for the
definition  of the  projective and injective operator space tensor
products. The operator space projective tensor product is
projective, i.e. $E/F\wet X$ is a quotient of $E\wet X/F\wet X$.
However, if $N$ is an injective von Neumann algebra, we also have
an isometric inclusion
\begin{eqnarray}
 N_*\wet E_1 &\subset&  N_*\wet E_2
\end{eqnarray}
whenever $E_1\subset E_2$ (completely isometrically). (This can
easily be deduced from Wittstock's extension theorem and using
$(N_*\wet E)^*=CB(E,N)$.)

As for Banach spaces, we have  direct sums $E\oplus_p F$ of given
operator spaces $E$ and $F$. However, we should be careful in
defining the operator space structure. Indeed, for $p=\8$ and a
matrix $[x_{kl}]$ with $x_{kl}=(e_{kl},f_{kl})$, we have
 \[ \noo [x_{kl}]\rrm_{M_m[E\oplus_\8F]}
 \lel \max\{ \noo [e_{kl}]\rrm_{M_m(E)},  \noo
 [f_{kl}]\rrm_{M_m(E)}\} \pl .\]
The operator space $E\oplus_1F$ is defined by its canonical
inclusion in $(E^*\oplus_\8 F^*)^*$. For $1\le p\le \8$, we
define the operator space structure  by complex interpolation
 \[ E\oplus_p F \lel [E\oplus_\8 F,E\oplus_1
 F]_{\frac1p} \pl .\]
We will recall some basic facts about Haagerup's $L_p$ spaces and
refer to  \cite{Haa1,Haa2} and \cite{Te} for references on
operator valued weights and basic properties of the Haagerup $L_p$
spaces. Let $\phi$ be a normal, semifinite faithful weight on $N$
with modular automorphism group $\si^{\phi}_t$. The crossed
product $N\rtimes_{\si_t^{\phi}} \rz$ is the von Neumann algebra
generated in $B(L_2(\rz,H))$ by $\pi(N)$, $(\la(t))_{t\in \rz}$
where
 \[ \la(s)\xi(t)\lel \xi(t-s) \quad \mbox{and} \quad \pi(x)\xi(t) \lel \si_{-t}^{\phi}(x)\xi(t) \pl .\]
Then the modular automorphism group satisfies
 \[ \pi(\si_t(x)) \lel \la(t)\pi(x)\la(t)^* \pl \]
for all $t\in \rz$ and $x\in M$. The dual action
$\theta:N\rtimes_{\si_t^{\phi}}\rz\to N\rtimes_{\si_t^{\phi}}\rz$
is given by
 \[ \theta_s(x)\lel W(t)xW(t)^* \quad \mbox{where}\quad
 W(s)(\xi)(t)\lel e^{-ist}\xi(s) \pl \]
is the unitary implemented by the Fourier transformation.  Then
$N$ appears as the fixpoint algebra with respect to the action
$\theta$
 \[ \pi(N) \lel \{ x\in M\rtimes \rz \pl |\pl \theta_s(x)=x \mbox{ for all } s\in \rz\} \pl .\]
A crucial ingredient in the Haagerup $L_p$ space is the operator
valued weight $T:N\rtimes_{\si_t^{\phi}}\rz\to N$ given by
 \[ T(x)\lel \intt_{\rz} \theta_s(x) ds \pl .\]
It turns out that $\si_t^{\phi\circ T}$ is implemented by the
unitary  $\la(t)$ and therefore (see \cite{PT}) there exists a
normal semifinite faithful trace $\tau$ on
$N\rtimes_{\si_t^{\phi}}\rz$ and a positive selfadjoint operator
$h$ affiliated with $N\rtimes_{\si_t^{\phi}}\rz$ such that
\[ \phi\circ T(x)\lel \tau(hx) \]
for all $x\in N\rtimes_{\si_t^{\phi}}\rz$ and $h^{it}=\la(t)$. Let
us note that $\tau$ is the uniquely determined trace satisfying
$\theta_s\circ \tau=e^{-s}\tau$. Following Haagerup, we recall
the following definition
 \[ L_p(N)  \lel L_{p}(N,\phi) \lel \{ x \pl|\pl  x \mbox{ $\tau$-measurable and
 } \theta_s(x)\lel \exp(-\frac{s}{p})x \} \pl .\]
The (quasi-) norm on $L_{p}(N,\phi)$ is given by
 \[ \noo x\rrm_p \lel tr(|x|^p)^{\frac12} \pl .\]
It is shown in \cite{Te} that $L_{p}(N,\phi)$ and $L_{p}(N,\psi)$
are isomorphic whenever $\phi$ and $\psi$ are normal, faithful,
semifinite weights. Let us note that H\"{o}lder's inequality
 \[ x\in L_p(N) \quad ,\quad y\in L_q(N)  \quad \Rightarrow xy\in
 L_r(N) \quad \frac1r\lel \frac1p+\frac1q \]
holds in the context of Haagerup $L_p$ spaces:
 \[ \noo xy\rrm_r\kl \noo x\rrm_p \noo y\rrm_q \pl .\]
For $p=2$, we obtain a Hilbert space $L_2(N)$ with scalar product
$(x,y)=tr(x^*y)$. We will use the notation $L_2^{s}(N)$ instead
of $L_2(N)^s$ for $s\in \{c,r,oh\}$.

In the theory of operator spaces it is customary to use the
duality
 \[ \langle [x^*_{ij}],[x_{ij}]\rangle \lel
 \summ_{i,j=1}^n x^*_{ij}(x_{ij}) \]
between matrices $[x^*_{ij}]\in S_1^n\wet X^*$ and $[x_{ij}]\in
M_n(X)$. Unfortunately this is not consistent with the trace $tr$
on $n\times n$-matrices, which correspond to
 \[ <\langle  [x^*_{ij}],[x_{ij}]\rangle>  \lel
 \summ_{i,j=1}^n x^*_{ij}(x_{ji}) \pl . \]
This forces us to define the \emph{operator space structure on}
$L_1(N)$ by its action on $N^{op}$. Since $N$ and $N^{op}$
coincides as Banach spaces, we may consider $\iota:L_1(N)\to
(N^{op})^*$ defined by
 \[ \iota(d)(y) \lel tr(dy) \lel \phi_{d}(y) \pl .\]
Here $\phi_d$ is the linear functional associated to the density
$d$ in $L_1(N)$. Clearly, $\iota(L_1(N))\subset N^{op}_*$. If
$\phi$ is a normal semifinite weight then $\phi_n=tr\ten \phi$ is
a normal semifinite weight on $M_n(N)$. Moreover, $tr\ten \tau$
is the unique trace satisfying $tr \ten \tau\circ
\theta_s=e^{-s}tr\ten \tau$ and $tr_n\ten tr:L_1(M_n(N),tr\ten
\phi)\to \cz$ still yields the evaluation at $1$. Therefore, we
get
 \for
 \noo [\iota(x_{ij})] \rrm_{S_1^n\wet N^{op}_*}
 &=& \sup_{\noo [y_{ij}]\rrm_{M_n(N^{op})}\le 1}
 \bet \summ_{ij=1}^n \iota(x_{ij})(y_{ij})\rag \\
 &=& \sup_{\noo [y_{ij}]\rrm_{M_n(N)}\le 1}
 \bet \summ_{ij=1}^n tr(y_{ji}x_{ij}) \rag\\
 &=& \sup_{\noo [y_{ij}]\rrm_{M_n(N)}\le 1}
 \bet tr_n\ten tr([y_{ij}][x_{ij}]) \rag\\
 &=& \noo [x_{ij}] \rrm_{L_1(M_n\ten N,tr\ten \tau)}
 \pl .
 \mel
Therefore, the use of $N^{op}$ enables us to `untwist'  the
duality bracket and we have
 \begin{eqnarray}\label{os1}
  S_1^n \wet L_1(N,\phi) &=&  L_1(M_n\ten N,tr_n \ten \phi)
 \end{eqnarray}
Here we distinguish between the predual $N_*^{op}$ and the
`concrete' realization of this space as space of operators in
$L_1(N,\phi)$. If $N$ happens to be semifinite, we have a
canonical map $i_p:L_p(N,\tau)\to L_p^{Haagerup}(N,\tau)$ (see
\cite{Te} for the connection in this case)  such that
 \[ tr(i_p(x)i_{p'}(y))  \lel \tau(xy) \pl .\]
In particular, for $p=1$ the map $\iota i_1:L_1(N,\tau)\to
N^{op}_*$ is given by
 \[ \iota i_1(x)(y) \lel \tau(xy) \pl .\]
Thus
 \[ S_1^n \wet L_1(N,\tau) \lel   L_1(M_n\ten N,tr_n \ten \tau) \]
remains true if we read $L_1(N,\tau)$ and $L_1(M_n\ten N,tr_n
\ten \tau)$ according to their definition for semifinite von
Neumann algebras (see e.g. \cite{TAK}). This also applies for
$N=M_n$ and thus the mapping $I_n:S_1^n\to L_1(M_n,\tau_n)$ given
by $I_n(x)=nx$ and yields a complete isometry.

In harmonic analysis it is often more convenient to use the
anti-linear duality bracket
 \[ \lb y,x \rb  \lel \tau(y^*x) \pl .\]
This corresponds to an embedding $\bar{\iota}:L_1(N,\phi)\to
\bar{N}_*$. The map $\bar{\iota}:L_1(N,\tau)\to \bar{N}_*$ given
by
 \begin{eqnarray} \label{os2}
 \bar{\iota}(x)(y) &=& tr(y^*x)
 \end{eqnarray}
is a complete isometry onto $\bar{N}_*$. Indeed, this follows
immediately from the fact that $J:\bar{N}\to N^{op}$ defined by
$J(x)=x^*$ is an normal $^*$-isomorphism.

An important result of Effros and Ruan \cite{ER}  shows that the
projective tensor is compatible with von Neumann algebras.
Indeed, given von Neumann algebras $N$ and $M$ then
 \begin{eqnarray} \label{ER}
 (N_*\wet M_*)^* &=&  N\bar{\ten}M \pl .
 \end{eqnarray}
Using the opposite algebras this also provides a natural
isomorphism
 \begin{eqnarray}  \label{os3}
  L_1(M\ten N) &\cong &  L_1(M)\wet
  L_1(N)\pl .
 \end{eqnarray}
Note again, that this remains true for semifinite von Neumann
algebras. We will conclude these preliminaries by collecting
further facts about the operator space projective tensor product
in connection with $L_1$ spaces. First, we note that in the
category of $L_1$ spaces the row Hilbert space can be represented
by columns and vice versa: Namely, let $(e_j)$ be the standard
vector basis on $\ell_2$, then for a finite sequence $(x_j)\in
L_1(N)$
\begin{eqnarray} \label{cn}
 \noo \summ_{j} x_j\ten e_j \rrm_{L_1(N)\wet \ell_2^r}
  &=&   \noo (\summ_j x_j^*x_j)^{\frac12} \rrm_{L_1(N)} \pl .
\end{eqnarray}
Indeed, by duality we have
 \[ (L_1(N)\wet \ell_2^r)^*\lel CB(\ell_2^r,L_1(N)^*)
 \lel CB(\ell_2^r,N^{op}) \pl .\]
Thus, we deduce from (\ref{os2}) that
  \for
  \noo   \summ_{j} x_j\ten e_j \rrm_{L_1(N)\wet \ell_2^r}
  &=& \sup_{\noo u:\ell_2^r\to N^{op}\rrm_{cb} \le 1}  \bet
  \summ_j \langle \iota(x_j),u(e_j)\rangle
  \rag \\
  &=&
  \sup_{ \noo \summ_j y_j^*y_j\rrm_{N^{op}} \le 1} \bet \summ_j
  tr(y_jx_j) \rag \lel
    \sup_{ \noo \summ_j y_jy_j^*\rrm_N \le 1} \bet \summ_j
  tr(y_jx_j) \rag \\
  &=&
    \noo (\summ_j x_j^*x_j)^{\frac12} \rrm_{L_1(N)} \pl .
  \mel

Note that on the level of $L_1(N)$-spaces the tensor product with
respect to $\ell_2^r$, $\ell_2^c$ corresponds to the columns, rows
in $L_1(N\ten B(\ell_2))$, $L_1(N\ten B(\ell_2))$, respectively.
This `switch' is a well-known phenomenon for operator spaces on
the level of matrices. As an application, we obtain formulas for
norms in $H^c$ and $H^r$. Let $x\in L_1(N)\wet H^r$, $y \in
L_1(N) \wet H^c$, then
\begin{eqnarray} \label{hc}
 \noo x\rrm_{L_1(N)\wet H^r}&=&  \noo (x^*,x)  \rrm_{\frac12}^{\frac12}
  \quad \mbox{and} \quad \noo x\rrm_{L_1(N)\wet H^c} \lel  \noo (x,x^*)
  \rrm_{\frac12}^{\frac12}
\end{eqnarray}
Indeed, let $(e_i)$ be an orthonormal basis of $H$. By density,
we may assume $x\lel \sum_i x_i\ten e_i$. Then, we apply
\eqref{cn} and deduce
 \for
 \noo x\rrm_{L_1(N)\wet H^r} &= \noo \kla \summ_i x_i^*x_i
 \mer^{\frac12}\rrm_{L_1(N)} \lel \noo \summ_i x_i^*x_i
 \rrm_{L_{\frac12}(N)}^{\frac12} \pl .
 \mel
On the other hand, we may define the bilinear map
 \[ ( \pll ,\hspace{0.2cm} ):\overline{L_1(N)\wet H^r} \times L_1(N)\wet H^r \to L_{\frac12}(N) \]
on elementary tensors as $(a\ten h,b\ten k)\lel (h,k) a^*b$. For
finite sums $x=\sum_i x_i\ten e_i$, we have
 \begin{equation}\label{cs0}
  (x,y)\lel \summ_{i} x_i^*y_i \quad \mbox{and} \quad
  \noo \summ_i x_i^*y_i\rrm_{\frac12}\kl
  \noo \summ_i
 x_i^*x_i\rrm_{\frac12}^{\frac12} \noo \summ_i
  y^*_iy_i\rrm_{\frac12}^{\frac12} \pl .\end{equation}
Hence continuity implies \eqref{hc}. In the sequel, we will also
the cb-norm of the the multiplication maps $R_b:L_2(N)\to L_1(N)$
and $L_b:L_2(N)\to L_1(N)$  for $b\in L_2(N)$:
 \begin{eqnarray}
  \noo R_b:L_2^r(N)\to L_1(N)\rrm_{cb} &=& \noo b\rrm_2
  \lel \noo L_b:L_2^c(N)\to L_1(N)\rrm_{cb}  \pl .
 \end{eqnarray}
Indeed, we first note that $(S_1^m\wet H^r)^*=M_m(H^c)$ and
\eqref{cs0}  implies
 \begin{equation} \label{cs1} |\summ_{ij} (y_{ij},x_{ij})|
 \kl \noo [y_{ij}] \rrm_{M_m(H^c)} \noo [x_{ij}]\rrm_{S_1^m\wet
 H^r}  \pl .\end{equation}
Let $[x_{ij}]\subset  L_2^r(N)$ be a matrix. Then, we deduce from
$(1.7)$ and \eqref{cs1} that
 \for
 \noo [x_{ij}b]\rrm_{S_1^m\wet L_1(N)} &=& \sup_{\noo [y_{ij}]\rrm_{M_n(N)}\le 1}
 \bet \summ_{ij} tr(y_{ij}^*x_{ij}b) \rag
 \lel \sup_{\noo [y_{ij}]\rrm_{M_n(N)}\le 1}
 \bet \summ_{ij} (y_{ij}b^*,x_{ij}) \rag \\
  &\le&  \noo [y_{ij}b^*] \rrm_{M_m(L_2^c(N))}  \pl  \noo [x_{ij}]
  \rrm_{S_1^m\wet L_2^r(N)} \pl .
 \mel
However, let $\phi$ be the positive functional
$\phi(y)=tr(yb^*b)$. Then, we deduce from the fact that
$\phi:N\to \cz$ is completely  bounded that
 \begin{align}
  \noo [y_{ij}b^*]\rrm_{M_m(L_2^c(N))}^2
  &= \noo \big[ \summ_{k} tr(by_{ki}^*y_{kj}b^*) \big]\rrm_{M_m}
  \lel
   \noo \big[ \summ_{k} \phi(y_{ki}^*y_{kj}) \big]  \rrm_{M_m} \notag \\
  &\le \phi(1)  \noo y^*y\rrm_{M_m(N)} \kl
  \noo b\rrm_2^2   \noo y^*y\rrm_{M_m(N)} \label{comul}
      \pl .
  \end{align}
The estimate from below follows from $\noo R_b(b^*)\rrm_1= \noo
b^*b\rrm_1=\noo b\rrm_2^2$. The proof for the left multiplication
is identical.

%\begin{cor}\label{spur} Let $H$ be Hilbert space and $x\in M_m(H^c)$, $y\in S_1^m\wet H^r$, then
% \[ \bet \summ_{ij} \tau(y_{ij},x_{ij})\rag \kl \noo x \rrm_{M_m(H^c)} \noo y\rrm_{S_1^m\wet H^r} \pl .\]
%\end{cor}\lz

\section{A logarithmic characterization of C$^*$-algebras with
WEP}

We will show that the `logarithmic little Grothendieck'
inequality $(0.4)$ (or even $(0.2)$) only holds for C$^*$-algebras
with WEP. We will assume $(0.4)$ for the C$^*$-algebra $B(H)$, see
\cite{PS} (or section 4 for an independent proof).

\begin{lemma}\label{triv}   Let $A$ be a $C^*$-algebra with WEP and $u:A\to OH$.
Then
 \for
  \kla \summ_{k=1}^n \noo u(x_k)\rrm^2 \mer^{\frac12}
  &\le& C (1+ \ln n)^{\frac12} \pl \noo u\rrm_{cb} \noo \summ_{k=1}^n x_k\ten
  \bar{x}_k\rrm^{\frac12} \pl
  \mel
holds for all $x_1,...,x_n\in A$.
\end{lemma}\lz

\begin{proof}[\bf Proof:] Let $A\subset A^{**}\subset B(H)$ and $P:B(H)\to
A^{**}$ be a contraction such that $P|_A=id_A$. Let $j:A^*\to
A^{**}$ be the canonical embedding, then the contraction
$E=j^*P^{**}:B(H)^{**}\to A^{**}$ satisfies
$E_{A^{**}}=id_{A^{**}}$. Hence $E$ is a conditional expectation
by Sakai's theorem, see \cite{TAK}, and in particular completely
contractive. The restriction of $E$ to $B(H)$ coincides with $P$
and in particular $P$ is a complete contraction. Then
$v=u^{**}P:B(H)\to OH_n$ satisfies
 \[ \noo v\rrm_{cb}\kl \noo P\rrm_{cb} \pl \noo u^{**}\rrm_{cb}
 \kl \noo u\rrm_{cb} \pl .\]
We can apply $(0.4)$ to $v$ and deduce the assertion from the
injectivity of the minimal tensor norm, see $(1.4)$.\qd

\begin{lemma}\label{cbohne} Let $N$ be a von Neumann algebra with a faithful, normal semifinite weight $\phi$.
Let and $a,b\in L_4(N)$, then the map  $M_{ab}:N\to L_2^{oh}(N)$
defined by
 \[ M_{ab}(x)\lel axb \]
satisfies
 \[ \noo M_{ab}:N\to L_2^{oh}(N)\rrm_{cb} \kl \noo a\rrm_4 \pl
 \noo b\rrm_4 \pl .\]
\end{lemma}\lz

\begin{proof}[\bf Proof:] By homgeneity we can assume $\noo a\rrm_4=\noo
b\rrm_4=1$. Using the left and right action of $N$ and polar
decomposition, we can assume that $a$, $b$ are positive. Let $e$
and $f$ be the support projections of $a$, $b$ respectively. We
consider $M_2(N)$ and the canonical matrix units
$e_{11},e_{12},e_{21}$ and $e_{22}$. Then the support $s(\phi)$
and the density $D$ of the state
$\phi(x)=tr(a^4e_{11}xe_{11}+b^4e_{22}xe_{22})$, $x\in M_2(N)$,
is given by
 \[ s(\phi) \lel \kla \begin{array}{cc} e &0\\0&f\end{array}\mer \quad,\quad
 D \lel \kla \begin{array}{cc} a^4&0\\0&b^4\end{array}\mer \pl .\]
By the invariance of $D^{4it}M_2(N)D^{-4it}\subset
s(\phi)M_2(N)s(\phi)$, we deduce that $a^{it}xb^{-it}\in N$ for
every $x\in eNf$. Let $x\in M_m(eNf)$. We have to show that
 \[ \noo id_{M_m}\ten v(x)\rrm_{M_m(L_2^{oh}(N))}
 \kl \noo x\rrm_{M_m(N)} \pl .\]
We define the complex function $f:\{0\le Re z \le 1\}\to M_m\ten
L_2(eNf)$ by
 \[ f(z) \lel (1\ten a^{2z})x(1\ten b^{2(1-z)}) \pl .\]
Observe that for $z=a+ib$ and $x_{kl}\in eNf$
  \[  a^{2z}x_{kl}b^{2(1-z)} \lel a^{2a}
  a^{2ib}x_{kl}b^{-2ib}b^{2-2a} \in L_2(N) \pl. \]
It follows from \cite{Te} that $z\mapsto a^{2z}x_{kl}b^{2(1-z)}$
is complex differentiable in the measure topology of
$\tilde{N}=N\rtimes_{\si_t}\rz$. Moreover, using spectral
projections $e_n=1_{[\frac1n,n]}(a)$ and
$f_n=1_{[\frac1n,n]}(b)$, we see that $e_n(a\log a)$ and $(b\log
b)f_n$ are bounded and hence $e_na^{2z}x_{kl}b^{2(1-z)}f_n$ is
even differentiable in $L_{2,\8}(\tilde{N})$ in the interior of
the strip. By approximation $z\mapsto a^{2z}x_{kl}b^{2(1-z)}$ is
analytic with values in $L_{2,\8}(\tilde{N})$.
 Since, the topology of
$L_2(N)$ is induced by its embedding into $L_{2,\8}(\tilde{N})$,
we see that $f$ is analytic in the interior of the strip. We note
that the functional $\phi_b:N\to \cz$, $\phi_b(x)=tr(xb^4)$
induces a contractive functional, hence a scalar valued
completely contractive map. Since $x\mapsto a^{2it}xb^{-2it}$ is
a contraction on $L_2(eNf)$ and thus by homogeneity this map is
also a complete contraction on $L_2^c(eNf)$. We deduce from
\eqref{comul} \for
 \noo f(it)\rrm_{M_m(L_2^c(N))} &=&
 \noo [a^{2it}x_{kl}b^2b^{-2it}]\rrm_{M_m(L_2^c(N))}
  \kl    \noo [x_{kl}b^2]\rrm_{M_m(L_2^c(N))} \\
  &\le&  tr(b^4)   \noo [x_{kl}]\rrm_{M_m(N)} \kl  \noo
  x\rrm_{M_m(N)}\pl .
 \mel
Similarly, we have
 \[ \noo f(1+it)\rrm_{M_m(L_2^r(N))} \kl \noo a^4\rrm_4 \noo
 x\rrm_{M_m(N)} \kl \noo  x\rrm_{M_m(N)} \pl . \]
The three line lemma implies
 \begin{align*}
 \noo (1\ten a)x(1\ten b)\rrm_{M_m(L_2^{oh}(N))}&\le  \noo x\rrm_{M_m(N)} \pl . \qedhere
 \end{align*}
\qd

\begin{proof}[\bf Proof of Theorem 8:] According to Lemma \ref{triv} every
$C^*$-algebra with WEP satisfies the `little Grothen- dieck
inequality' with factor $c\sqrt{1+\ln n}$. Conversely, we assume
that $A$ is a C$^*$-algebra such that for some $C>0$ and $\al>0$
 \begin{eqnarray}
  \kla \summ_{k=1}^n \noo u(x_k)\rrm^2 \mer^{\frac12}
  &\le&  C (1+\ln n)^{\alpha} \pl \noo u\rrm_{cb} \noo \summ_{k=1}^n x_k \ten
  \bar{x}_k\rrm_{A\ten_{min}  \bar{A}}
  \end{eqnarray}
holds for  all $\nen$  and $x_1,...,x_n\in A$ and all $u:A\to OH$.
Let us consider the von Neumann algebra $N=A^{**}$. Recall that
$N$ is in standard form on $L_2(N)$, in particular $N$ acts on
$L_2(N)$ by left multiplication $\pi(x)h=xh$ and $J(h)=h^*$ is an
anti-linear isometry $J$ such that $N'\lel JNJ$, see \cite{Te}.
Let $x_1,...,x_n\in A$ and $h\in L_2(N)$ be a unit vector. Then,
we can find $a,b\in L_4(N)$ of norm $1$ such that $h=ab$.
According to Lemma \ref{triv},  the maps $M_{a^*a}:N\to
L_2^{oh}(N)$ and $M_{bb^*}:N\to L_2^{oh}(N)$ are complete
contractions and therefore we deduce from $(2.1)$ that
 \begin{eqnarray}
 \begin{minipage}{14cm}\vspace{-0.3cm}
 \for
 (h,\summ_{j=1}^n x_iJx_iJh) &=&  (h,\summ_{j=1}^n x_ihx_i^*) \lel
 \summ_{i=1}^n tr(h^*x_ihx_i^*) \nonumber \\
 &=& \summ_{i=1}^n tr(a^*x_iabx_i^*b^*) \nonumber\\
 %  \end{eqnarray} \begin{eqnarray}
  &\le& \kla \summ_{i=1}^n \noo a^*x_ia\rrm_2^2 \mer^{\frac12}
 \kla \summ_{i=1}^n \noo bx_ib^*\rrm_2^2 \mer^{\frac12}
 \nonumber\\
  &\le& c^2(1+\log n)^{2\al} \noo \summ_{k=1}^n x_k\ten
 \bar{x}_k\rrm_{A\ten_{min}\bar{A}} \pl .
 \mel
 \end{minipage}
 \end{eqnarray}
In order to eliminate the log-term, we use Haagerup's trick and
consider the positive  operator
  \for
 \lefteqn{ [(\summ_{i=1}^n x_iJx_iJ)^*(\summ_{i=1}^n
 x_iJx_iJ)]^m}\\
   & & \lel  \summ_{i_1,...,i_{2m}=1}^n   x_{i_1}^*(Jx_{i_1}^*J)x_{i_2}(Jx_{i_2}J) \cdots
     \cdots     x_{i_{2m-1}}^*(Jx_{i_{2m-1}}^*J)      x_{i_{2m}}(Jx_{i_{2m-1}}J) \\
   & & \lel  \summ_{i_1,...,i_{2m}=1}^n   x_{i_1}^*x_{i_2}\cdots
   x_{i_{2m-1}}^*x_{i_{2m}}  (Jx_{i_1}^*x_{i_2}\cdots
   x_{i_{2m-1}}^*x_{i_{2m}}
   J)  \\
      \mel
We apply $(2.2)$ for the finite family
$x_{i_1,....,i_{2m}}=x_{i_1}^*x_{i_2}\cdots
   x_{i_{2m-1}}^*x_{i_{2m}}$  and deduce
\for
 \noo (\summ_{i=1}^n x_iJx_iJ)]^m(h)\rrm^2
 &=& \summ_{i_1,...,i_{2m}=1}^n (h,x_{i_1}^*x_{i_2}\cdots
   x_{i_{2m-1}}^*x_{i_{2m}}Jx_{i_1}^*x_{i_2}\cdots
   x_{i_{2m-1}}^*x_{i_{2m}}
   Jh) \\
 &\le& c^2 (1+\ln n^{2m})^{\al} \pl \noo \summ_{i_1,....,i_{2m}=1}^n
 x_{i_1,....,i_{2m}} \ten \bar{x}_{i_1,....,i_{2m}}
 \rrm_{A\ten_{min} \bar{A}} \\
 &=&  c^2 (1+\ln n^{2m})^{\al} \pl \noo [(\summ_{i=1}^n x_i\ten
 \bar{x}_i)^* (\summ_{i=1}^n x_i\ten \bar{x}_i)]^m \rrm_{A\ten_{min} \bar{A}} \\
 &=& c^2 (1+2m \ln n)^{\al} \noo \summ_{i=1}^n x_i\ten
 \bar{x}_i\rrm^{2m}_{A\ten_{min} \bar{A}} \pl .
 \mel
Taking the supremum over $\noo h\rrm\le 1$, we get
 \[ \noo \summ_{i=1}^n x_iJx_iJ \rrm \kl c^{\frac1m}
 (1+2m)^{\frac{\al}{2m}} (1+ \ln n)^{\frac{\al}{m}}
 \pl \noo \summ_{i=1}^n x_i\ten
 \bar{x}_i\rrm_{A\ten_{min} \bar{A}} \pl .\]
Taking the limit for $m\to \8$, we obtain  (in the language of
\cite{P3}) that
 \for
 \noo \summ_{i=1}^n L_{x_i}R_{x_i^*}\rrm &=&
  \noo \summ_{i=1}^n x_iJx_iJ \rrm \kl \noo \summ_{i=1}^n x_i\ten
 \bar{x}_i\rrm_{A\ten_{min}\bar{A}} \\
 &=&
  \noo \summ_{i=1}^n x_i\ten
  \bar{x}_i\rrm_{B(H)\ten_{min}\bar{B}(H) }
  \kl  \noo \summ_{i=1}^n x_i\ten
  \bar{x}_i\rrm_{B(H)\ten_{max}\bar{B}(H) } \pl .
 \mel
(Note that the last inequality is indeed an equality). Let us
recall an equality proved by  Pisier \cite[Theorem 2.1]{P3} (and
\cite[Theorem 2.9]{Haa} in the non-semifinite case)
 \begin{eqnarray}
 \noo \summ_{i=1}^n x_i \ten \bar{x}_i\rrm_{A\ten_{max} \bar{A}}
 &=&
 \noo  (x_1,...,x_n)\rrm_{[R_n(A),C_n(A)]_{\frac12}}^2 \lel
   \noo \summ_{i=1}^n L_{x_i}R_{x_i^*}\rrm \pl .
 \end{eqnarray}
This implies
 \for
 \noo \summ_{i=1}^n x_i \ten \bar{x}_i\rrm_{A\ten_{max} \bar{A}}
 &=& \noo \summ_{i=1}^n L_{x_i}R_{x_i^*}\rrm
 \kl  \noo \summ_{i=1}^n x_i\ten
  \bar{x}_i\rrm_{B(H)\ten_{max}\bar{B}(H) }  \pl . \nonumber
  \mel
According to \cite[Theorem 3.7]{Haa} $A$ has WEP.\qd

\begin{rem} {\rm In case of a von Neumann algebra $A=N$, we deduce
from Pisier's characterization \cite[Theorem 2.9]{Poh} that $N$
is injective if and only if
  \[ \noo (x_1,...,x_n)\rrm_{[C_n(N),R_n(N)]_{\frac12}} \kl \noo
   \summ_{i=1}^n x_i \ten \bar{x}_i
   \rrm_{N\ten_{min}\bar{N}}^{\frac12}\pl .
   \]
holds for all $x_1,...,x_n\in N$. In \cite{Poh} Pisier proved
this for semifinite von Neumann algebras, for  the general case
we refer to \cite{P3} that this characterization remains true for
arbitrary von Neumann algebras. Thus modulo the results in
\cite{P3}, the argument above provides a complete proof for the
fact that the logarithmic `little Grothendieck' characterizes
injective von Neumann algebras without appealing to the deep
results of Haagerup in \cite{Haa}.}
\end{rem}

\section{Pusz/Woronowicz' formula and the operator space $OH$}

We will recall the Pusz/Woronowicz formula for square roots of
sesquilinear forms and derive a dual version. This is the key
element in identifying $OH$ as a subspace of a quotient of
$L_2^c\oplus L_2^r$. Given two sesquilinear forms
 \[ \al(x,x)\lel (Ax,x) \quad \mbox{and} \quad \beta(x,x)\lel (Bx,x) \]
for positive {\it commuting} operators $A,B$ on a Hilbert space
$H$, then according to \cite{Wo}
 \begin{eqnarray}
    ((AB)^{\frac12}x,x) &=&
   \inf_{x=a(t)+b(t)} \intt_0^1 \frac{\al(a(t),a(t))}{t}+
  \frac{\al(b(t),b(t))}{1-t} \pl \frac{dt}{\pi \sqrt{t(1-t)}}
  \pl .
 \end{eqnarray}
Here the infimum is taken over piecewise constant functions in
$H$, see \cite[Appendix]{Wo} for a proof. Pusz/Woronowicz use the
symbol
\[ \sqrt{\al\beta}(x,x) \lel ((AB)^{\frac12}x,x) \pl \]
and show (in a sense which is not obvious from our presentation)
that the symbol $\sqrt{\al\beta}$ is independent of the
`representation' of $\al$, $\beta$ by $A$, $B$, respectively. Let
us denote by $H_{\sqrt{\al \beta}}$ the Hilbert space $H$
equipped with the scalar product $\sqrt{\al \beta}$. In the
following, we will use the symbol $\mu$ for the probability
measure
 \[ d\mu(t)\lel \frac{dt}{\pi  \sqrt{t(1-t)}} \pl . \]
Let us rewrite equality (3.1)  in terms of subspaces and
quotients of Hilbert spaces. We define
 \begin{equation} \label{nu}
 d\nu_1(t)\lel t^{-1}d\mu(t)\quad \mbox{and} \quad
 d\nu_2(t) \lel (1-t)^{-1} d\mu(t)\end{equation}
and the Hilbert space $H_1=L_2(\nu_1,H_\al)$, where $\noo
x\rrm_{H_\al}=\al(x,x)^{\frac12}$. Similarly,
$H_2=L_2(\nu_2,H_\beta)$ with $\noo
x\rrm_\beta=\beta(x,x)^{\frac12}$.  If $A$ and $B$ are
invertible, the canonical inclusion map $H_{\al}\subset H$ and
$H_{\beta}\subset H$ are continuous and we may  define the linear
 map $Q:H_1\oplus_2 H_2\to L_0(\mu;H)$ by $Q(f,g)(t)\lel
f(t)+g(t)\in H$. We denote by
 \[ K \lel H_1\oplus_2 H_2/kerQ \]
the quotient space and
 \[ E\lel \{(f,g)+kerQ \pl|\pl Q(f,g) \mbox{ constant a.e.}
 \} \]
In this terminology the Pusz/Woronowicz formula reads as follows.
\lz

\begin{lemma} \label{W1} If $A$ and $B$ are invertible, then  $H_{\sqrt{\al \beta}}$ is isometricaly isomorphic to the subspace
$E\subset K$.
\end{lemma}\lz

Using the Hahn-Banach theorem, we will derive a second formula
for the square root $\sqrt{\al \beta}$. We will first
characterize the linear functionals on $E$.\lz

\begin{lemma}\label{lf} Let $A$ and $B$ be boounded
and  invertible.  Linear  functionals  $\phi:E\to \cz$ are given
by a pair $(f,g)\in H_1\oplus H_2$ such that
 \begin{equation} \label{caracf}  \frac{Af(t)}{t} \lel \frac{Bg(t)}{1-t} \quad \mu \mbox{a.e.} \quad
 \mbox{and} \quad
  \phi(x) \lel \intt_0^1 (A\frac{f(t)}{t},x) \pl
 d\mu(t) \pl .\end{equation}
Moreover,
 \[ \noo \phi\rrm \lel \inf (\noo f\rrm_{H_1}^2 +\noo
 g\rrm_{H_2}^2)^{\frac12} \pl ,\]
where the infimum is taken over all pairs satisfying
\eqref{caracf}.
\end{lemma}\lz

\begin{proof}[\bf Proof:] By the Hahn-Banach's theorem, the norm
one functionals on $E$ are in one to one correspondence with the
restrictions of $\hat{\phi}|_E$ of norm one functionals
$\hat{\phi}:K\to \cz$. Moreover,
 \[ \noo \hat{\phi}\rrm \lel \noo  \hat{\phi}Q:H_1\oplus_2H_2\to \cz\rrm \pl .\]
Therefore, norm a norm one functional $\phi:E\to \cz$ is given by
by a pair $(f,g) \in H_1\oplus_2 H_2$ such that
 \begin{eqnarray}
  \intt_0^1 \al(f(t),f(t)) \frac{d\mu(t)}{t} +
  \intt_0^1 \beta(g(t),g(t)) \frac{d\mu(t)}{1-t}
  &=&  \noo\hat{\phi} \rrm^2 \lel
  \noo \phi \rrm^2 \lel 1 \pl
 \end{eqnarray}
and
 \[ \hat{\phi}Q(h_1,h_2)  \lel  \intt_0^1 \al(f(t),h_1(t)) \frac{d\mu(t)}{t} +
  \intt_0^1 \beta(f(t),h_2(t)) \frac{d\mu(t)}{1-t} \pl. \]
Since $Q$ vanishes on all pairs $(h,-h)$, we deduce
 \for
 0&=&  \intt_0^1 \al(f(t),h(t)) \frac{d\mu(t)}{t} +
  \intt_0^1 \beta(g(t),-h(t)) \frac{d\mu(t)}{1-t} \\
  &=&  \intt_0^1 (Af(t),h(t)) \frac{d\mu(t)}{t} +
  \intt_0^1 (Bg(t),-h(t)) \frac{d\mu(t)}{1-t} \\
  &=&
   \intt_0^1 \big(\frac{Af(t)}{t}-\frac{Bg(t)}{1-t},h(t)\big) \pl d\mu(t)  \pl .
 \mel
Since $H_\al=H_{\beta}=H$ as Banach spaces,  this implies
\begin{eqnarray}
  \frac{Af(t)}{t}&=&    \frac{Bg(t)}{1-t} \quad
  \mu  \pl \mbox{a.e.} \label{ABea}
 \end{eqnarray}
We will now calculate how $\phi$ acts on $E$. Indeed, let $x\in
H$ and note that $x\lel \sqrt{t}x+(1-\sqrt{t})x$. It is easily
checked that $a(t)=\sqrt{t}$, $b(t)=1-\sqrt{t}$ are elements in
$L_2(\nu_1)$, $L_2(\nu_2)$, respectively. Thus, we deduce from
\eqref{ABea} that
 \begin{align}
  \phi(x) &= \hat{\phi}Q(x\sqrt{t},x(1-\sqrt{t})) \lel
   \intt_0^1 \big(A\frac{f(t)}{t},\sqrt{t}x\big) d\mu(t) +
   \intt_0^1 \big(B\frac{g(t)}{1-t},(1-\sqrt{t})x\big) d\mu(t) \notag \\
  &= \intt_0^1 \big(A \frac{f(t)}{t},x\big) d\mu(t) \pl .
  \label{ff}
  \end{align}
Conversely, any pair $(f,g)$ satisfying \eqref{ABea} induces a
state $\hat{\phi}$ which vanishes on $kerQ$. The restriction of
$\hat{\phi}$ to $E$ is given by \eqref{ff}. \qd

\begin{lemma} \label{W2} Let $A,A^{-1}$, $B,B^{-1}$ be bounded. Let $y\in H$, then
 \[ ((AB)^{\frac12}y,y) \lel \inf \intt_0^1 \frac{(Af(t),f(t))}{t} d\mu(t)  +
  \intt_0^1 \frac{(Bg(t),g(t))}{1-t} d\mu(t) \pl ,\]
where the infimum is taken over all tuples $(f,g)$ of $H$-valued
measurable functions satisfying
$\frac{Af(t)}{t}=\frac{Bg(t)}{1-t}$ $\mu$ a.e. and
 \[ B^{\frac12}y \lel \intt_0^1 \frac{A^{\frac12}f(t)}{t} \pl d\mu(t) \pl .\]
\end{lemma}\lz

\begin{proof}[\bf Proof:] Since $\sqrt{\al \beta}$ is a positive
sesquilinear form, we have
 \[ \noo y \rrm_{H_{\sqrt{\al \beta}}}  \lel ((AB)^{\frac12}y,y)^{\frac12}
 \lel \sup_{\sqrt{\al \beta}(x,x)\le 1} |((AB)^{\frac12}y,x)|    \pl .\]
Thus the norm of $y$ in $H_{\sqrt{\al\beta}}$ coincides with the
norm of the linear functional
 \[ \phi_{y}(x)\lel ((AB)^{\frac12}y,x) \pl .\]
According to Lemma \ref{lf}, we  can find $(f,g)$ such that
 \begin{eqnarray}
   ((AB)^{\frac12}y,x) &=&   \big(\intt_0^1  \frac{Af(t)}{t} d\mu(t),x
   \big)\pl,
 \end{eqnarray}
and
 \for
   ((AB)^{\frac12} y,y) &=&  \noo \phi_{y} \rrm^2
    \lel
  \intt_0^1  (Af(t),f(t)) \frac{d\mu(t)}{t} +
   \intt_0^1  (Bg(t),g(t)) \frac{d\mu(t)}{1-t} \pl .
   \mel
and $A \frac{f(t)}{t}=B \frac{g(t)}{1-t}$ $\mu$ a.e. According to
Using $(3.7)$ we get
 \[ B^{\frac12}y \lel \intt_0^1 \frac{A^{\frac12}f(t)}{t} \pl d\mu(t)   \pl .\]
Conversely, any pair satisfying these conditions induces the
functional $\phi_{y}$ and thus provides an upper estimate for the
norm of $y$. The assertion is proved.\qd

Motivated by these formulae, we introduce the following operator
space $F_n$ ($n\in \nz\cup \{\infty\}$) as a subspace of the
quotient of $L_2^c(\nu_1,\ell_2^n)\oplus_1
L_2^r(\nu_2,\ell_2^n)$. Indeed, let
 \[ G_n \lel  L_2^c(\nu_1,\ell_2^n) \oplus_1  L_2^r(\nu_2,\ell_2^n)/ker Q \pl ,\]
where $Q:L_2(\nu_1,\ell_2^n)\oplus_1 L_2(\nu_2,\ell_2^n)\to
L^0(\mu;\ell_2^n)$ is given by
 \[ Q(f,g)(t) \lel f(t)+g(t)\in \ell_2^n  \pl .\]
Then $F_n$ consists of the equivalence classes $(f,g)+kerQ\in
G_n$ such that $Q(f,g)$ is $\mu$-almost everywhere a constant
element in $\ell_2^n$. Let us denote by $f_1,..,f_n$ the
canonical unit vector basis in $F_n$ given by
$f_k=(\sqrt{t}e_k,(1-\sqrt{t})e_k)+kerQ$. Here $e_1,..,e_n$ is the
standard unit vector basis in $\ell_2^n$. For $n=\infty$, we read
$\ell_2^{\infty}=\ell_2$ and denote the outcome by $G$ and $F$,
respectively. The following lemma is proved using the polar
decomposition and the density of invertible matrices.\lz

\begin{lemma} \label{t1}  Let $x_1,...,x_n\in M_m$, then
\[ \noo \summ_{k=1}^n x_k \ten \bar{x}_k \rrm_{M_m\ten_{min} M_m}^{\frac12}
 \kl  \sup_{\noo a\rrm_4\le 1, \noo b\rrm_4\le 1, a> 0, b>0} \kla \summ_{k=1}^n \noo bx_ka\rrm_2^2 \mer^{\frac12} \pl. \]
Here $a>0$ means that $a\ge 0$ and $a$ is invertible.
\end{lemma} \lz

\begin{lemma} \label{1s} The natural inclusion map $id:F_n\to OH_n$ has $cb$-norm less $\sqrt{2}$.
\end{lemma} \lz

\begin{proof}[\bf Proof:] Let $x_1,...,x_n\in M_m$ and assume
$ \noo \sum_{k=1}^n x_k\ten f_k\rrm_{M_m(F_n)}\pl <\pl 1$. Then,
there are elements $f\in M_m(L_2^c(\nu_1;\ell_2^n))$ and $g\in
M_m(L_2^r(\nu_1;\ell_2^n))$ such that
 \[ x_k \lel f_k(t) + g_k(t) \quad \mbox{$\mu$ a.e. }\]
and
\[
   \max\left\{\noo \intt_0^1 \summ_{k=1}^n f_k(t)^*f_k(t) \frac{d\mu(t)}{t}
 \rrm_{M_m}^{\frac12},
   \noo \intt_0^1 \summ_{k=1}^n g_k(t)g_k(t)^* \frac{d\mu(t)}{1-t} \rrm_{M_m}^{\frac12}\right\}  \kl 1 \pl .\]
Let $a,b$ be  positive, invertible norm one elements in $S_4^m$.
On the Hilbert space $H=\ell_2^n(S_2^m)$ with scalar product
 \[ ((x_k),(y_k)) \lel \summ_{k=1}^n tr(x_k^*y_k) \]
we  define $A(x_k)=(x_ka^4)$ and $B(x_k)=(b^4x_k)$. Clearly,
these operators commute and we deduce from $(3.1)$ that
 \for
 \summ_{k=1}^n \noo b x_ka\rrm_2^2
 &=& \summ_{k=1}^n tr(a^*x_k^*b^*bx_ka)
 \lel  \summ_{k=1}^n tr(a^2x_k^*b^2x_k) \lel
  ( (AB)^{\frac12}(x_k),(x_k)) \mel
  \for
 &\le& \intt_0^1 (A(f_k(t)),(f_k)) \frac{d\mu(t)}{t} +
 \intt_0^1 (B(g_k(t)),(g_k)) \frac{d\mu(t)}{1-t} \\
 &=&  \intt_0^1 \summ_{k=1}^n tr(a^4f_k^*(t)f_k(t)) \frac{d\mu(t)}{t} +
 \intt_0^1 \summ_{k=1}^n tr(g_k(t)^*b^4g_k(t)) \frac{d\mu(t)}{1-t} \\
 &=& tr\Bigg(a^4 \intt_0^1 \summ_{k=1}^n f_k^*(t)f_k(t) \frac{d\mu(t)}{t}\Bigg) +
 tr\Bigg(b^4\intt_0^1 \summ_{k=1}^n g_k(t)g_k(t)^*
 \frac{d\mu(t)}{1-t}\Bigg) \\
 % \mel \for
 &\le&   \noo a^4\rrm_1 \noo
 \intt_0^1 \summ_{k=1}^n f_k^*(t)f_k(t)) \frac{d\mu(t)}{t}\rrm_{M_m}
 +
 \noo b^4\rrm  \noo \intt_0^1 \summ_{k=1}^n g_k(t)g_k(t)^*) \frac{d\mu(t)}{1-t}\rrm_{M_m}
 \kl 2 \pl .
%  \!\! &\le&   \! 2 \p .
 \mel
According to Lemma \ref{t1}, we deduce
 \begin{align*}
  \noo id: F_n\to OH_n\rrm_{cb} &\le\sqrt{2} \pl .\qedhere
  \end{align*}\qd

\begin{lemma} \label{2s}  The natural inclusion map $id: F_n^*\to OH_n$ has $cb$-norm less than $\sqrt{2}$.
\end{lemma} \lz

\begin{proof}[\bf Proof:] We have to consider a norm one element  $z\in
M_m(F_n^*)=CB(F_n,M_m)$. Let us denote by  $u_z:F_n\to M_m$ be the
corresponding complete contraction. Then we obtain coefficients
$z_k=u_z(f_k)\in M_m$ of $z$ with respect to the natural dual
basis $f_1^*,...f_n^*\in F_n^*$ satisfying
$f_k^*(f_j)=\delta_{kj}$. According to Wittstock's theorem there
exists a completely contraction $v:G_n\to M_m$. Since $G_n$ is a
quotient space $vQ:L_2^c(\nu_1,\ell_2^n)\oplus_1
L_2^r(\nu_2,\ell_2^n)\to M_m$ is a complete contraction. Thus
there are $x\in M_m(L_2^r(\nu_1,\ell_2^n))$ and $y\in
M_m(L_2^c(\nu_2,\ell_2^n))$ of norm less than $1$ such that
 \[ vQ(h^1,h^2) \lel \intt_0^1 \summ_{k=1}^n x_k(t)h^1_k(t) \frac{d\mu(t)}{t}
 + \intt_0^1 \summ_{k=1}^n y_k(t)h^2_k(t) \frac{d\mu(t)}{1-t}  \pl \]
for all $h^1\in L_2(\nu_1,\ell_2^n)$, $h^2\in
L_2(\nu_2,\ell_2^n)$.  Again, we can use the fact that $v$
vanishes on $kerQ$ and get
 \[   \frac{x_k(t)}{t} \lel \frac{y_k(t)}{1-t} \quad \mbox{$\mu$ a.e.}\]
for all $k=1,...,n$. In order to identify our original map $u_z$,
we compute
 \for
  z_k &=& u_z(x_k) \lel  vQ(e_k\sqrt{t},e_k(1-\sqrt{t})) \lel
  \intt_0^1 \sqrt{t} x_k(t) \pl \frac{d\mu(t)}{t} +
  \intt_0^1 (1-\sqrt{t}) y_k(t) \pl \frac{d\mu(t)}{1-t}\\
  &=&  \intt_0^1 x_k(t) \pl \frac{d\mu(t)}{t} \pl .
  \mel
Now, let us consider invertible positive elements $a,b\in S_4^m$.
As above, we define $A(z_k)_{k=1}^n=(z_ka^4)_{k=1}^n$ and
$B(z_k)_{k=1}^n=(b^4z_k)_{k=1}^n$ and
$\tilde{x}_k(t)=b^2x_k(t)a^{-2}$. Then, we get
 \for
  B^{\frac12}((z_k)_{k=1}^n) &=& (b^2z_k)_{k=1}^n \lel \intt_0^1 (b^2x_k(t))_{k=1}^n \pl  \frac{d\mu(t)}{t}
  \lel
  \intt_0^1 (b^2z_k(t)a^{-2})_{k=1}^n a^2  \pl  \frac{d\mu(t)}{t} \\
  &=&
  \intt_0^1 A^{\frac12}((\tilde{x}_k(t))_{k=1}^n)  \pl  \frac{d\mu(t)}{t}  \pl .
  \mel
On the other hand  for $\tilde{x}=(\tilde{x}_k)_{k=1}^n$, we
deduce
 \begin{align*}
 \intt_0^1 (A \tilde{x}(t),\tilde{x}(t)) \frac{d\mu(t)}{t}
 &=  \intt_0^1 (\tilde{x}(t)),A \tilde{x}(t))  \frac{d\mu(t)}{t}
 \lel
  \intt_0^1 \summ_{k=1}^n tr(\tilde{x}_k^*(t) \tilde{x}_k(t)a^4) \frac{d\mu(t)}{t} \\
 &=  \intt_0^1 \summ_{k=1}^n tr(a^{-2}x_k^*(t)b^2b^2 x_k(t)a^{-2}a^4) \frac{d\mu(t)}{t}\\
%   \mel \for
 &=  tr\bigg(b^4\intt_0^1 \summ_{k=1}^n x_k(t)x_k^*(t) \pl \frac{d\mu(t)}{t}\bigg) \\
 &\le \noo b^4\rrm_1 \pl \noo \intt_0^1 \summ_{k=1}^n x_k(t)x_k^*(t) \pl
 \frac{d\mu(t)}{t}\rrm_{M_m}\\
 &\le  \noo x\rrm_{M_m(L_2^r(\nu_1,\ell_2^n))}^2 \kl 1 \pl .
 \end{align*}
Similarly, we define $\tilde{y}_k(t)=b^{-2}y_k(t)a^2$,
$\tilde{y}(t)=(\tilde{y}_k(t))_{k=1}^n$ and get
 \for
  \frac{A((\tilde{x}_k(t))_{k=1}^n)}{t} &=& \frac{(b^2x_k(t)a^2)_{k=1}^n}{t}
  \lel \frac{(b^2y_k(t)a^2)_{k=1}^n}{1-t}
  \lel \frac{B((\tilde{y}_k(t))_{k=1}^n)}{1-t}  \quad \mbox{$\mu$-a.e.} \pl .
  \mel
As above we deduce
 \begin{align*}
 \intt_0^1 (B \tilde{y}(t),(\tilde{y}(t)) \frac{d\mu(t)}{1-t}
 &=  \intt_0^1 (\tilde{y}(t)),B\tilde{y}(t) )
 \frac{d\mu(t)}{1-t} \lel
  \intt_0^1 \summ_{k=1}^n tr(\tilde{y}_k^*(t) b^4\tilde{y}_k(t))
  \frac{d\mu(t)}{1-t} \\
  &=  \intt_0^1 \summ_{k=1}^n tr(a^{2}y_k^*(t)b^{-2}b^4b^{-2}
  y_k(t)a^{2})
 \frac{d\mu(t)}{1-t} \\
 &=
  tr\bigg(a^4\intt_0^1 \summ_{k=1}^n y_k^*(t)y_k(t) \pl
 \frac{d\mu(t)}{1-t}\bigg)\\
 &\le \noo a^4\rrm_1 \pl \noo \intt_0^1 \summ_{k=1}^n y_k^*(t)y_k(t) \pl \frac{d\mu(t)}{1-t}\rrm_{M_m} \\
 &\le \noo y\rrm_{M_m(L_2^c(\nu_1,\ell_2^n))}^2 \kl 1 \pl .
 \end{align*}
Therefore Lemma \ref{W2} implies
  \for
 \summ_{k=1}^n \noo bz_ka\rrm_2^2 &=& \noo (\sqrt{AB}(x_k))\rrm_2^2
 \kl  \intt_0^1 (A\tilde{x}(t),\tilde{x}(t)) \frac{d\mu(t)}{t} +
 \intt_0^1 (B\tilde{y}(t),\tilde{y}(t)) \frac{d\mu(t)}{1-t} \kl 2 \pl .
 \mel
Since, $a,b$ is arbitrary, we deduce from Lemma \ref{t1} that
 \[ \noo \summ_{k=1}^n z_k\ten \bar{z}_k\rrm^{\frac12} \kl \sqrt{2} \pl \noo \summ_{k=1}^n z_k\ten f_k^*\rrm_{M_m(F_n^*)} \pl .\]
This completes the proof.\qd

\begin{cor} $F_n$ is $2$ completely isomorphic to $OH_n$.
\end{cor} \lz

\begin{proof}[\bf Proof:] Since $OH_n$ is selfdual, see \cite{Poh}, it
suffices to observe that by Lemma \ref{W1} and Lemma \ref{W2}
 \begin{align*}
  \noo id:F_n\to OH_n\rrm_{cb} \noo id:OH_n\to F_n\rrm_{cb}
 &=  \noo id:F_n\to OH_n\rrm_{cb} \noo id:F_n^*\to OH_n\rrm_{cb}
 \kl  2  \pl . \qedhere
 \end{align*}\qd

\begin{rem}{\rm A similar result holds in the
context of $L_p$ spaces. Let us recall the operator spaces
$H^{c_p}=[H^c,H^r]_{\frac1p}$ and $H^{r_p}=[H^r,H^c]_{\frac1p}$.
Let $2\le p\le \8$ and $p'\le q\le p$ and
$\frac{1}{r}+\frac{1}{2p}=1$.  We consider the quotient space
$G_n=L_2^{c_p}(\nu_1;\ell_2^n)\oplus_q
L_2^{r_p}(\nu_2;\ell_2^n)/kerQ$ and $F_n^p=Q^{-1}(\{constants\})$
the subspaces corresponding to constant functions. Using
 \[ \noo \summ_{k=1}^n a_k\ten e_k\rrm_{S_p^m[OH_n]}
  \lel \sup_{\noo a\rrm_{2r}\le 1, \noo b\rrm_{2r}\le 1}
  \kla \summ_{k=1}^n \noo bx_ka\rrm_2^2 \mer^{\frac12} \]
we deduce as in in Lemma \ref{1s}  and in Lemma \ref{2s} that
 \[ \noo id:F_n^p\to OH_n\rrm_{cb} \kl
 2^{\frac12-\frac1p} \quad, \quad \noo id:(F_n^p)^*\to OH_n\rrm_{cb} \kl 2^{\frac12-\frac1p} \pl . \]
This yields
 \[ d_{cb}(F_n^p,OH_n) \kl 2^{1-\frac{2}{p}} \pl .\]
This is a concrete embedding of $OH$ as a subspace of a quotient
of $S_p=L_p(B(\ell_2),tr)$ (with $q=p$) and of $S_{p'}$ (with
$q=p'$ and using $H^{c_p}=H^{r_{p'}}$). For less concrete
realizations of such embedings we refer to \cite{P}. }\end{rem}
\hz

\begin{rem} {\rm A slight modification of this approach
yields the space $C_p=[C,R]_{\frac1p}$. Indeed, let
$\al=1-\frac{1}{p}$. Using $u=t^{-1}-1$, we have
 \[ \intt_0^1 \frac{1}{(1-t)+tB} \pl \frac{dt}{t^{\al}(1-t)^{1-\al}}  \lel \intt_0^\infty
 \frac{1}{u+B} \pl \frac{du}{u^{1-\al}}\lel
 \frac{c(\al)}{B^{1-\al}} \pl .\]
Following \cite{Wo}, we find for arbitrary commuting operators
$A$, $B$ that
 \begin{align*}
 &(x,A^{1-\al}B^{\al}x) \lel \inf_{x=f(t)+g(t)} \\
 &\quad \intt_0^1
 \frac{(f(t),Af(t))}{t} \frac{dt}{c(\al)t^\al(1-t)^{1-\al}}
 + \intt_0^1
 \frac{(g(t),Bg(t))}{1-t} \frac{dt}{c(\al)t^\al(1-t)^{1-\al}} \pl
 .\end{align*}
Similar as in Lemma \ref{1s}, it then easily follows that  the
space $E_{\al}$ of constants in the natural quotient of
$L_2^c(t^{-1}\mu_{\al})\oplus_1 L_2^r((1-t)^{-1}\mu_{\al})$
satisfies $\noo id:E_{\al}\to [C,R]_{\frac1p}\rrm_{cb}\le
\sqrt{2}$. The analogue Lemma \ref{2s} is slightly more involved.
The dual Pusz/Woronowicz formula has to be used for
$A(x_k)=(a_0^{\gamma_1}x_kb_0^{\beta_1})$ and
$B(x_k)=(a_0^{\gamma_2}x_kb_0^{\beta_2})$. The appropriate  powers
$\gamma_1,\gamma_2$ and $\beta_1,\beta_2$ are obtained by solving
solving linear equations. We leave the details to the interested
reader.}
\end{rem}\lz

Assuming Theorem 2, we can now obtain an embedding of $OH$.\lz

\begin{theorem}\label{oh}  $OH$ embeds into the predual of a von Neumann algebra with QWEP.
\end{theorem}\lz

\begin{proof}[\bf Proof:] Since $\bigcup_n OH_n$ is dense in $OH$, we
deduce that the subspace
 \[ F\subset L_2^c(\nu_1;\ell_2) \oplus_1   L_2^r(\nu_2;\ell_2)/kerQ\lel G \]
of all equivalence classes $(f,g)+kerQ$ such that $Q(f,g)$ is
$\mu$-almost everywhere constant is cb-isomorphic to $OH$.
Theorem 2 implies the assertion.\qd

\begin{rem}{\rm We refer to  \cite{JO} for an embedding in the predual of the hyperfinite III$_1$ factor.
}
\end{rem}\lz

\begin{cor} \label{ohf} There exists a constant $C>0$ such that for all $\nen$, there is an injective  linear map $u:OH_n\to S_1^m$ such that
\[ \noo u\rrm_{cb}\noo u^{-1}:u(OH_n)\to OH_n\rrm_{cb} \kl C \pl .\]
\end{cor}\lz

\begin{proof}[\bf Proof:] This is an immediate consequence of the strong
principle of local reflexivity in \cite{EJR} and the fact that
$N$ is QWEP.\qd

\section{The projection constant of the operator space $OH_n$}

In this section, we will provide the prove of Theorem 4 modulo the
probabilistic result Proposition \ref{kick00} (see section 7 for
the proof). The main tool used in this section is the trace
duality relation between $1$-summing norms and maps factorizing
though $B(H)$. Let us recall that for a  linear map $v:F\to E$
this factorization norm is defined as
 \[ \Gamma_\8(v) \lel \inf \noo \al \rrm_{cb}  \noo \beta
 \rrm_{cb}  \pl ,\]
where the infimum is taken over all $\iota_Ev=\al\beta$,
$\beta:F\to B(H)$ and $\al:B(H)\to E^{**}$. Here $\iota_E:E\to
E^{**}$ is the natural inclusion map. For finite rank maps
$v:F\to E$, we also define
 \[ \gamma_\8(v) \lel \inf \noo \al \rrm_{cb}  \noo \beta
 \rrm_{cb}  \pl ,\]
where the infimum is taken over all $m\in \nz$, $\beta:F\to  M_m$
and $\al:M_m\to E$ such that $v=\al\beta$. The $\gamma_\8$ norm
is related to the $1$-summing norm via trace duality.\lz

\begin{lemma}\label{trace} Let $u:E\to F$, then
 \[ \pi_1^o(u) \lel \sup\{  tr(vu) \pl |\pl \gamma_\8(v)\le 1\}
 \pl .\]
\end{lemma}\lz

\begin{proof}[\bf Proof:] Indeed, we have
 \begin{align*}
  \pi_1^o(u) &= \sup\left\{  |\langle id\ten u(x),y\rangle| \pl|\pl
  \noo x\rrm_{S_1^m\ten_{min} E}\le 1 \pl ,\pl \noo
  y\rrm_{M_m(F^*)}\le 1\right\}\\
  &= \sup \left\{ |tr(T_y^*uT_x)| \pl|\pl \noo T_y:M_m^*\to
  F^*\rrm_{cb} \le 1 \pl ,\pl \noo T_x:M_m\to E\rrm_{cb}\le 1\right\} \\
  &= \sup \left \{ |tr(vu)| \pl|\pl \gamma_\8(v)\le 1 \right\} \pl
  . \qedhere
 \end{align*}\qd

The following results from \cite{EJR} turns out to be quite
useful in our context.\lz

\begin{lemma}[EJR]\label{webfc} Let $E$ and $F$ be finite dimensional operator
spaces and $v:F\to E$, then
 \[ \gamma_\8(v)\lel \Gamma_\8(v) \pl .\]
\end{lemma}\lz

For subspaces of $L_1$ the one summing norm is closely related to
operator spaces projective norm and this connection is useful for
concrete estimates. For a finite dimensional operator space we
define
 \for
 d_{SL_1}(E) &=&
 \inf_{w:E\to E_1\subset S_1^m}  \noo w\rrm_{cb}\noo
 w^{-1}\rrm_{cb} \pl ,
 \mel
We extend  this to infinite dimensional spaces by
 \for
 d_{SL_1}(Y) &=& \sup_{E\subset Y}  d_{SL_1}(E) \pl
 ,
 \mel
where the supremum is taken over all finite dimensional
subspaces. The following fact follows immediately from the
definition.  \lz

\begin{lemma}\label{ptfd} Let $x\in X\ten Y$ and $\eps>0$, then
then there  exist finite dimensional spaces $E\subset X$ and
$F\subset Y$ such that
 \[ \noo x\rrm_{E\wet F} \kl (1+\eps) \noo
 x\rrm_{X\wet Y} \pl .\]
\end{lemma}\lz

\begin{lemma} \label{linf1}  Let $X$ and $Y$ be operator spaces, $E\subset X$ and $F\subset Y$
be finite dimensional subspaces. Let  $x\in E\ten F$ be a tensor
with associated linear map $T_x:E^*\to Y$, then
 \[ \pi_1^{o}(T_x)\kl d_{SL_1}(Y)\noo x\rrm_{X\wet
 Y} \pl. \]
\end{lemma}\lz

\begin{proof}[\bf Proof:] Let $w:F\to F_1\subset S_1^n$ be a linear map
and $w^{-1}:F_1\to E$ be completely contractive inverse. Let
$\beta:F\to M_m$ and $\al:M_m\to E^*$ be complete contractions. By
Wittstock's extension theorem, there is a complete contraction
$\hat{\beta}:S_1^n\to M_m$ such that $\hat{\beta}|_{F_1}=\beta
w^{-1}$. Then $\al \hat{\beta}$ corresponds to an element $z\in
E^*\ten_{min} M_n$ of norm less than one. The injectivity of the
projective tensor product on $S_1^n$ yields
 \for
  |tr(T_x\al\beta)| &=& |tr(\al \beta T_x)|
  \lel |tr(\al\hat{\beta}wT_x)| \lel |\langle z, id\ten w(x) \rangle|\\
  &\le& \noo z\rrm_{E^*\ten_{min} M_n} \noo id_E\ten w(x)  \rrm_{E\wet S_1^m } \\
  &=& \noo z\rrm_{E^*\ten_{min} M_n} \noo id_X\ten w(x)
  \rrm_{X\wet S_1^n} \\
  &\le& \noo w\rrm_{cb} \noo x\rrm_{X\wet F} \pl .
 \mel
Taking the infimum over all $w$, we may replace $\noo w\rrm_{cb}$
by $d_{SL_1}(F)$. Since $F\subset Y$ is arbitrary, we deduce the
assertion from Lemma \ref{ptfd}.\qd

We have  a  partial converse to this observation using the tools
from \cite{EJR,Fub}.
 \lz

\begin{lemma} \label{linf2}
Let $M$ and $N$ be von Neumann algebras such that $N$ is QWEP. Let
$X$ and $Y$ be operator spaces  and $w:Y\to M_*$ be a complete
contraction. Let $E$ be a finite dimension  subspace of $X$ and
$x\in E\ten Y$ with associated linear map $T_x:E^*\to Y$ and
$u:E\to N_*$ be a complete contraction, then
 \[ \pi_1^o(T_x:E^*\to Y) \gl \noo u\ten
 w(x)\rrm_{N_*\wet M_*} \pl .\]
\end{lemma}{\lz}

\begin{proof}[\bf Proof:] We will follow the ultraproduct approach from
\cite{Fub}. Let $N\subset B(H)$ be an embedding  and $\E_*:N_*\to
B(H)^*$ be a complete contraction such $\E_*(\phi)(x)=\phi(x)$ for
every functional $\phi\in N_*$ and $x\in N$. Using the strong
principle of local reflexivity we can find for every subset
$E_1\subset N_*$ and $F_\8\subset N$ a map $v_{E_1,F_\8}:E_1\to
S_1$ such that $\noo v_{E_1,F_\8}\rrm_{cb} \le 1+(dim(E))^{-1}$
and
 \[ v_{E_1,F_\8}(\phi)(x)\lel \phi(x) \]
for all $\phi\in E_1$ and $x\in F_\8$. Let $\U$ be an ultrafilter
refining the natural filtration by inclusion. Then, we define
$v:N_*\to \prod_\U S_1$ and $v(x)=(v_{E_1,F_\8}(x))_{(E,F)}$
whenever $x\in E_1$ and obtain a complete contraction. Let $\pi:
M\bar{\ten} N \to \prod M\ten B(\ell_2)$ be the diagonal
embedding. Let $\phi \in N_*$ and  $\psi\in M^*$. Let $x\in
N\bar{\ten}M$ and $T_x:M_*\to N$ be the associated linear map,
then we have
 \for
 \langle id\ten v(\psi\ten \phi), \pi(x)\rangle
 &=&  \lim_{E_1,F_\8} \langle \psi \ten v_{E_1,F_\8}(\phi),x\rangle \\
 &=&  \lim_{\phi\in E_1, T_x(\phi)\in F_\8} v_{E_1,F_\8}(\phi)(T_x(\psi))
  \lel \phi(T_x(\psi)) \\
 &=& \langle \psi\ten \phi, x\rangle \pl .
 \mel
Hence, $\pi^*(id_{M_*}\ten v)$ coincides with the natural
embedding of $M_*\wet N_*=(M\bar{\ten}N)_*$ into its bidual
$(M\bar{\ten}N)^*$. Therefore  $id\ten v:N_*\wet M_*\to \prod
S_1\wet M_*$ is (completely) isometric. Let us consider the
finite dimensional subspace $E_1=u(F)$. Given $\eps>0$, we can
find an index $(E_1',F_\8)$ with $E_1\subset E_1'$  and a map
$v=v_{E_1',F_\8}$ such that
 \[ \noo u\ten w(x)\rrm_{N_*\wet M_*} \kl
 (1+\eps)  \noo vu\ten w(x)\rrm_{S_1\wet M_*} \pl .\]
Hence, we can find a norm one element $y\in B(\ell_2)\bar{\ten}M$
such that
 \[ (1+\eps) |\langle y,vu\ten w(x)\rangle| \gl
 \noo vu\ten w(x)\rrm_{S_1\wet M_*} \pl .\]
However, by \cite{ER}, we can associate to $y$ a linear complete
contraction $T_y:M_*\to B(\ell_2)$ such that $T_y(\phi)(a)\lel
\langle y,\phi\ten a\rangle$. By testing with rank one tensors is
easy to see that
 \[ \langle y,vu\ten w(x)\rangle \lel tr( (vu)^*
 T_ywT_x) \pl .\]
Thus  with Lemma \ref{webfc} we conclude that
 \for
 |\langle y,vu\ten w(x)\rangle| &=&  |tr( (vu)^* T_ywT_x)|
 \kl \gamma_\8((vu)^* T_yw|_{Im(T_x)} )\pi_1^o(T_x)
 \\
 &\le& \Gamma_\8((vu)^* T_yw) \pi_1^o(T_x)
 \kl \noo v\rrm_{cb}\noo u\rrm_{cb} \noo T_y\rrm_{cb} \noo
 w\rrm_{cb} \pi_1^o(T_x) \pl .
 \mel
Since $u$, $w$ and $T_y$ are complete contractions and $\eps>0$ is
arbitrary, we deduce the assertion.\qd

Using free  probability we will be provide  an embedding theorem
for $G=L_2^c(\nu_1;\ell_2)\oplus_1 L_2^r(\nu_2;\ell_2)/ker Q$
into the predual of von Neumann algebra (see section 7). \lz

\begin{prop}\label{kick00}
There exists a von Neumann algebra $N$ with QWEP and a von
Neumann algebra $M$ and contractions $u:G\to N_*$ and $w:G\to
M_*$ such that for all $x\in G\ten G$
 \[ \noo u\ten w(x) \rrm_{N_*\wet  M_*} \gl \frac{1}{9} \noo
 x\rrm_{G\wet G} \pl \]
for all $x\in G\ten G$. Moreover, $\noo u^{-1}:u(G)\to
G\rrm_{cb}\le 3$.
\end{prop}\lz

\begin{cor} Let $E\subset G$ and $F\subset G$ be
finite dimensional subspaces and $x\in E\ten F$, then
 \[ \frac{1}{9}\noo x\rrm_{G\wet G}\kl \pi_1^o(T_x:E^*\to G) \kl 3 \noo x\rrm_{G\wet G}  \]
\end{cor}\lz

\begin{proof}[\bf Proof:] Since $G$ is $3$-cb isomorphic to a subspace of
$N_*$ and $N$ is QWEP, the strong principle of local reflexivity
in \cite{EJR} implies $d_{SL_1}(G)\le 3$ and thus by Lemma
\ref{linf1}
 \[  \pi_1^o(T_x:E^*\to G) \kl 3 \noo x\rrm_{G\wet G}
 \pl .
 \]
Conversely, Lemma \ref{linf2} applies due to Proposition
\ref{kick00} and hence
 \[ \noo x\rrm_{G\wet G} \kl 9 \noo u\ten w(x)\rrm_{N_*\wet M_*}
 \kl 9  \pl  \pi_1^o(T_x:E^*\to F)
 \pl .\]
The proof is complete.\qd

We are now well-prepared for the proof of Theorem 4.

\begin{proof}[\bf Proof of Theorem 4:] Let $F_n\to G$ be the subspace
constructed in section 3 and $u:OH_n\to OH_n$ be a linear map
with matrix $(a_{ij})$, then we deduce from Lemma \ref{1s} and
Lemma \ref{2s} that
 \for
  \pi_1^o(u:OH_n\to OH_n) &=&
  \pi_1^o(u:OH_n^*\to OH_n) \kl 2 \pl \pi_1^o(u:F_n^*\to
  F_n) \\
  &\le& 6 \pl \noo \summ_{i,j=1}^n a_{ij} f_i\ten
  f_j \rrm_{G\wet G} \pl .
  \mel
Conversely,
 \for
 \noo \summ_{i,j=1}^n a_{ij} f_i\ten
  f_j \rrm_{G\wet G}
  &\le& 9 \pl \pi_1^o(u:F_n^*\to F_n) \kl 18 \pl
  \pi_1^o(u:OH_n\to OH_n) \pl .
  \mel
The assertion is proved.\qd

The following norm calculations in $G_n\wet G_n$ will be postponed
to the next section.\lz

\begin{prop}\label{ncalc} Let $[a_{ij}]$ be a $n\times n$ matrix. Then
 \[ \noo \summ_{i,j=1}^n a_{ij} f_i\ten f_j \rrm_{G_n\wet G_n} \kl 16
 \pl
 \sqrt{1+\ln n} \kla \summ_{i,j=1}^n |a_{ij}|^2 \mer^{\frac12} \pl
 .\]
Moreover
 \[ \noo \summ_{i=1}^n  f_i\ten f_i \rrm_{G_n\wet G_n} \gl
 (8\pi)^{-1} \pl \sqrt{n(1+\ln n)} \pl .\]
\end{prop}\lz

As an application, we derive an independent proof of $(0.5)$. \lz

\begin{cor} Let $u:B(H)\to OH$ be a completely bounded map or rank at most $n$, then
 \[ \kla \summ_{k=1}^n \noo u(x_k)\rrm_2^2 \mer^{\frac12} \kl 96
 \pl
 \sqrt{1+\ln n} \pl \noo  u\rrm_{cb} \pl \noo \summ_{k=1}^n x_k
 \ten \bar{x}_k\rrm_{B(H)\ten_{min} \bar{B}(H)}^{\frac12} \pl .\]
\end{cor}\lz

\begin{proof}[\bf Proof:] We may assume that $u(x_1),...,u(x_n)$ are
contained in an $n$-dimensional subspace $K\subset OH$. Since $OH$
is homogeneous, it suffices to consider $Pu$, where $P$ is the
orthogonal projection on $K$. Choosing an orthonormal basis, we
may as well assume $K=OH_n$ and $u:B(H)\to OH_n$. We define
$w:OH_n\to B(H)$ by $w(e_k)=x_k$. Therefore, we deduce from Lemma
\ref{webfc}  and (\ref{ohnorm}) that $v=uw:OH_n\to OH_n$ satisfies
 \for
 \gamma_\8(v)&=&  \Gamma_\8(uw)\kl \noo u\rrm_{cb} \noo
 w\rrm_{cb} \lel
 \noo u\rrm_{cb} \pl \noo \summ_{k=1}^n x_k
 \ten \bar{x}_k\rrm_{B(H)\ten_{min} \bar{B}(H)}^{\frac12} \pl .
 \mel
Then we find  $(a_{ij})\in \ell_2^{n^2}$ of norm one such that
 \[ \kla \summ_{k=1}^n \noo u(x_k)\rrm_2^2 \mer^{\frac12}
 \lel \summ_{k,l} a_{kl}(e_l,uw(e_k)) \lel tr(auw) \pl .
 \]
Hence, we get
 \begin{align*}
 \kla \summ_{k=1}^n \noo u(x_k)\rrm^2 \mer^{\frac12}  &=
 |tr(av)| \kl \pi_1^o(a)  \gamma_\8(v)  \\
 &\le 6\pl 16 \pl \sqrt{1+\ln n} \pl \noo a\rrm_2 \pl
 \noo u\rrm_{cb} \pl \noo \summ_{k=1}^n x_k
 \ten \bar{x}_k\rrm_{B(H)\ten_{min} \bar{B}(H)}^{\frac12} \\
 &\le 96 \pl  \sqrt{1+\ln n} \pl \noo u\rrm_{cb} \pl \noo \summ_{k=1}^n x_k
 \ten \bar{x}_k\rrm_{B(H)\ten_{min} \bar{B}(H)}^{\frac12} \pl
 .\qedhere
 \end{align*}
\qd

Now, we apply a typical trace duality argument.\lz

\begin{cor} Let $\nen$. Then
 \[ \gamma_\8(id_{OH_n}) \kl 144\pi  \sqrt{\frac{n}{1+\ln n}} \pl .\]
\end{cor}\lz

\begin{proof}[\bf Proof:] Let $n\ge 4$. According to  Lemma \ref{trace},
Theorem 4. and Proposition \ref{ncalc}, we can find $v:OH_n\to
OH_n$ such that $\gamma_\8(v)\le 1$ and
 \for
 (8\pi)^{-1} \sqrt{n(1+\ln n)} &\le& \noo \summ_{i=1}^n f_i\ten f_i
 \rrm_{G\wet G} \kl 18 \pl \pi_1^o(id)  \pl \lel 18 \pl |tr(v)| \pl .
 \mel
Let $U_n$ be the unitary group in $\cz^n$ with normalized Haar
measure $\si$. Then, we consider
 \[  w\lel \intt_{U_n} u^*vu d\si_n(u) \lel \frac{tr(v)}{n} \pl  id \pl . \]
By homogeneity, we deduce
 \[ \frac{|tr(v)|}{n}\gamma_{\infty}(id) \lel  \gamma_\8(w)\kl \intt_{U_n} \gamma_\8(u^*vu)d\si_n(u) \kl
 \gamma_\8(v) \kl 1 \pl .\]
Therefore, we obtain the assertion
 \begin{align*}
 \gamma_\8(id) &= n|tr(v)|^{-1} \gamma_\8(w)
  \kl n \frac{8\pi}{\sqrt{n(1+\ln n)}} \pl 18
  \lel 144 \pi \pl  \sqrt{\frac{n}{1+\ln n}} \pl . \qedhere
 \end{align*} \qd

\begin{cor}\label{proj} Let $OH_n\subset B(\ell_2)$, then there exists
a projection $P:B(\ell_2)\to OH_n$ such that
 \[ \noo P\rrm_{cb} \kl 144 \pi \pl  \sqrt{\frac{n}{1+\ln n}} \pl .\]
\end{cor}\lz

\begin{proof}[\bf Proof:] It suffices to assume $n\ge 4$. Write
$id_{OH_n}=vw$, with $w:OH_n\to B(H)$ and  $v:B(H)\to OH_n$, and
 \[ \noo v\rrm_{cb} \noo w\rrm_{cb} \kl 144 \pi \pl  \sqrt{\frac{n}{1+\ln n}}
 \pl .\]
According to Wittstock's extension theorem, we can find a lifting
$\hat{w}:B(\ell_2)\to B(H)$ with the same cb-norm as $w$. Then
$P=v\hat{w}$ is the corresponding projection.\qd

\begin{cor} The order $\sqrt{(1+\ln n)}$  is  best possible in
$(0.4)$
\end{cor} \lz

\begin{proof}[\bf Proof:] Let $\iota:OH_n\to B(\ell_2)$ be a completely
isometric embedding and $x_k=\iota(e_k)$. According to Corollary
\ref{proj}, we can find a projection $P:B(\ell_2)\to OH_n$ of
cb-norm
 \[ \noo P\rrm_{cb} \kl C \sqrt{\frac{n}{1+\ln n}} \pl ,\]
where $C=96$. We define
 \[ c_n \lel \sup\left \{ \kla \summ_{k=1}^n \noo u(y_k)\rrm_2^2
 \mer^{\frac12} \pl \Bigg | \pl   \noo u\rrm_{cb}\le 1\pl ,\pl
 \noo \summ_{k=1}^n y_k\ten \bar{y}_k \rrm^{\frac12}_{B(H)\ten_{min}
 \bar{B}(H)}\le 1 \right \} \pl .\]
Let us consider the complete contraction $u=\frac{P}{\noo
P\rrm_{cb}}$. Since $\iota$  is a contraction, we deduce from
(\ref{ohnorm}) that
 \for
 \sqrt{n} &=& \kla \summ_{k=1}^n \noo e_k\rrm_2^2 \mer^{\frac12}
 \lel \kla \summ_{k=1}^n \noo P\iota(e_k)\rrm_2^2 \mer^{\frac12}
 \lel  \noo P\rrm_{cb}  \kla \summ_{k=1}^n \noo u(x_k) \rrm_2^2
 \mer^{\frac12} \mel \for
 &\le& C \sqrt{\frac{n}{1+\ln n}} \pll c_n \pl .
 \mel
Thus, we have $\sqrt{1+\ln n} \kl C c_n$ and the assertion is
proved.\qd

\section{Norm calculations in a quotient space}

We will now provide the  norm calculations in $G_n\wet G_n$ used
in the previous section.\lz

\begin{lemma}\label{start}
$G_n\stackrel{\wedge}{\ten} G_n$ is  isometrically  isomorphic the
to quotient space  of
\[ L_2(\nu_1\ten \nu_1,\ell_2^{n^2}) \oplus_1
  L_2(\nu_1,\ell_2^n)\ten_{\pi}  L_2(\nu_2,\ell_2^n)
  \oplus_1
  L_2(\nu_2,\ell_2^n)\ten_{\pi}  L_2(\nu_1,\ell_2^n)
  \oplus_1
  L_2(\nu_2\ten \nu_2,\ell_2^{n^2}) \]
with respect to
 \[ S \lel \{(f,g,h,k) \pl |\pl f(t,s)
 +g(t,s)+h(t,s)+k(t,s)\lel 0 \pll \mu\ten \mu \mbox{
 a.e.}\} \pl .\]
\end{lemma} \lz

\begin{proof}[\bf Proof:] By the properties of the projective operator
space tensor product, we have
 \for
 G_n\wet G_n &=& (L_2^c(\nu_1;\ell_2^n)\oplus_1
 L_2^r(\nu_1;\ell_2^n))\wet (L_2^c(\nu_1;\ell_2^n)\oplus_1
 L_2^r(\nu_1;\ell_2^n)) / ker(Q\ten Q) \pl .
 \mel
We note that $H^c \stackrel{\wedge}{\ten}K^c\lel H\ten_2 K\lel H^r
\stackrel{\wedge}{\ten}K^r$ and $H^c \stackrel{\wedge}{\ten}K^r
\lel H\ten_{\pi} K\lel H^r \stackrel{\wedge}{\ten}K^c$. Therefore
the properties  of $\oplus_1$ and $\wet$ imply that
 \for
  \lefteqn{ (L_2^c(\nu_1;\ell_2^n)\oplus_1
 L_2^r(\nu_1;\ell_2^n))\wet (L_2^c(\nu_1;\ell_2^n)\oplus_1
 L_2^r(\nu_1;\ell_2^n))} \\
  & & \lel  (L_2^c(\nu_1;\ell_2^n)\wet L_2^c(\nu_1;\ell_2^n)) \oplus_1
       (L_2^c(\nu_1;\ell_2^n)\wet L_2^r(\nu_2;\ell_2^n)) \\
  & & \pll \pll \pll     \oplus_1 (L_2^r(\nu_2;\ell_2^n)\wet L_2^c(\nu_1;\ell_2^n)) \oplus_1
    (L_2^r(\nu_2;\ell_2^n)\wet L_2^r(\nu_2;\ell_2^n)) \\
  & & \lel L_2(\nu_1\ten \nu_1;\ell_2^{n^2})\oplus_1
  (L_2(\nu_1;\ell_2^n)\ten_{\pi} L_2 (\nu_2;\ell_2^n)) \\
  & & \pll \pll \pll \oplus_1  (L_2(\nu_2;\ell_2^n) \ten_{\pi}  L_2(\nu_1;\ell_2^n)) \oplus_1
  L_2^r(\nu_2\ten \nu_2;\ell_2^{n^2}) \pl .
 \mel
Using $(1.1)$, we observe  that all the four components can be
represented by $\mu\ten \mu$ measurable functions. Applying
$Q\ten Q$ yields the assertion. \qd

\begin{cor}\label{ffunct} Let  $a=[a_{ij}]$ be a $n\times n$-matrix, then
 \for
  \noo \summ_{i,j=1}^n  a_{ij}f_i\ten f_j \rrm_{G_n\stackrel{\wedge}{\ten} G_n}
  &\ge& \frac{|tr(a)|}{\sqrt{n}}  \pl \sup \bet \intt f(t,s) \frac{d\mu(t)}{t}
  \frac{d\mu(s)}{s} \rag \pl ,
  \mel
where the supremum  is taken over all measurable functions
$(f,g,h,k)$ such that
 \[ \frac{f(t,s)}{ts} \lel \frac{g(t,s)}{(1-t)(1-s)}
  \lel \frac{h(t,s)}{t(1-s)} \lel
  \frac{k(t,s)}{(1-t)s} \pll  \mu\ten \mu \mbox{  a.e. }\]
and
 \begin{align}
 \max\{  \| f\|_{L_2(\nu_1\ten \nu_1)},
 \| g\|_{L_2(\nu_2\ten \nu_2)}\} &\le 1 \pl ,\label{sqn1}\\
  \max\{ \noo h\rrm_{L_2(\nu_1)\ten_{\eps}L_2(\nu_2)},
  \| k\|_{L_2(\nu_2)\ten_{\eps}L_2(\nu_1)}\} &\le \sqrt{n} \pl .
  \label{sqn}
 \end{align}
\end{cor}\lz

\begin{proof}[\bf Proof:] Let $(f,g,h,k)$ be given as above. Consider a
decomposition of $a$ in matrix valued functions
$a=a^1(t)+a^2(t)+a^3(t)+a^4(t)$ such that
 \for
 \lefteqn{ \noo a^1\rrm_{L_2(\nu_1\ten\nu_1;\ell_2^{n^2})}
  +  \noo
  a^2\rrm_{L_2(\nu_2\ten\nu_2;\ell_2^{n^2})}}\\
  & &  + \noo a^3\rrm_{L_2(\nu_1;\ell_2^n)\ten_{\pi}
 L_2(\nu_2;\ell_2^n)} +
 \noo a^4\rrm_{L_2(\nu_2;\ell_2^n)\ten_{\pi} L_2(\nu_1;\ell_2^n)} \kl
 (1+\eps) \pl \noo \summ_{i,j=1}^n a_{ij} f_i\ten f_j \rrm_{G\stackrel{\wedge}{\ten} G}
 \mel
Then the Cauchy Schwartz inequality implies
 \for
 \lefteqn{ \bet  \intt \summ_{i=1}^{n} a^1_{ii}(t,s) f(t,s)
 \frac{ d\mu(t)}{t} \frac{d\mu(s)}{s} \rag }\\
 & & \kl \kla \intt \summ_{i=1}^{n} |a^1_{ii}(t,s)|^2
 \frac{ d\mu(t)}{t} \frac{d\mu(s)}{s} \mer^{\frac12} \pl
 \kla  \intt \summ_{i=1}^{n} |f(t,s)|^2
 \frac{ d\mu(t)}{t} \frac{d\mu(s)}{s}
  \mer^{\frac12}\\
 & & \kl  \noo a^1\rrm_{L_2(\nu_1\ten
    \nu_1;\ell_2^{n^2})} \pl\sqrt{n} \pl  \noo f\rrm_{L_2(\nu_1\ten
    \nu_1;\ell_2^{n^2})} \kl \sqrt{n} \pl \noo
    a^1\rrm_2\pl .
    \mel
Similarly,
 \[
  \bet  \intt \summ_{i=1}^{n} a^2_{ii}(t,s) g(t,s)
    \frac{ d\mu(t)}{1-t} \frac{d\mu(s)}{1-s} \rag
    \kl \sqrt{n}\noo a^2\rrm   \pl .\]
For every operator $h:L_2\to L_2$, we recall $\noo h\ten
id_{\ell_2^n}\rrm=\noo h\rrm$. Hence, we deduce from trace
duality and \eqref{sqn}
 \begin{align*}
   \bet  \intt \summ_{i=1}^{n} a^3_{ii}(t,s) h(t,s)
     \frac{ d\mu(t)}{t} \frac{d\mu(s)}{1-s} \rag
    &\le   \noo h\ten id_{\ell_2^n}\rrm_{
    L_2(\nu_1;\ell_2^n)\ten_{\eps}
    L_2(\nu_2;\ell_2^n)} \pl
      \noo
    a^3\rrm_{L_2(\nu_1;\ell_2^n)\ten_{\pi}
    L_2(\nu_2;\ell_2^n)} \\
    &\le   \noo
    h\rrm_{L_2(\nu_1)\ten_{\eps}L_2(\nu_1)}  \pl \noo a_3\rrm  \kl
     \sqrt{n} \noo a^3\rrm \pl .
    \end{align*}
Similarly,
 \[    \bet  \intt \summ_{i=1}^{n} a^4_{ii}(t,s) k(t,s)
    \frac{ d\mu(t)}{1-t} \frac{d\mu(s)}{s} \rag \kl \noo k\rrm \pl
    \noo a^4\rrm  \kl \sqrt{n}\pl
    \noo a^4\rrm  \pl .\]
Therefore, we get
 \begin{align*}
   \bet \intt  \summ_{i=1}^n a_{ii} f(t,s)
 \frac{d\mu(t)}{t}\frac{d\mu(s)}{s}\rag
% & & \lel  \intt \summ_{i=1}^n f(t,s)
% \frac{d\mu(t)}{t}\frac{d\mu(s)}{s}\\
 &\le \bet \intt \summ_{i=1}^n a^1_{ii}f(t,s)
 \frac{d\mu(t)}{t}\frac{d\mu(s)}{s} \rag +
 \bet \intt \summ_{i=1}^n a^2_{ii}g(t,s)
 \frac{d\mu(t)}{1-t}\frac{d\mu(s)}{1-s} \rag \\
 & \pll +\bet \intt \summ_{i=1}^n a^3_{ii}h(t,s)
 \frac{d\mu(t)}{t}\frac{d\mu(s)}{1-s} \rag +
 \bet \intt \summ_{i=1}^n a^4_{ii}g(t,s)
 \frac{d\mu(t)}{1-t}\frac{d\mu(s)}{s} \rag \\
 &  \le  \sqrt{n}
 (1+\eps)
 \noo \summ_{i,j=1}^n a_{ij} f_i\ten f_j \rrm_{G\stackrel{\wedge}{\ten} G} \pl .
 \end{align*}
The assertion follows by letting $\eps\to 0$. \qd

Given a measure $\nu$ and positive measurable densities $g,h$, we
denote by
 \[ L_2(g\nu)+_p L_2(h\nu) \lel L_2(g\nu)\oplus_p L_2(h\nu)/ker Q\]
where $Q(f_1,f_2)=f_1+f_2$. Given a measurable function $k$, we
define the norm of $k$ in $L_2(g\nu)+_p L_2(h\nu)$ as the norm of
the equivalence class $(k,0)+ker Q$.  For $p=2$, this is a again a
Hilbert space and we can find an explicit formula.\lz

\begin{lemma} \label{L2q}  Let $\nu$ be a measure and $g,h$
strictly  positive measurable functions. For a measurable
function $k$ the norm of in $L_2(g\nu)+_2L_2(h\nu)$ is given by
 \[  \noo  k\rrm_{L_2(g\nu)+_2L_2(h\nu)}
  \lel \kla \intt \frac{|k|^2}{(g^{-1}+h^{-1})} d\nu
  \mer^{\frac12} \pl .\]
\end{lemma}\lz

\begin{proof}[\bf Proof:] By the  Hahn-Banach theorem (see the proof of
Lemma \ref{W2}), we have
 \[ \noo k\rrm_{L_2(g\nu)+L_2(h\nu)}
 \lel \sup |\intt k f_1g  d\nu| \pl \]
where the supremum is taken over $(f_1,f_2)$ with
 \[ f_1 g \lel f_2 h \quad \mbox{and}
 \quad \intt |f_1|^2 g\nu + \intt |f_2|^2h d\nu \kl 1\pl
 .\]
We introduce $v=f_1g$ and find $f_1=g^{-1}v$, $f_2=h^{-2}v$ and
thus
 \for
  \noo k\rrm_{L_2(g\nu)+L_2(h\nu)} &=&
  \sup_{\intt |v|^2g^{-1}d\nu + \intt |v|^2h^{-1}d\nu\le
  1} \bet \intt k v d\nu \rag \\
 &=&   \sup_{\intt |w|^2 d\nu \le 1} \bet \intt k
  \frac{w}{(g^{-1}+h^{-1})^{\frac12}} d\nu \rag \lel
 \kla \intt \frac{k^2}{g^{-1}+h^{-1}}   d\nu
 \mer^{\frac12} \pl .
 \mel
The last equality is the equality case in H\"older's
inequality.\qd

\begin{cor} \label{calc} Let $0<\delta<\frac12$, then
 \begin{eqnarray}
 \noo 1_{[\delta,\frac12]}\ten
  1_{[\frac12,1-\delta]}\rrm_{L_2(\nu_1\ten
  \nu_1)+_1L_2(\nu_2\ten \nu_2)} &\le&  4\sqrt{2} \pi^{-1} (-\ln \delta)^{\frac12} \\
 \noo 1_{[\frac12,1-\delta]}\ten
  1_{[\delta,\frac12]}\rrm_{L_2(\nu_1\ten
  \nu_1)+_1L_2(\nu_2\ten \nu_2)} &\le&
  4\sqrt{2} \pi^{-1} (-\ln \delta)^{\frac12} \\
  \noo
  1_{[0,\frac12]}\ten
  1_{[0,\frac12]}\rrm_{L_2(\nu_1\ten
  \nu_1)+_1L_2(\nu_2\ten \nu_2)} &\le& 2\sqrt{2}\\
 \noo
  1_{[\frac12,1]}\ten
  1_{[\frac12,1]}\rrm_{L_2(\nu_1\ten
  \nu_1)+_1L_2(\nu_2\ten \nu_2)} &\le& 2\sqrt{2} \pl
  .
 \end{eqnarray}
\end{cor}

\begin{proof}[\bf Proof:] Recall that $d\nu_1(t)=\frac{1}{t}d\mu(t)$ and
$d\nu_2(t)=\frac{1}{1-t}d\mu(t)$ and $d\mu(t)=\frac{dt}{\pi
\sqrt{t}\sqrt{1-t}}$. Applying Lemma \ref{L2q}, we get
  \for
  \lefteqn{ \noo 1_{[\delta,\frac12]}\ten
  1_{[\frac12,1-\delta]}\rrm_{L_2(\nu_1\ten
  \nu_1)+_1L_2(\nu_2\ten \nu_2)}  }\\
  & &    \kl   \sqrt{2} \kla \intt_{\delta}^\frac12
  \intt_{\frac12}^{1-\delta} \frac{1}{ts+(1-t)(1-s)}
  d\mu(t)d\mu(s) \mer^{\frac12} \\
  & & \lel  \sqrt{2}
  \kla \intt_{\delta}^\frac12
  \intt_{\delta}^{\frac12} \frac{1}{t(1-s)+(1-t)s}
  d\mu(t)d\mu(s) \mer^{\frac12} \\ %\mel  \for
   & & \kl  2 \pl
  \kla \intt_{\delta}^\frac12
  \intt_{\delta}^{\frac12} \min(\frac1t,\frac1s)
  d\mu(t)d\mu(s) \mer^{\frac12}  \kl
  2 \kla \frac{4}{\pi^2} \intt_{\delta}^{\frac12}
  \intt_{\delta}^t \frac{ds}{\sqrt{s}} \pl
  \frac{dt}{t\sqrt{t}} \mer^{\frac12} \\
  & & \kl
   2  \kla \frac{8}{\pi^2}
   \intt_{\delta}^{\frac12}\sqrt{t} \pl
  \frac{dt}{t\sqrt{t}} \mer^{\frac12} \kl
  4\sqrt{2} \pi^{-1} (-\ln \delta)^{\frac12}
   \pl .
   \mel
By symmetry, we deduce $(5.2)$. Equation $(5.3)$ and $(5.4)$
follow from
 \begin{align*}
  \intt_0^{\frac12} \intt_0^{\frac12}
  \frac{1}{ts+(1-t)(1-s)} d \mu(t) d\mu(s)
  &\le 4 \intt_0^{1} \intt_0^{1} d\mu(t) d\mu(s)
  \lel 4 \pl .\qedhere
  \end{align*}\qd

The next estimate yields the upper estimate in the logarithmic
`little Grothendieck inequality'.\lz

\begin{lemma}\label{logest} Let $a$ be a $n\times n$
matrix, then
 \[ \noo \summ_{ij=1}^n a_{ij} f_i\ten f_j \rrm_{
 G_n\stackrel{\wedge}{\ten} G_n} \kl 16  \pl \sqrt{1+\ln n}
  \pl
  \noo a\rrm_2\pl .\]
\end{lemma}

\begin{proof}[\bf Proof:] Given $a\in \ell_2^{n^2}$ and
$0<\delta<\frac12$, we decompose
 \[ a \lel a^1(t,s) + a^2(t,s) \pl ,\]
where
 \for
 \lefteqn{  a^1(t,s) \lel  a\ten \bigg(1_{[0,\frac12]}(t)
  1_{[0,\frac12]}(s)+ 1_{[\delta,\frac12]}(t) 1_{[\frac12,1-\delta]}(s)}\\
  & & \pll \pll \pll  \pll \pll \pll  \pll \pll \pll
  +  1_{[\frac12,1-\delta]}(t)
 1_{[\delta,\frac12]}(s)+
 1_{[\frac12,1]}(t)
  1_{[\frac12,1]}(s)\bigg)
  \mel
and
 \[ a^2(t,s) \lel a\ten 1-a^1(t,s)  \pl .\]
According to Corollary \ref{calc}, we get
 \for
  \noo a^1\rrm_{L_2(\nu_1\ten \nu_1;\ell_2^{n^2})+_1
  L_2(\nu_2\ten \nu_2;\ell_2^{n^2})}
  &\le& (4\sqrt{2} + 8\sqrt{2} \pi^{-1} (-\ln
  \delta)^{\frac12}) \noo a\rrm_2 \pl .
  \mel
In order to estimate $a^2$, we note that
 \for
  \| 1_{[0,\delta]}\ten 1_{[\frac12,1]}
  \|_{L_2(\nu_2)\ten_{\pi} L_2(\nu_1)} &=&
   \| 1_{[0,\delta]} \|_{L_2(\nu_2)} \|
   1_{[\frac12,1]}\|_{L_2(\nu_1)}  \kl
   \frac{2^{\frac54}\delta^{\frac14}}{\sqrt{\pi}}  \pl
   2^{\frac12} \pl .
   \mel
Similarly, we get
 \for
 & &\max\{  \|1_{[\delta,\frac12]}\ten 1_{[1-\delta,1]}
  \|_{L_2(\nu_2)\ten_{\pi} L_2(\nu_1)},  \| 1_{[\frac12,1]}\ten 1_{[0,\delta]}
  \|_{L_2(\nu_1)\ten_{\pi} L_2(\nu_2)}, \\
 & & \pll \pll\pll \pl  \| 1_{[1-\delta,1]}\ten 1_{[\delta,\frac12]}
   \|_{L_2(\nu_1)\ten_{\pi} L_2(\nu_2)}\}
 \kl  \frac{2\pl 2^{\frac34}}{\sqrt{\pi}} \pl\delta^{\frac14} \pl .
 \mel
Using $L_2(\ell_2^n)\ten_\pi L_2(\ell_2^n)= (L_2^c\wet L_2^r)\wet
((\ell_2^n)^c\wet (\ell_2^n)^r)$, we get
 \for
  \noo a_2\rrm_{L_2(\nu_1;\ell_2^n)\ten_{\pi}
 L_2(\nu_2;\ell_2^n)+
 L_2(\nu_2;\ell_2^n)\ten_{\pi}
 L_2(\nu_1;\ell_2^n)}
 &\le& \frac{8\pl 2^{\frac34}}{\sqrt{\pi}} \pl \delta^{\frac14}
 \noo a\rrm_1 \kl
 \frac{8\pl 2^{\frac34}}{\sqrt{\pi}} \pl \delta^{\frac14}
 \pl \sqrt{n} \noo a\rrm_2 \pl .
 \mel
We can choose $\delta=\frac{1}{e^2n^2}$ and deduce the
assertion.\qd

Although this decomposition seems to be of technical nature, it
turns out to be essentially  optimal for the identity matrix. \lz

\begin{lemma} \label{lest}
 \[ \noo \summ_{i=1}^n f_i\ten f_i\rrm_{G_n\stackrel{\wedge}{\ten} G_n} \gl
 \frac{1}{8\pi}
 \sqrt{n(1+\ln n)} \pl \]
holds for all $\nen$.
\end{lemma}

\begin{proof}[\bf Proof:]
Let $\delta>0$ to be determined later. We consider the interval
$I=[\delta,\frac12]\times [\frac12,1-\delta]$ and the function
 \[ v(t,s)\lel  \frac{1}{ts+(1-t)(1-s)}  \pl 1_I \pl .\]
We are looking for functions $(f,g,h,k)$ as in Corollary
\ref{ffunct}. We define  $f(t,s)\lel ts v(t,s)$ and
$g(t,s)=(1-t)(1-s)v(t,s)$ and observe
 \begin{align*}
  \intt_I f(t,s)^2 \frac{d\mu(t)}{t}\frac{d\mu(s)}{s}
 &+ \intt_I g(t,s)^2
 \frac{d\mu(t)}{1-t}\frac{d\mu(s)}{1-s}
  \lel    \intt_I v(t,s)^2[ts+(1-t)(1-s)] d\mu(t)d\mu(s)\\
 &=   \intt_I \frac{1}{ts+(1-t)(1-s)}  d\mu(t)d\mu(s) \\
 & \le   4\pi^{-2} \intt_{\delta}^{\frac12}
 \intt_{\frac12}^{1-\delta} \min(t^{-1},(1-s)^{-1})
 \frac{dt}{\sqrt{t}} \frac{ds}{\sqrt{1-s}} \\
 %\mel\for
 & \le   4\pi^{-2} \intt_{\delta}^{\frac12}
 \intt_{\delta}^{\frac12} \min(t^{-1},s^{-1})
 \frac{dt}{\sqrt{t}} \frac{ds}{\sqrt{s}}
  \lel  8\pi^{-2} \intt_{\delta}^{\frac12}
 \intt_{\delta}^t \frac{ds}{\sqrt{s}}  \pl
 \frac{dt}{t\sqrt{t}} \\
 & \le  16 \pi^{-2} \intt_{\delta}^{\frac12} \sqrt{t}
 \frac{dt}{t\sqrt{t}} \kl 16 \pi^{-2} (-\ln \delta)
 \pl .
 \end{align*}
In order to  estimate the norm for $h(t,s)=t(1-s)v(t,s)$, we use
$(1.1)$ and hence it suffices to estimate the $L_2$-norm:

 \begin{align*}
  \noo h\rrm_{L_2(\nu_1\ten \nu_2)}^2 &\le  4 \intt_I
   \min(t^{-2},(1-s)^{-2}) \pl \pl t^2(1-s)^2 \pl
   d\nu_1(t) d\nu_2(s)  \\
  &\le 8\pi^{-2} \intt_{\delta}^{\frac12} \intt_{\delta}^{\frac12}
  \min(t^{-2},s^{-2}) ts \frac{dt}{\sqrt{t}}
  \frac{ds}{\sqrt{s}} \\
  &= 16\pi^{-2} \intt_\delta^{\frac12}
  \intt_{\delta}^t  \sqrt{s} ds \p
  \frac{dt}{t\sqrt{t}} \kl
 \frac{32}{3\pi^2} \intt_{\delta}^{\frac12}
  t^{\frac32} \frac{dt}{t\sqrt{t}}
  \lel \frac{16}{3\pi^2} \pl .
  \end{align*}
Finally, we need the $L_2$-norm estimate of $k(t,s)=(1-t)sv(t,s)$.
 \for
  \noo k\rrm^2_{L_2(\nu_2\ten \nu_1)}
  &\le& 4 \intt_I  \min(t^{-2},(1-s)^{-2}) \pl \pl  (1-t)^2 s^2 \pl
  \frac{d\mu(t)}{(1-t)} \frac{d\mu(s)}{s}\\
  &\le&  8 \pi^{-2} \pl \intt_{\delta}^{\frac12} \intt_{\delta}^{\frac12}
   \min(t^{-2},s^{-2}) \frac{dt}{\sqrt{t}}
   \frac{ds}{\sqrt{s}} \kl
 16  \pi^{-2} \pl \intt_{\delta}^{\frac12} \intt_{\delta}^{t}
  \frac{ds}{\sqrt{s}}
    \frac{dt}{t^2\sqrt{t}} \\
    &\le&     32 \pi^{-2} \pl \intt_{\delta}^{\frac12}
  \sqrt{t}    \frac{dt}{t^2\sqrt{t}}
  \kl  32 \pi^{-2} \delta^{-1} \pl .
  \mel
We note that
 \for
 \lefteqn{ \intt_I f(t,s) \frac{d\mu(t)}{t}\frac{d\mu(s)}{s}
 \lel
  \intt_{\delta}^{\frac12}\intt_{\frac12}^{1-\delta} \frac{1}{ts+(1-t)(1-s)} d\mu(t)d\mu(s)
  }\\
 & & \gl  \frac{1}{2\pi^{2}}  \pl  \intt_{\delta}^{\frac12}\intt_{\delta}^{\frac12}
 \min(t^{-1},s^{-1})  \frac{dt}{\sqrt{t}} \frac{ds}{\sqrt{s}} \lel
  \frac{2}{\pi^{2}}  \pl
 \intt_{\delta}^{\frac12} (\sqrt{t}-\sqrt{\delta})
    \frac{dt}{t\sqrt{t}} \\
  & &   \gl  \frac{1}{\pi^{2}}  \pl  \intt_{4\delta}^{\frac12}
 \frac{dt}{t} \lel \frac{(-\ln 8\delta)}{\pi^2} \pl
 .
 \mel
We define $\delta=\frac{1}{8ne}$, $C=\frac{4}{\pi}$ and
$\tilde{f}=\frac{f}{C\sqrt{-\ln\delta}}$. Then \eqref{sqn1} and
\eqref{sqn} are satisfied for corresponding quadrupel
$(\tilde{f},\tilde{g},\tilde{k},\tilde{h})$ and the assertion
follows from $\ln8\le 3$ which implies $\sqrt{-\ln\delta}\le
2\sqrt{1+\ln n}$ and thus $\frac{(-\ln
8\delta)}{\pi^2C\sqrt{\ln(-\delta)}}\gl \frac{\sqrt{1+\ln
n}}{8\pi}$.
 \qd

\begin{proof}[\bf Proof of Proposition \ref{ncalc}:] Combine Lemma
\ref{logest} and Lemma \ref{lest}.\qd

\section{K-functionals associated to states and conditional expectations}

In this section, we will clarify certain basic properties of
quotient spaces $K(N,d)$, $\kz_t(N,d)$ and $\kz_n(\M,\E)$. These
spaces play an important role in our investigation and provide a
link between the space $G$ considered in the previous sections
and the free probability techniques in the following sections. It
is possible to define these K-functionals in the more general
context of Haagerup's $L_p$ spaces, but semifinite von Neumann
algebras are sufficient for our applications. Therefore in the
following, we will assume that $\tau$ is a normal, faithful,
semifinite \emph{trace} on a von Neumann algebra $N$. We denote by
$L_0(N,\tau)$ the space $\tau$-measurable operators affiliated to
$N$. Given a positive density $d\in L_0(N,\tau)$ we define
$Q:L_2^c(N,\tau)\oplus_1 L_2^r(N,\tau)\to L_0(N,\tau)$ by
 \[ Q(x,y)\lel dx+yd  \pl .\]
On the vector space $Im(Q)\subset L_0(N,\tau)$ we define the
operator space structure of the quotient
 \[ K(N,d) \cong L_2^c(N,\tau)\oplus_1
 L_2^r(N,\tau)/kerQ \pl .\]
Note that $kerQ$ is closed by continuity of the left and right
multiplication in $L_0(N,\tau)$ and the continuity of the
inclusion $L_2(N,\tau)\subset L_0(N,\tau)$. \lz

\begin{lemma}\label{dlimit1} Let $(e_\la)$ be an increasing family
of projection converging strongly to $1$ and such that
  \[ e_\la d \lel e_\la d e_\la \lel  de_\la  \pl .\]
Then $K(e_\la Ne_\la,e_\la de_\la)$ is completely complemented in
$K(N,d)$ and  $K(N,d)$ is a direct limit of the $K(e_\la
Ne_\la,e_\la de_\la)$'s.
\end{lemma}\lz

\begin{proof}[\bf Proof:] Let $N_\la=e_{\la}Ne_{\la}$ and $d_\la= e_\la
de_\la$ and $Q_\la:L_2^c(N_{\la})\oplus_1 L_2^r(N_{\la})\to
K(N_{\la},d_{\la})$ be the corresponding quotient map. By
homogeneity, we see that $T_\la: L_2^c(N)\oplus_1 L_2^r(N)\to
L_2^c(N_{\la})\oplus_1 L_2^r(N_{\la})$ given by
$T_{\la}(x_1,x_2)=(e_{\la}x_1e_{\la},e_{\la}x_2e_{\la})$ is a
complete contraction. Moreover, $Q_{\la}T_{\la}$ vanishes on $ker
Q$ and thus induces a complete contraction
$\hat{T}_{\la}:K(N,d)\to K(N_{\la},d_{\la})$. Let
$\iota_{\la}:L_2^c(N_{\la},\tau)\oplus_1 L_2^r(N_{\la},\tau)\to
L_2^c(N,\tau)\oplus_1 L_2^r(N,\tau)$ be the canonical inclusion
map. Then $Q\iota_{\la}$ vanishes on $kerQ_{\la}$ and thus
induces a complete contraction
$\hat{\iota}_{\la}:K(N_{\la},d_{\la})\to K(N,d)$. It is easily
checked that $\hat{P}_{\la}\hat{\iota}_{\la}=id$ and this
completes the proof of  the first assertion. Since $e_\la\to 1$
strongly, we see that $\bigcup L_2(e_\la N e_{\la},\tau)$ is norm
dense in $L_2(N,\tau)$. By continuity of $Q$ this implies that
$\bigcup_{\la} Im(Q\iota_{\la})$ is dense in $K(N,d)$. Moreover,
if $e_{\la}\le e_{\la'}$ we have $Im(Q\iota_{\la})\subset
Im(Q\iota_{\la'})$. Hence $K(N,d)$ is a direct limit of the
$K(N_{\la},d_{\la})$'s. \qd

The next technical lemma will allow us to use approximation
arguments. Let us say that \emph{a  density $d$ is in
$L_2^{ap}(N,\tau)$} if $d$  is $\tau$ measurable and there exists
an increasing net $e_\al$ such that $e_{\al}d=de_{\al}$ and
$e_{\al}d\in L_2(N,\tau)\cap N$, $(e_{\al}d)^{-1}$ a bounded
element in $e_{\al}Ne_{\al}$. \lz

\begin{lemma}\label{wstar} Let $N$ be a von
Neumann algebra with normal faithful trace $\tau$ and $d$ a
positive density in  $L_2^{ap}(N,\tau)$ and $N_d$ be the subspace
of elements $n\in N$ such that $\tau(d^2(n^*n+nn^*))<\infty$.
Then $n\in \bar{N}_d$ defines a continuous functional $\phi_n$ by
\begin{eqnarray} \label{dform}
  \phi_n((a,b)+kerQ)) &=&  \tau(n^*(da+bd)) \pl .
\end{eqnarray}
For every von Neumann algebra $M\ten \bar{N_d}$ is weak$^*$ dense
in the unit ball of $M \ten_{min} K(N,d)^*$. If in addition
$d,d^{-1}$ are bounded and $d$ in $L_2(N,\tau)$, then $N=N_d$ and
$M\ten \bar{N}$ is norm dense in the unit ball $B$ of
$M\ten_{min} K(N,d)^*$.
\end{lemma}

\begin{proof}[\bf Proof:] According to our assumption and Lemma
\ref{dlimit1}, we know that $\bigcup_{\la}
K(e_{\la}Ne_{\la},e_{\la}d)$ is norm dense in $K(N,d)$. In view
of $M\ten_{min}K(N,d)^*\subset (M_*\wet K(N,d))^*$, it therefore
suffices to show the second assertion under the additional
assumption that $d\in N\cap L_2(N,\tau)$ and $d^{-1}$ is bounded.
By approximation it suffices to consider a finite rank tensor
$x=\sum_{j=1}^n m_j \ten \gamma_j$ of norm less than one in $M
\ten_{min} K(N,d)^*$. By changing the $m_j$'s and $\gamma_j$'s, we
may assume that the $(m_j)$'s are part of a biorthogonal basis,
i.e. there are functionals $(m_i)^*_{i=1,...n}$ such that
$m_i^*(m_j)=\delta_{ij}$ and $\|m_i^*\|=\|m_j\|=1$ for all $1\le
i,j\le n$. Using the anti-linear duality $\tau(x^*y)$, we have
$\overline{L_2(N,\tau)}^*=L_2(N,\tau)$. Therefore, we may assume
that the functional $\gamma_j\in K(N,d)^*$ is given by elements
$(x_j,y_j)$ such that
 \[ \gamma_j((a,b)+kerQ)) \lel \tau(x_j^*a) + \tau(y_j^*b) \pl.\]
Since $\gamma_j$ vanished on $kerQ$, we deduce from $da+bd=0$
that $\gamma_j((a,b)+kerQ)=0$. Given  $z\in L_2(N,\tau)$ we
define  $a=zd$, $b=-dz$ and obtain
 \[ 0\lel \tau(x_j^*zd)- \tau(y_j^*dz) \lel
 \tau((dx_j^*-y_j^*d)z) \pl .\]
Since $z$ is arbitrary, we deduce
 \begin{equation} \label{eqf3}
 x_jd \lel dy_j \pl . \end{equation}
for all $j=1,..,n$. Since $K(N,d)=L_2^c(N,\tau)\oplus_1
L_2^r(N,\tau)$ is a quotient space, it follows from the
definition and $H_c^*=H_r$, $H_r^*=H_c$ that
 \begin{eqnarray}\label{dualnorm}
 \pll \pll \pl \noo x\rrm_{M\ten_{min}K(N,d)^*} \!\!\!&=& \!\!\! \max\left\{\noo  \summ_j m_j\ten
  x_j\rrm_{M\ten_{min} L_2^r(N,\tau)},
  \noo  \summ_j m_j\ten
  y_j\rrm_{M\ten_{min} L_2^c(N,\tau)} \right\} \pl .
 \end{eqnarray}
Now, we define  $w_j=d^{-1}x_j=y_jd^{-1}$ (see \eqref{eqf3}).
Given $\eps>0$, we can find elements $v_j\in N$ such that $\|
v_j-w_j\|_2\le \frac{\eps}{n\|d\|_\8}$ for $j=1,..,n$. Let us
define
\begin{eqnarray}
 \hat{\gamma}_j(a,b) &=&  \tau(v_j^*da)+\tau(dv_j^*b)
  \lel \tau(v_j^*(da+bd)) \pl .
\end{eqnarray}
Clearly, $\gamma_j$ vanishes on $kerQ$ and moreover,
 \for
 \lefteqn{ \noo \summ_{j=1}^n m_j\ten
 (\gamma_j-\hat{\gamma}_j)\rrm_{M\ten_{min}
 K(N,d)^*}}\\
  & & \kl  \summ_{j=1}^n \noo m_j\rrm
  \max\{\noo x_j^*-v_j^*d\rrm_2,\noo y_j^*-dv_j^*\rrm_2\} \\
  & & \kl  \summ_{j=1}^n
  \max\{\noo w_j^*d-v_j^*d\rrm_2,\noo dw_j^*-dv_j^*\rrm_2\}
  \kl \eps \pl .
 \mel
Thus $(1-\eps)^{-1}\summ_{j=1}^n m_j\ten \hat{\gamma}_j$ is a
good approximation in $M\ten_{min} \bar{N}$. \qd

\begin{rem} \label{wstar2} {\rm Obviously,
$M\ten_{min}K(N,d)^*$ is a subspace of
 \[ (M_*\wet K(N,d))^* \lel CB(M_*,K(N,d)^*) \pl .\]
If moreover,  $M$ is finite dimensional these spaces coincide.
Let $B$ be the unit ball of the space  $CB(M_*,K(N,d)^*)$. For
general $M$, it is still true that $M\ten K(N,d)^*\cap 3B$ is
weak$^*$ dense in $B$. Indeed,  by the results in section 7,
$K(N,d)$ is $3$ completely complemented in $\M_*$ for some von
Neumann algebra $\M$ and thus $CB(M_*,K(N,d)^*)$ is $3$
complemented in $M\bar{\ten} \M$. By Kaplansky's density theorem
(see \cite{EJR} for this type of application), the unit ball of
$M\ten_{min}\M$ is strongly dense in the unit ball of
$M\bar{\ten}\M$. Using the complementation, the assertion
follows.}
\end{rem}\lz

\begin{lemma} \label{dykprep} Let  $N$ be a hyperfinite, semifinite von Neumann with normal semifinite trace $\tau$ and
$d\in L_2^{ap}(N,\tau)$. Then $K(N,d)$ is completely
contractively complemented in
\[   \prodd_{\al,\U} K(A_{\al},d_\al)\]
where $A_\al$ are finite dimensional algebras. If moreover, $N$
is given as $N=L_\8(\Om,\mu;M_n)$ and $(\Om,\mu)$ is a
$\si$-finite measure space, then every density $d\in L_0(N,\tau)$
is in $L_2^{ap}$ and  one can choose $A_\al$'s  of the form
$\ell_\8^{m_\al}(M_n)$.
\end{lemma}\lz

\begin{proof}[\bf Proof:] Again it suffices to prove the assertion under
the assumption that $d$ is an invertible bounded element  with
$d\in L_2(N,\tau)$. Let $E_{\al}$ be a family of conditional
expectations onto finite dimensional subalgebras $A_{\al}\subset
N$. Then, we define $d_\al=E_{\al}(d)$ and denote by $Q_\al$ the
corresponding quotient map. Let $\U$ be an ultrafilter refining
the natural order given by inclusion. We define
 \[ v:L_2(N,\tau)\oplus_1 L_2(N,\tau)\to \prodd_{\al,\U}
 K(A_{\al},d)  \]
by
 \[ v(x,y)(\al) \lel Q_{\al}(E_{\al}(x),E_{\al}(y)) \lel d_{\al}E_{\al}(x)+E_{\al}(y)d_{\al} \pl .\]
Since $x$ is  in $L_2(N,\tau)$, we have norm convergence of
$E_{\al}(x)$, $E_{\al}(y)$, $E_{\al}(d)$  to $x$, $y$, $d$
respectively. Thus, we get
 \[ \noo dx+yd\rrm_1 \lel \lim_{\al} \noo
 d_{\al}E_{\al}(x)+E_{\al}(y)d_{\al} \rrm_1 \pl .\]
This yields  $ker v \lel kerQ$ and $v$ induces a map
 \[ \hat{v}:K(N,d)\to \prodd_{\al,\U} K(A_{\al},d_{\al})\pl .\]
By construction it is clear that $v$ and $\hat{v}$ are complete
contractions. By density, it suffices to construct a map
$T:\bar{N}\to \to \prod_{\al} K(A_{\al},d_{\al})^*$ such that
 \[ \langle v(a,b), T(\gamma)\rangle \lel \gamma(Q(a,b)) \pl
 .\]
We define
 \[  T(n)_{\al}((a,b)+kerQ_{\al})\lel \tau(n^*(d_{\al}a+bd_{\al}))
 \]
Obviously, $\gamma_\al$ vanishes on $ker Q_{\al}$ and thus is a
linear functional on $K(A_\al,d_{\al})$. Given a matrix
$w_{kl}\in M_m(\bar{N})$, we deduce from \eqref{dualnorm} that the
functionals
 \[ \gamma_{\al}([a_{kl}],[b_{kl}]) \lel \summ_{kl}
  \tau(w_{kl}^*(d_{\al}a_{kl}+b_{kl}d_\al)) \quad
 \gamma([a_{kl}],[b_{kl}]) \lel \summ_{kl}
  \tau(w_{kl}^*(d_{\al}a_{kl}+b_{kl}d_\al))
  \]
satisfy
 \for
 \lefteqn{ \lim_{\al} \noo \gamma_{\al}
 \rrm_{M_m(K(A_\al,d_\al)^*)} }\\
 & & \lel  \lim_{\al} \max\left\{ \noo \summ_{kl} e_{kl} \ten
 d_\al w_{kl}\rrm_{M_m(L_2^r(N))},
 \noo \summ_{kl} e_{kl} \ten
 w_{kl}d_\al \rrm_{M_m(L_2^c(N))} \right\} \\
 & & \lel  \max\left\{ \noo \summ_{kl} e_{kl} \ten
 d w_{kl}\rrm_{M_m(L_2^r(N))},
 \noo \summ_{kl} e_{kl} \ten
 w_{kl}d\rrm_{M_m(L_2^c(N))} \right\} \\
 & & \lel   \noo \gamma \rrm_{M_m(K(N,d)^*)} \pl .
  \mel
In the last line, we use the norm  convergence of $d_{\al}$ to
$d$ in $L_2$. Hence, $T=(T_{\al}):\bar{N}\to \prod_{\al}
K(A_\al,d_\al)^*$ is a complete contraction satisfying
\[ \langle v(a,b),T(\gamma)\rangle \lel \lim_{\al}
 \tau(w^*(d_{\al}a+bd_{\al})) \lel \tau(w^*da)+\tau(bdw^*)
 \lel \gamma(a,b) \pl.\]
Since $\bar{N}$ is norm dense, we can extend $T$ to $K(N,d)^*$
and  the first assertion is proved. We will now assume that
$N=L_{\infty}(\Om,\mu;M_n)$ and $(\Om,\mu)$ is $\si$-finite. Let
$(A_k)$ be in increasing sequence of measurable set of finite
measure and $f_k=1_{A_k}$ the corresponding projections in the
center of $N$. Given  $d\in L_0(N)$, we may consider the spectral
projections of
$e_{l,k}=1_{[\frac1l,l]}(f_kdf_k)=1_{[\frac1l,l]}(d)f_k$ and
observe that $e_{l,k}d$ is in $N\cap L_2(N)$ and invertible  in
$e_{l,k}Ne_{l,k}$. Thus $d\in L_2^{ap}$. In order to show the
last assertion, we may choose a family sequence $\Si_{\al}$ of
subalgebras with finitely many atoms $\Om$ such that
$\bigcup_{\al} \Si_{\al}=\Si$. Then, we may use the conditional
expectations $E_{\al}:N\to L_\8(\Om,\Si_{\al};\mu;M_n)$ such that
$E_{\al}(x)$ converges to $x$ strongly. This provides the
particular form of $A_{\al}\cong \ell^{m_\al}(M_n)$ given by the
number of atoms $m_{\al}$ in $\Si_{\al}$.\qd

The space $G$ used to identify a copy of $OH$ involves two
densities. Therefore, given  two $\tau$-measurable positive
operators  $d_1$ and $d_2$ we define
 \[ K(N,d_1,d_2) \cong L_2^c(N,\tau)\oplus_1
 L_2^r(N,\tau)/kerQ_{d_1,d_2}  \]
where
 \[ Q_{d_1,d_2}(x,y)\lel d_1x+yd_2 \pl
 .\]
Again, we may think of the vector space $K(N,d_1,d_2)$ as a
subspace of $L_0(N,\tau)$. \lz

\begin{lemma}\label{2matrix} Let $d=\kla
\begin{array}{cc}d_1&0\\0&d_2\end{array}\mer \in
L_0(M_2(N),tr\ten \tau)$. Then $K(N,d_1,d_2)$ is completely
isometrically complemented in $K(M_2(N),d)$.
\end{lemma}\lz

\begin{proof}[\bf Proof:] Let $P:L_2(M_2(N),tr\ten \tau)\to L_2(N,\tau)$
be the projection onto the right upper corner, i.e.
  \[ P \kla \begin{array}{cc} x_{11}& x_{12}  \\x_{12}& x_{22} \end{array}\mer
   \lel x_{12} \pl .\]
Then
 \[ (P,P): L_2^c(M_2(N),tr\ten \tau)\oplus_1
 L_2^c(M_2(N),tr\ten \tau) \to L_2^c(N,\tau)\oplus_1
  L_2^r(N,\tau) \]
defined by $(P,P)(x,y)\lel (P(x),P(y))$ is a complete contraction
such that $Q_{d_1d_2}P$ vanished on $ker Q_d$. Here $Q_d$ denotes
the quotient map induced by $d$. Indeed,
 \[ d  \kla \begin{array}{cc} x_{11}& x_{12}  \\x_{12}& x_{22} \end{array}\mer
 +  \kla \begin{array}{cc} y_{11}& y_{12}  \\y_{12}& y_{22}
 \end{array}\mer d \lel 0\]
implies $d_1x_{12}+y_{12}d_2=0$. Let $\hat{P}:K(M_2(N),d)\to
K(N,d_1,d_2)$ be the induced contraction. Similarly, we have a
completely isometric  embedding $\iota:L_2^c(N,\tau)\oplus_1
L_2^r(N,\tau)\to L_2^c(M_2(N),tr\ten \tau)\oplus_1
L_2^r(M_2(N),tr\ten \tau)$ given by
 \[ \iota(x,y) \lel \kla
 \kla \begin{array}{cc} 0& x  \\ 0 & 0
 \end{array}\mer,
  \kla \begin{array}{cc} 0& y  \\ 0 & 0
  \end{array}\mer\mer  \pl .\]
Then $Q_{d}\p \iota$ vanishes on $kerQ_{d_1d_2}$. The induced
contraction $\hat{\iota}:K(N,d_1,d_2)\to K(N,d)$ satisfies
$\hat{P}\hat{\iota}=id$ and the assertion is proved.\qd

We should warn the reader that the definition of $K(N,d_1,d_2)$
involves a slightly different identification then the one used
for $G$. However, this can be adjusted by a change of density as
follows.\lz

\begin{lemma} \label{G} Let  $N=\ell_\8(L_\8([0,1],\mu))$ equipped  with
the normal, semifinite faithful trace $\tau(x)=\sum_{\nen} \intt
x_n(t) d\mu(t)$. Let $d_1(t,n)=\sqrt{t}$ and
$d_2(t,n)=\sqrt{1-t}$. Then $G$ is completely  isometrically
isomorphic to $K(N,d_1,d_2)$.
\end{lemma}\lz

\begin{proof}[\bf Proof:] Let $I_1:L_2(\nu_1;\ell_2)\to L_2(\mu;\ell_2)$
by the isometry $I_1(f)(t)=f(t)t^{-\frac12}$ and
$I_2:L_2(\nu_2;\ell_2)\to L_2(\mu;\ell_2)$ be given by
$I_2(f)=f(t)(1-t)^{-\frac12}$. Then, we note that
 \[ Q_{d_1d_2}(I_1(x),I_2(y)) \lel
  0 \]
if and  only if $x_n=y_n$ $\mu$ almost everywhere for all $\nen$.
Thus $Q_{d_1d_1}(I_1,I_2)$ induces a complete contraction
$\widehat{Q_{d_1d_1}(I_1,I_2)}:G\to K(N,d_1,d_2)$. The inverse is
given by the map $\widehat{Q(I_1^{-1},I_2^{-1})}$ induced by
$Q(I_1^{-1},I_2^{-1})$, where $Q:L_2^c(\nu_1;\ell_2)\oplus_1
L_2^r(\nu_1;\ell_2)\to G$ is the canonical quotient map.\qd

Voiculescu's inequality for sums of free independent random
variables  involves three terms. Therefore, we investigate these
modified spaces and the matrix-valued analogues. Given a positive
$\tau$-measurable operator $d\in L_2(N,\tau)$ and $t>0$, we define
 \[ \kz_t(N,d) \lel L_1(N,\tau)\oplus_1
 L_2^c(N,\tau)\oplus_1 L_2^r(N,\tau)/kerQ\]
on the  the vector space
 \[ V\lel L_1(N,\tau)+dL_2(N,\tau)+ L_2(N,\tau)d \subset L_1(N,\tau) \]
by the matrix norms
 \for
 \lefteqn{  \noo x\rrm_{S_1^m \wet \kz_t(N,d)} \lel }\\
  & & \!\! \!\!
  \inf_{x=x_1+(1\ten d) x_2+x_3(1\ten d)}
  t\noo x_1\rrm_{L_1(M_m\ten
  N,tr\ten \tau)}
    + \sqrt{t}\noo x_2\rrm_{S_1^m\wet
  L_2^c(N,\tau)} + \sqrt{t} \noo x_2\rrm_{S_1^m\wet
  L_2^r(N,\tau)}  .
  \mel
This simply defines a new operator space structure on the space
$L_1(N,\tau)$. \lz

\begin{lemma}\label{limit1} Let $d\in
L_2^{ap}(N,\tau)$. Let
 $x\in S_1^m\ten K(N,d)$, then
 \[ \liminf_{t\to \8} t^{-\frac12} \noo x\rrm_{S_1^m\wet
 \kz_t(N,d)} \lel  \pl \noo x\rrm_{S_1^m\wet K(N,d)} \pl .\]
Moreover, $K(N,d)$ is completely contractively complemented in
$\prod_{n,\U} \kz_{t_n}(N,d)$ where $(t_n)$ is a sequence
converging to $\8$ and $\U$ is a free ultrafilter.
\end{lemma}\lz

\begin{proof}[\bf Proof:] For every $(x_2,x_3)\in S_1^n\wet
L_2^c(N,\tau)\oplus_1 L_2^r(N,\tau)$, we have
 \for
 t^{-\frac12}\noo (1 \ten d)(x_2)+ x_3(1 \ten d)
 \rrm_{S_1^n\wet \kz_t(N,d)}
 & \le& \noo x_2\rrm_{S_1^n\wet L_2^c(N,\tau)}+
 \noo x_3\rrm_{S_1^n\wet L_2^r(N,\tau)} \pl .
 \mel
Thus the upper estimate is obvious and we may define $v:K(N,d)\to
\prod_{n,\U} \kz_{t_n}(N,d)$ by
$v((x,y)+kerQ)=(t_n^{-\frac12}(dx+yd))$. Clearly, $v$ is a
complete contraction. Let us assume that $(t_n)$ converges to
$\8$. We may assume by approximation that $d\in L_2(N,\tau)\cap
N$ and $d$ is invertible. In this case, we have shown in Lemma
\ref{wstar} that $\bar{N}$ is indeed norm dense in $K(N,d)^*$.
Given $w\in N$, we define $T(w)_n\in \kz_{t_n}(N,d)^*$ by
 \[ T(w)_n(x_1,x_2,x_3) \lel t_n^{\frac12}\tau(w^*(x_1+dx_2+x_3d))
 \pl .\]
Now, if $t_n\noo x_1\rrm_1 \le C $ for all $n\in A\in \U$, we have
 \[ \lim_{n,\U} |t_n^{\frac12} \tau(w^*x_1(n))| \kl C \lim_{n,\U}
 t_{n}^{-\frac12} \noo w\rrm_\8  \lel 0 \pl .\]
Thus, we have
 \[ \langle v(x_1,x_2)),T(w)\rangle \lel
 \tau(w^*(dx_1+x_2d)) \pl .\]
Given a matrix $[w_{ij}]$ such that $w_{ij}\in \bar{N}$, we
deduce from \eqref{dualnorm}that
 \for
  \lefteqn{\lim_{n,\U} \noo [T(w_{ij})_n]\rrm}  \\
  & & \lel  \lim_n \max\{  t_n^{-\frac12} \noo w_{ij}\rrm_{M_m(N)},
  \noo
  [dw_{ij}]\rrm_{M_m(L_2^r(N,\tau))},
  \noo  [w_{ij}d]\rrm_{M_m(L_2^c(N,\tau))}\}\\
  & & \lel \noo [w_{ij}]\rrm_{M_m(K(N,d)^*)} \pl .
 \mel
Thus $T$ is a complete isometry, which extends by density to
$K(N,d)^*$. Using the reflexivity of the Hilbert space $K(N,d)$,
we get $Tv=id$. This implies the lower estimate and the
factorization through the ultraproduct.\qd

\begin{rem} Given an arbitrary von Neumann algebra $\M$ and  $x\in \M_*\ten
K(N,d)$ then  still
 \[ \frac{1}{3} \noo x\rrm_{\M_*\wet K(N,d)} \kl  \liminf_{t\to \8} t^{-\frac12} \noo x\rrm_{\M_*\wet
 \kz_t(N,d)} \kl   \noo x\rrm_{\M_*\wet K(N,d)} \pl .\]
We will not use this result in the course of this paper.
\end{rem}\lz

\begin{proof}[\bf Proof:] According to Remark \ref{wstar2}, finite rank
tensors $y\in \M\ten \bar{N}\cap 3B$ are weak$^*$ in the unit
ball $B$ of $CB(\M_*,K(N,d)^*)$. Hence, given $x\in \M_*\ten
K(N,d)$, we deduce with the   Hahn-Banach theorem that there is a
net $y_\al$ satisfying $\noo y_\al\rrm\le 3$ such that
 \[ \noo x\rrm \lel \lim_{\al} |\langle
 x,y_\al\rangle | \pl .\]
Let $(t_n)$ tend to $\8$ and $y_\al=\sum_{i=1}^m a_i\ten w_i\in
\N\ten N$.  Then we may define a functional $y_n\in (\M_*\wet
\kz_n(N,d))^*$ by
 \for
 \lefteqn{  \langle (\phi_1\ten x_1,\phi_2\ten x_2,\phi_3\ten
 x_3), y_n\rangle}\\
 & &  \lel t_n^{\frac12} \summ_{i=1}^m
  (\phi_1(a_i)\tau(w_i^*x_1)+\phi_2(a_i)\tau(w_i^*dx_2)+\phi_3(a_i)\tau(w_i^*x_3d)
  \mel
The argument above shows that for a $x\in \M_*\wet K(N,d)$, we
have
 \for
 |\langle x,y_{\al}\rangle| &=& |\lim_n \langle
 t_n^{-\frac12} x,y_n\rangle|\\
 &\le& \liminf_n \noo t_n^{-\frac12}x\rrm_{\M_*\wet
 \kz_{t_n}(N,d)}\pl \limsup_n \noo y_n\rrm_{(\M_*\wet
 \kz_{t_n}(N,d))^*} \\
 &\le& \liminf_n \noo t_n^{-\frac12}x\rrm_{\M_*\wet
 \kz_{t_n}(N,d)} \pl \noo y_{\al}\rrm_{CB(\M_*,K(N,d)^*)}
 \\
 &\le& 3 \pl \liminf_n \noo t_n^{-\frac12}x\rrm_{\M_*\wet
 \kz_{t_n}(N,d)} \pl .
 \mel
Taking the supremum over $\al$ yields the assertion for finite
rank tensors $x$. By density it holds for all elements $x\in
\M_*\wet K(N,d)$.  \qd

At the end of this section, we will investigate $K$-functionals
associated to a conditional expectation instead of a state. Let
$\phi$ be a normal faithful state  on a von Neumann $\N$ algebra
and  $\E:\N\to N$ a normal  conditional expectation  onto a von
Neumann  subalgebra $N$ such that $\phi=\phi\circ E$. Let $D \in
L_1(\N)$ be the density of $\phi$. The space $L_1^c(\N,\E)$ is
defined as the completion of $\N D$ with respect to the norm
 \[ \noo xD\rrm_{L_1^c(\N,\E)} \lel \noo
 D\E(x^*x)D\rrm_{\frac12}^{\frac12} \pl .\]
Let us consider $\eps>0$ and $b=(D\E(x^*x)+\eps D)^{\frac14}\in
L_2(N)$, $y=xDb^{-1}$ and $z=xDb^{-2}$. Then $z$ is bounded
because $zz^*=xDb^{-4}Dx^*\le \eps^{-1}xx^*$ and thus $y=zb$ is
in $L_2(M)$ and satisfies
 \for
 \noo y\rrm_2^2 &=&  tr(b^{-1}Dx^*xDb) \lel
 tr(b^{-1}D\E(x^*x)Db^{-1}) \kl tr(b^2) \lel \noo b^2\rrm_2^2 \\
 &=& \noo D\E(x^*x)+\eps D\rrm_{\frac12}^{\frac12}
 \kl \noo D\E(x^*x)D\rrm_{\frac12}^{\frac12} +
 \eps^{\frac12} \pl .
 \mel
Thus, every element $xD$ can be written as a product  in
$L_2(M)L_2(N)$. To prove the converse, we recall the notations
$\N_a$, $N_a$ for the strongly dense subalgebras  of analytic
elements in $\N$, $N$, respectively. Recall that $x$ is called
analytic if the map $t\mapsto \si_t^{\phi}(x)$ extends to an
analytic function on $\cz$ with values in the underlyong von
Neumann algebra. Let $y=aD^{\frac12}bD^{\frac12}$ such that $a\in
\N$ and $b\in \N_a$. Then $y=a\si_{-\frac{i}{2}}(b)D$ and by
H\"{o}lder's inequality
 \for
  \noo y\rrm_{L_1^c(\N,\E)}^2
  &=& \noo D\E\big((a\si_{-\frac{i}{2}}(b))^*
  a\si_{-\frac{i}{2}}(b)\big)D\rrm_{\frac12}^{\frac12}
  \lel \noo D \si_{-\frac{i}{2}}(b)^*\E(a^*a) \si_{-\frac{i}{2}}(b)D\rrm_{\frac12}^{\frac12} \\
  &=& \noo D^{\frac12}b D^{\frac12} \E(a^*a)D^{\frac12} bD^{\frac12} \rrm_{\frac12}^{\frac12}
  \kl  \noo bD^{\frac12}\rrm_2^2
  tr(D^{\frac12}\E(a^*a)D^{\frac12}) \\
  &=& \noo bD^{\frac12}\rrm_2^2 \pl \noo
  aD^{\frac12}\rrm_2^2 \pl .
  \mel
Thus $L_1^c(\N,\E)=L_2(\N)L_2(N)$. This observation is  the key
ingredient in the proof of the following Lemma (see also
\cite{JD} and \cite{JX}).\lz

\begin{lemma}\label{bpE}
$L_1^c(\N,\E)$ is the quotient space of $L_2(M)\ten_{\pi}L_2(N)$
under the map $q(x\ten y)\lel xy$. The inclusion map
$L_1^c(\N,\E)\subset L_1(\N)$ is contractive and  $\N_aD$ is
dense in $L_1^c(\N,\E)$. The unit ball of the dual of
$L_1^c(\N,\E)$ is the strong operator closure of the unit ball
of  $\N p_\E$, where $p_\E:L_2(\N)\to L_2(\N)$ is the orthogonal
projection given by $p_\E(xD^{\frac12})=\E(x)D^{\frac12}$.
\end{lemma}\lz

In particular, we obtain a nice duality relation between
$L_1^c(\N,\E)$ and the space $L_\8^c(\N,\E)$ defined on $\N$ by
the norm
 \[ \noo x\rrm_{L_\8(\N,\E)} \lel
 \noo\E(x^*x)\rrm_\N^{\frac12} \pl .\]
Indeed, given $y=ba$, we have
 \for
 |tr(x^*y)| &=& |tr(x^*ba)| \lel  |tr(bax^*)|
 \kl \noo b\rrm_2  tr(ax^*xa^*)^{\frac12} \\
 &=& \noo b\rrm_2 tr(a\E(x^*x)a^*)^{\frac12}
 \kl \noo b\rrm_2 \noo a\rrm_2 \noo
 \E(x^*x)\rrm_{\8}^{\frac12} \pl .
 \mel
Hence, the complex conjugate space  $\overline{L_\8^c(\N,\E)}$
of  $L_\8^c(\N,\E)$ embeds naturally in $L_1^c(\N,\E)^*$. However,
for $x\in \N$, we have $\noo \E(x^*x)\rrm^{\frac12}=\noo
xp_{\E}\rrm$ and therefore Lemma \ref{bpE} implies
 \begin{eqnarray} \label{moddual}
 \noo
 y\rrm_{L_\8^c(\N,\E)}
 &=& \sup\big\{  |tr(y^*)| \pl \big | \pl     \noo x\rrm_{L_1^c(\N,\E)}\le
 1\big\} \pl .
 \end{eqnarray}

Similarly, we define for $p\in \{1,\8\}$ the corresponding row
version
 \[ \noo x\rrm_{L_p^r(\N,\E)} \lel \noo
 x^*\rrm_{L_p^c(\N,\E)} \pl .\]
The 3-term quotient space associated to $\kz_{n}(\N,\E)$ on
$L_1(\N)$ is given as follows
 \for
 \noo x\rrm_{\kz(\N,\E)} &=& \inf_{x=x_1+x_2+x_3}  \kla n\noo x_1\rrm_{L_1(\N)}
    \right. \\
 & & \hspace{2.7cm} \left.  +\sqrt{n } \noo \E(x_2^*x_2)^{\frac12}\rrm_{L_1(\N)}
  +\sqrt{n} \noo \E(x_3x_3^*)^{\frac12}\rrm_{L_1(\N)}
  \mer
  \pl .
 \mel
Here, we allow  $x_2\in L_1^c(\N,\E)\subset L_1(\N)$ and $x_3\in
L_1^r(\N,\E)\subset L_1(\N)$. We will also use the symbol
$\kz_{n}(\N,\E)$ for  $L_1(\N)$ equipped with this norm. Using
the anti-linear duality bracket
 \[ \lb y,x\rb_n \lel n tr(y^*x) \pl. \]
we will now calculate the dual space of $\kz_n(\N,\E)$.\lz

\begin{lemma}\label{dualk} The dual space of $\kz_n(\N,\E)$
with respect to the duality bracket $\lb \pl \pl ,\pl \rb_n$ is
$\bar{\N}$ and the norm is given by
 \for
  \noo y\rrm_{\kz_n(\N,\E)^*} &=&  \max\{ \noo y\rrm_{\N}, \sqrt{n} \noo
  \E(y^*y)\rrm_N^{\frac12},  \sqrt{n} \noo   \E(yy^*)\rrm_N^{\frac12} \}
  \pl .\
 \mel
\end{lemma}\lz

\begin{proof}[\bf Proof:] By continuity of the inclusions
$L_1^c(\N,\E)\subset L_1(\N)$ and $L_1^r(\N,\E)\subset L_1(\N)$
it is clear that $\kz_n(\N,\E)=L_1(\N)$ topologically. Thus the
dual space with respect to $\lb \pl ,\pl  \rb_n$ is $\bar{\N}$.
For an element $y\in \N$, we have
 \[ \sup_{n\noo x\rrm_1\le 1} |ntr(y^*x)| \lel  \noo y\rrm_\8 \pl
 .\]
According to \eqref{moddual}, we also have
 \[ \sup_{\sqrt{n} \noo x\rrm_{L_1^c(\N,\E)}\le 1} |ntr(y^*x)|
 \lel \sup_{\noo x\rrm_{L_1^c(\N,\E)}\le 1} \sqrt{n}
 |tr(y^*x)| \lel \sqrt{n} \noo \E(y^*y)\rrm_N^{\frac12}
 \pl .\]
The calculation for the row term is similar and thus the
assertion follows from the definition of $\kz_{n}(\N,\E)$ as a
quotient of three spaces. \qd

\begin{lemma}\label{ex1} Let $\N=M\ten N$,
$\phi_M$ be a normal faithful state on $M$ with density $d_M\in
L_1(M)$, $\phi_N$ be a normal faithful state on $N$ with density
$d_N\in L_1(N)$ and and $\phi=\phi_M\ten \phi_N$ with density
$D\in L_1(\N)$. The conditional expectation $\E$ is given by
 \[ \E(x\ten y)\lel \phi_N(y)x \]
Let $d$ be  the density of $\phi_N$ in $L_1(N)$. The densely
defined  map $T^{right}:L_1(M)\wet L_2^r(N)\to L_1(\N)$ defined by
 \[  T^{right}(md_M\ten nd_N^{\frac12}) \lel (m\ten n)D \]
extends to an isometry between $L_1(M)\wet L_2^r(N)$ and
$L_1^c(\N,\E)$. Similarly, $T^{left}(d_Mm\ten d_N^{\frac21}n)\lel
D(m\ten n)$ extends to an isometry between $L_1(M)\wet L_2^c(N)$
and $L_1^r(\N,\E)$.
\end{lemma}\lz

\begin{proof}[\bf Proof:]
By Kaplansky's density theorem $(N\ten M)D^{\frac12}$ is dense in
$L_2(\N)=L_2(M\ten N)$. Since $\N_{a}D$ is dense in
$L_1^c(\N,\E)$, we deduce that $L_2(\N)D^{\frac12}$ is dense in
$L_1^c(\N,\E)$ and thus $(N\ten M)D^{\frac12}$ is also dense in
$L_1^c(\N,\E)$. This shows that $T^{right}$ has dense range.
Moreover, given an element $z=\sum_{j=1}^k m_jD_M\ten
n_jd_N^{\frac12} \in L_1(\N)\ten L_2^r(N)$, we deduce from
\eqref{hc} and the isometric isomorphism
$i_{\frac12}:L_{\frac12}(M)\to L_{\frac12}(\N)$ given by
$\i_{\frac21}(d_Mmd_M)=DmD$ that
 \for
 \lefteqn{  \noo T^{right}(z)\rrm_{L_1^c(\N,\E)}^2 \lel \noo
  \summ_{j=1}^k   (n_j\ten m_j)D
  \rrm_{L_1^c(\N,\E)}^2 \lel
   \noo D\E(\summ_{j,l} n_l^*n_j\ten m_l^*m_j)D
   \rrm_{\frac12}^{\frac12}}
 \\
 &=&   \noo D\summ_{j,l} \phi(n_l^*n_j) m_l^*m_jD
 \rrm_{\frac12}^{\frac12} \lel
  \noo d_M\summ_{j,l} \phi(n_l^*n_j) m_l^*m_jd_M \rrm_{\frac12}^{\frac12}
  \\
  &=&
  \noo (z^*,z)\rrm_{L_\frac12(M)}^{\frac12}  \lel
 \noo z\rrm_{L_1(M)\wet L_2^r(N)} \pl .
  \mel
Hence, $T$ is isometric and has dense range in $L_1^c(\N,\E)$.
Therefore it extends to an isometric isomorphism. The proof of
the second assertion follows by passing to adjoints.\qd

As an application, we can identify the spaces $\kz_n(M\ten
N,\E)$.\lz

\begin{exam}\label{mainex} Let $N$ be a semifinite von Neumann algebra with trace $\tau$,
$d\in L_1(N,\tau)$ a positive element with full support and
$\phi_N(x)=\tau(dx)$. Let $M$ be a further $\si$-finite von
Neumann algebra with normal faithful state $\phi_M$ and
 \[ \E: M\ten N \to M  \quad, \quad \E(x\ten y) \lel
 \phi_N(y)\pl x \pl .\]
Then $\kz_n(M\ten N,\E)$ and $L_1(M)\wet \kz_n(N,d^{\frac12})$ are
isometrically isomorphic.
\end{exam}\lz

\begin{proof}[\bf Proof:] Since $d^{\frac12}\in L_2(N,\tau)$ we may consider
$\kz_n(N,d^{\frac12})$ as a subspace of $L_1(N,\tau)\subset
L_0(N,\tau)$. We will also need the family of $^*$-preserving
isomorphisms $(i_p)$ between $L_p(N,\tau)$ and $L_p(N)$. Indeed,
let $d_N\in L_1(N)$ be the of $\phi_N$ in $L_1(N)$. Since $d$ has
full support, we know that $\phi_N$ is a faithful normal state
and therefore $i_p$ is given by
$i_p(nd^{\frac1p})=nd_N^{\frac1p}$. Note that moreover
 \[ i_p(x_1x_2)\lel i_{p_1}(x_1)i_{p_2}(x_2) \]
holds for all $x_1\in L_{p_1}(N,\tau)$, $x_2\in L_{p_2}(N,\tau)$
such that $\frac1p=\frac{1}{p_1}+\frac{1}{p_2}$. In particular,
 \[ i_1(\kz_n(N,d^{\frac12})) \lel i_1(L_1(N,\tau))\lel L_1(N)\]
coincide as vector spaces and topologically. Let $d_M$ be the
density of $\phi_M$.  As we have seen in (1.9), we have a natural
isometric isomorphism $I_1$ between $L_1(M)\wet L_1(N)$ given by
$I_1(md_M\ten nd_N)\lel (m\ten n)D$, where $D\in L_1(M\ten N)$ is
the density of $\phi_M\ten \phi_N$ (see \cite{Fub} for details).
Therefore, it suffices to show that on $L_1(M\ten N)$  the
corresponding norms defined on their domain coincide.  This is
obvious for the $L_1$-norm. Given
 \[ z\lel \sum_j m_jd_M\ten n_jd_N \in L_1(M)\ten \kz_n(N,d)  \]
we deduce from Example \ref{mainex} that
 \begin{align*}
  \noo \summ_j m_jd_M\ten n_jd^{\frac12}\rrm_{L_1(M)\wet
 L_2^r(N,\tau)} &=
  \noo \summ_j m_jd_M\ten n_jd_N^{\frac12}\rrm_{L_1(M)\wet
 L_2^r(N)} \\
 &= \noo I_1(z)\rrm_{L_1^c(M\ten N,\E)} \pl .
 \end{align*}
By density of $(M\ten N)D$ in $L_1^c(M\ten N,\E)$, we deduce that
$I_1\circ(id\ten i_1)\circ Q$ provides an isometric isomorphism
between $L_1(M)\wet L_2^r(N,\tau)$ and $L_1^c(M\ten N,\E)$. Using
adjoints, we see that $L_1(M)\wet L_2^c(N,\tau)$ isometrically
isomorphic to $L_1^r(M\ten N,E)$. Finally, since $\wet$ preserves
quotient maps, we know that $L_1(M)\wet \kz_n(N,d^{\frac12})$ is
a quotient of
 \[ L_1(M)\wet L_2^c(N,\tau)\oplus_1 L_1(M)\wet L_2^r(N,\tau)\wet
 L_1(M)\wet L_1(N,\tau) \]
(with the corresponding weights) and hence $I_1\circ(id\ten
i_1)\circ Q$ provides an isometric isomorphism between
$L_1(M)\wet \kz_n(N,d^{\frac12})$ and the image $\kz_n(M\ten
N,\E)$.\qd

\section{Sums of free mean zero variables}

We will establish the probabilistic estimates using Voiculescu's
concept of free probability. This concept also plays a crucial
role in the proof of Grothendieck's theorem for operator spaces
in \cite{PS}. We extend Voiculescu's inequality for the norm  of
free mean $0$ variables. We are indept to U. Haagerup for the
collaboration on this extension (see Proposition \ref{Voi}).
Dualizing this result, we find a complemented copy of
$\kz_n(N,\E)$ in the predual of a von Neumann algebra. Let us
recall the notion of operator valued free probability needed in
this context. We assume that $(\M,\phi)$ is a
$W^*$-noncommutative probability space, meaning that $\M$ is a
von Neumann algebra and $\phi$ is a faithful normal state. Let
$B\subset \M$ be a von Neumann subalgebra and $E:\M\to B$ be a
normal conditional expectation satisfying $E\circ \phi=\phi$. A
family of subalgebras $(A_i)_{i\in I}$ is called {\it freely
independent over $E$} if $E(a_j)=0$ for all $j$   implies
 \[ E(a_1 \cdots a_n) \lel 0 \]
whenever $a_1\in A_{i_1}$,...,$a_{i_n}\in A_{i_n}$ and $i_1\neq
i_2\neq \cdots i_{n-1}\neq i_{n}$. In the original scalar valued
case one simply replaces $E$ by $\phi$ and $B=\cz1$. In the
following, we will use standard notation from free probability and
use ${\rm \AA}_{i}$ for the subspace of elements $a\in A_i$ with
$E(a)=0$. The following example is well-known. Since it is
crucial for this paper, we include a proof for the convenience of
the reader.\lz

\begin{exam}\label{matmod} Let $(M,\phi)$ be a probability space and $(A_i)$ freely independent over
$\phi$. Let $\M=N\ten M$ and $E(x\ten y)\lel \phi(y)x\ten 1$.
Then $(N\bar{\ten} A_i)_i$ are freely independent over $E$.
\end{exam}

\begin{proof}[\bf Proof:] Let $a=\summ_{j=1}^m x_j\ten y_j \in N\ten A_i$,
then we observe
 \[ a-E(a) \lel \summ_{j=1}^n x_j \ten
 (y_j-\phi(y_j)1) \in N\ten {\rm \AA}_{i} \pl .\]
Now, given $a_1,...,a_n$ such that $a_i\in N\ten {\rm \AA}_{i_k}$
with $i_1\neq i_2\cdots i_{n-1}\neq i_n$, we may write
$a_k=\sum_{j=1}^{m} x_{kj}\ten y_{kj}$ with $y_{kj}\in {\rm
\AA}_{i_k}$. (The same $m$ is achieved by adding $0$'s). Then, we
deduce from the freeness of the $A_i$'s that
 \for
 E(a_1\cdots a_n) &=& \summ_{j_1,...,j_n=1}^m
 x_{1k_1}\cdots x_{nk_n}
 \phi(y_{1,k_1}\cdots y_{n,k_n}) \lel 0 \pl .
 \mel
By Kaplansky's density theorem the unit ball of $N\ten_{min}A_i$
is strongly and strongly$^*$ dense in $N\bar{\ten}A_i$. Moreover,
if $a\in N\bar{\ten} A_i$ and $E(a)=0$ and $a^\al$ is a bounded
net converging in the strong topology to $a$, then
$a^\al-E(a^\al)$ converges to $a-E(a)=a$. Moreover, for bounded
nets $a^\al_1$, ... , $a^\al_n$, we still have strong convergence
of $ \lim_{\al_1}\cdots \lim_{\al_n}
(a_1^{\al_1}-E(a_1^{\al_1}))\cdots (a_n^{\al_n}-E(a_n^{\al_n}))$
to $a_1\cdots a_n$. By continuity of $E$ with respect to the weak
operator topology, we deduce $E(a_1\cdots a_n)=0$ for arbitrary
elements with $E(a_i)=0$, $a_k\in A_{i_k}$, $i_1\neq i_2\neq
\cdots \neq i_n$.\qd

Similar as in  the scalar case, operator valued freeness admits a
natural Fock space representation and this can be used to
construct the free amalgamated product. Let us review some very
basic and simple properties of this construction. We first assume
that $B,A_1,...,A_n$ are $C^*$-algebras such that $B\subset A_j$
for $j=1,...,n$ and  $E_j:A_j\to B$. Following \cite{Voir,Dyk2},
we consider the Hilbert $C^*$-module
 \[ {\rm \AA}_{i_1}\ten_B \cdots  \ten_B {\rm \AA}_{i_n} \]
with the $E$-valued scalar product
 \[ (a_1\ten \cdots \ten a_n,b_1\ten \cdots \ten b_n) \lel E_{i_n}(a_n^*E_{i_{n-1}}(a_{n-1}^*E_{i_{n-1}}(\cdots
 a_2^*E_{i_1}(a_1^*b_1)b_2\cdots )b_{n-1})b_n) \pl \]
which has an obvious multilinear extension. Then the usual Fock
space is replaced by the Hilbert $C^*$-module
 \[ \Hh_B \lel B \oplus \summ_{i_1\neq \cdots \neq i_n} \oplus \pll  {\rm \AA}_{i_1}\ten_B \cdots  \ten_B {\rm
 \AA}_{i_n}
  \pl .\]
Here $\oplus$ means that these components are all mutually
orthogonal. We recall that $\L(\Hh_B)$ is  the algebra of
adjoinable maps on $\Hh_B$.  A linear right $B$ module map
$T:\Hh_B\to \Hh_B$ is called {\it adjoinable} if there is a
linear map $S:\Hh_B\to \Hh_B$ such that
 \[ (x,Ty)\lel (Sx,y) \]
holds for all $x,y \in \Hh_B$. $\L(\Hh_B)$ is a C$^*$-algebra,
see \cite{La} for more details. Let us recall how elements in
$A_i$ act on $\Hh_B$. Every element $B\in B$ acts on $\Hh_B$ by
left multiplication on all the components (diagonal action).
However, for $c\in {\rm \AA}_i$, we have
 \[ c(a_{i_1}\cdots a_{i_n}) \lel \begin{cases} ca_{i_1}\cdots
 a_{i_n} & \mbox{ if $i\neq i_1$ } \\
  \big(ca_{i_1}-E(ca_{i_1})\big)a_{i_2}\cdots a_{i_n} \oplus
  E(ca_{i_1}) a_{i_2}\cdots a_{i_n}
  & \mbox{ if $i=i_1$ } \end{cases}  \pl .\]
The last definition also applies if $a_{i_2}\cdots a_{i_n}$ is
the empty word. Of course the action of $c$ is not diagonal. The
$C^*$-free product with amalgamation $C^*(\ast_B A_i)$ is defined
as the $C^*$-closure of linear combinations $a_{i_1}\cdots
a_{i_n}$. Following Voiculescu's work \cite{Voir} it is obvious
that whenever $(A_i)\subset \M$ are free over $E$, the
$C^*$-closure $\M_0$ generated by the $A_i$'s is isomorphic to
$C^*(\ast_B A_i)$ with respect to a $E$ preserving isomorphism.

Let us now explain the modifications needed for von Neumann
algebras. We will assume that the $B$ and $A_j$'s are von Neumann
algebras. Although this is not needed in general, we will also
assume that $\varphi:B\to \cz$ is a normal faithful state. Then
the induced states $\varphi_j\lel \varphi \circ E$ are faithful
and satisfy  $E\circ \si_t^{\varphi_j}\lel \si_t^{\varphi}\circ E$
for all $j=1,..,n$ and $t\in \rz$. Using $\varphi$, we can
associated to ${\rm \AA}_{i_1}\ten_B \cdots  \ten_B {\rm
 \AA}_{i_n}$ the Hilbert space $L_2({\rm \AA}_{i_1}\ten_B \cdots  \ten_B {\rm
 \AA}_{i_n},\varphi)$ given by the scalar product
 \[ \langle a_1\ten \cdots \ten a_n,b_1\ten \cdots \ten b_n\rangle \lel
 \varphi((a_1\ten \cdots \ten a_n,b_1\ten \cdots \ten b_n))
 \pl .\]
We obtain the Hilbert space
 \[ H_{\varphi}\lel \cz \oplus \summ_{i_1\neq \cdots \neq i_n}
  L_2({\rm \AA}_{i_1}\ten_B \cdots  \ten_B {\rm
 \AA}_{i_n},\varphi)  \pl .\]
Since $\varphi$ is faithful, we see that the representation of
$\L(\Hh_B)$ on $H_{\varphi}$ is faithful. Let us denote by
$\M(\Hh_B)$ the corresponding von Neumann algebra. Then the von
Neumann algebra $\ast_B A_i$ is defined as the closure of
$C^*(\ast_B A_i)$ with respect to the weak operator topology. We
use the notation  $Q_{\emptyset}$, $Q_{i_1\cdots i_n}$ for the
orthogonal projections onto $L_2(B,\varphi)$, $L_2({\rm
\AA}_{i_1}\ten_B \cdots  \ten_B {\rm \AA}_{i_n},\varphi)$,
respectively.  Then the conditional expectation $E$ onto $B$  is
given by
 \[ E(x)\lel Q_{\emptyset}xQ_{\emptyset} \in B \subset
 B(L_2(B,\phi)) \pl .\]
Given an  element $x\in \ast_B A_i$, we may consider the elements
in the first column
 \[ x_{i_1\cdots i_n}\lel Q_{i_1\cdots i_n}xQ_{\emptyset} \pl .\]
We note that $E(x^*x)=Q_{\emptyset}x^*xQ_{\emptyset}=0$ implies
$xQ_{\emptyset}=0$ and hence $x_{i_1\cdots i_n}=0$. In order to
show that $E$ is faithful on $\ast_BA_i$, we have to show that
$x=0$. Let us first consider a tuple $(j_1,....j_m)$  with
$j_1\neq i_n$. By approximation with elements in $C^*(\ast_B
A_i)$, we deduce that for $h\in Q_{j_1\cdots j_m}$
 \[ Q_{i_1\cdots i_m}x(h) \lel x_{i_1\cdots i_n}\ten h \pl .\]
Thus $E(x^*x)=0$ implies that entry is $0$. Now, we assume
$j_1=i_n$ and that there are no other coefficients. Let us
consider $Q_{A_j}\lel Q_{\emptyset}+Q_j$. Again by approximation
with elements in $C^*(\ast_B A_i)$ we see that
$Q_{A_j}x^*xQ_{A_j}$ is a positive element in $A_j$ such that
 \[ E(x^*x)\gl E(Q_{A_j}x^*xQ_{A_j}) \pl .\]
Since $E$ is faithful, we deduce $Q_{i_1\cdots i_m}xQ_{A_j}=0$.
Since this holds for all $i_1\cdots i_m$, we get $xQ_{a_j}=0$.
Given an arbitrary element $h=a_{j_1}\cdots a_{j_m}$ with
$j_1=i_n$, we may use the right action and find
 \[ x(h) \lel R_{a_{j_m}}\cdots R_{a_{j_2}} x(a_{j_1}) \lel
 R_{a_{j_m}}\cdots R_{a_{j_2}} xQ_{A_j}(a_{j_1}) \lel 0 \]
provided $E(x^*x)=0$. This shows that $E$ is faithful on
$\ast_BA_i$. In this argument we use the  conditional expectation
$E_j:\ast_B A_i\to A_j$ given by
 \begin{equation}\label{condex} E_j(x)\lel Q_{A_j}xQ_{A_j} \pl .
 \end{equation}
Clearly, $E_j$ defines a completely positive  map on
$B(H_{\phi})$. Given $a_1\cdots a_n\in {\rm \AA}_{i_1}\ten_B
\cdots  \ten_B {\rm \AA}_{i_n}$, we see that $E_j(x)=0$ if and
only if $n=1$ and $i_1=i_n=j$. By Kaplansky's density theorem we
deduce $E_j(\ast_BA_i)\lel A_j$. We will use this opportunity to
clarify the action of the modular group of $\phi\lel \varphi\circ
E$ which is very similar to the classical scalar case. Indeed, for
every $H_{i_1\cdots i_n}$ we deduce from
  \begin{align*}
 \lefteqn{ \si_t^{\varphi}(a_1\ten \cdots \ten a_n,b_1\ten \cdots \ten
  b_n) }\\
 &=  \si_t^{\varphi}(E(a_n^*E(a_{n-1}^*E(\cdots
 a_2^*E(a_1^*b_1)b_2\cdots b_{n-1})b_n))\\
 &= E(\si_t^{\varphi}(a_n^*E(a_{n-1}^*E(\cdots
 a_2^*E(a_1^*b_1)b_2\cdots b_{n-1})b_n)))) \\
 &=
 E(\si_t^{\varphi}(a_n^*)E(\si_t^{\varphi}(a_{n-1}^*)E(\cdots
 E(\si_t^{\varphi}(a_1^*)\si_t^{\varphi}(b_1))\cdots \si_t^{\varphi}(b_{n-1}))\si_t^{\varphi}(b_n)))\\
 &= (\si_t^{\varphi}(a_1)\ten \cdots \ten \si_t^{\varphi}(a_n),
 \si_t^{\varphi}(b_1)\ten \cdots \ten \si_t^{\varphi}(b_n)) \pl
 \end{align*}
that $u_t^{i_1,...,i_n}(a_1\ten \cdots \ten a_n)=
\si_t^{\varphi}(a_1)\ten \cdots \ten \si_t^{\varphi}(a_n)$
extends to a unitary on $H_{i_1\cdots i_n}$. Therefore
$u_t=\sum_{i_1\neq \cdots \neq i_n}u_t^{i_1,...,i_n}$ is a unitary
on $H_{\phi}$ and it easily checked that the family $(u_t)_{t\in
\rz}$ is strongly continuous. Moreover, for $a\in A_j$, we find
 \[ u_tau_t^*\lel \si_t^{\varphi}(a) \pl .\]
In particular $u_tC^*(\ast_BA_i)u_t^*\subset C^*(\ast_BA_i)$ and
hence $u_t(\ast_BA_i)u_t^*\lel \ast_B A_i$. Let $a_1\in {\rm
\AA}_{i_1}$, $a_2\in {\rm \AA}_{i_2}$,..., $a_n\in {\rm
\AA}_{i_n}$ be analytic elements and $x=b_1\ten \cdots \ten b_m
\in {\rm \AA}_{i_1}\cdots {\rm \AA}_{i_n}$.
 \[ f(z)\lel \phi(x^*\si_{-iz}^{\varphi}(a_1)\ten \cdots \ten
 \si_{-iz}^{\varphi}(a_n)x) \]
is analytic. Obviously
 \[ f(it)\lel \phi(x^*\si_t^{\varphi}(a_1)\ten \cdots
 \ten \si_t^{\varphi}(a_n))\lel \phi(x^*u_t(a_1\cdots a_n)u_t^*) \pl .\]
Moreover, since $E$ is selfadjoint on $L_2(A_j)$ we find
 \begin{align*}
 f(1+it) &= \phi(b_{i_n}^*\cdots b_{i_1}^*\si_{-i+t}^{\varphi}(a_1)\cdots
 \si_{-i+t}^{\varphi}(a_n))\\
 &= \phi(b_n^*E(b_{n-1}^*\cdots
 E(b_1^*\si_{-i+t}^{\varphi}(a_1))\cdots)\si_{-i+t}^{\varphi}(a_n)) \\
 &= \phi(\si_{t}^{\varphi}(a_n)b_n^*E(b_{n-1}^*\cdots
 E(b_1^*\si_{-i+t}^{\varphi}(a_1))\cdots)) \\
 &= \phi\bigg(E\big(\si_t^{\varphi}(a_n)b_n^*\big)b_{n-1}^*E(b_{n-2}^*\cdots E(b_1^*\si_{i+t}^{\varphi}(a_1))
 \cdots si_{-i+t}^{\varphi}(a_{n-1}) \bigg)\\
 &= \cdots \\
 &= \phi(\si_t^{\varphi}(a_1)E(\si_t^{\varphi}(a_2)\cdots      E(\si_t(a_{n-1})E(\si_t(a_n)b_n^*)b_{n-1}^*)\cdots
 b_1^*) \\
 &=  \phi(u_t a_1\cdots a_nu_t^*
 b_{i_n}^*\cdots b_{i_1}^*)  \pl .
 \end{align*}
By approximation, we obtain a two point function for every
polynomial with analytic coefficients  in $C^*(\ast_BA_i)$. Using
Kaplansky's density theorem, we find a two point function
obtained for arbitrary elements $x,y\in \ast_BA_i$ and thus by
modular theory (see e.g. \cite{Strat})
 \begin{equation} \label{modg}  \si_t^{\phi}(a_1\cdots a_n)
  \lel \si_t^{\varphi_{i_1}}(a_1)\cdots \si_t^{\varphi_{i_n}}(a_n)
\end{equation}
for  all $a_1\in A_{i_1}$,...,$a_n\in A_{i_n}$.

For our estimates, it will be crucial to follow Voiculescu
\cite{Voir} in defining the projections
 \[ P_i\lel \sum_{i=i_1\neq \cdots i_n}Q_{i_1\cdots i_n}\pl .\]
Note that $P_i$ is in general not in $\ast_B A_i$.\lz

\begin{lemma} \label{fre1} Let $c\in A_i$ such that $E(c)\lel 0$, then
\[ (1-P_i)c(1-P_i)\lel 0 \pl .\]
\end{lemma}

\begin{proof}[\bf Proof:] Given $a_{i_1}\cdots a_{i_n}\in
{\rm \AA}_{i_1}\ten_B \cdots  \ten_B {\rm
 \AA}_{i_n}$ and $i_1\neq i$, we observe that
 $c(a_{i_1}\cdots a_{i_n})=ca_{i_1}\cdots a_{i_n}$ is an element
 of ${\rm \AA}_i\ten_B{\rm \AA}_{i_1}\ten_B \cdots  \ten_B {\rm
 \AA}_{i_n}$ and thus $P_i(ca_{i_1}\cdots a_{i_n})=ca_{i_1}\cdots
 a_{i_n}$. By linearity this yields the assertion.\qd

\begin{cor} \label{fre2} For all $a\in A_i$,
\[ (1-P_i)a (1-P_i)\lel E(a)(1-P_i) \pl.\]
\end{cor}

\begin{proof}[\bf Proof:] This is obvious from Lemma \ref{fre1} by writing
$a=a-E(a)+E(a)$.\qd

The operator valued extension of Voiculescu's argument reads as
follows.\lz

\begin{prop} \label{Voi} Let {\rm $a_k\in {\rm \AA}_{i_k}$} such that the
$i_k,...,i_n$ are mutually different. Then
\[ \noo \summ_{k=1}^n a_{i_k}\rrm \kl \sup_{k=1,..,n} \noo a_{i_k}\rrm +
\noo \summ_{k=1}^n E(a_{i_k}^*a_{i_k}) \rrm^{\frac 12}+ \noo
\summ_{k=1}^n E(a_{i_k}a_{i_k}^*)\rrm^{\frac 12} \pl.
\]
\end{prop}

\begin{proof}[\bf Proof:] Let us use $P_k=P_{i_k}$. We deduce from Lemma
\ref{fre1} that
 \for
 \summ_{k=1}^n a_{k} &=& \summ_{k=1}^n P_ka_{k}P_k + \summ_{k=1}^n
 (1-P_k)a_{k}P_k + \summ_{k=1}^n
 P_ka_k(1-P_k) \\
 & & \pll \pll + \summ_{k=1}^n
 (1-P_k)a_k(1-P_k) \\
 &=& \summ_{k=1}^n P_ka_kP_k + \summ_{k=1}^n
 (1-P_k)a_kP_k + \summ_{k=1}^n
 P_ka_k(1-P_k)  \pl .
 \mel
Since $P_k$'s  are mutually orthogonal, we have
 \[ \noo \summ_{k=1}^n P_ka_kP_k \rrm \lel \sup_k \noo
 P_ka_kP_k\rrm \kl \sup_k \noo a_k\rrm \pl .\]
Now, we  consider the second term. By orthogonality, positivity
and the module property, we deduce from  Corollary \ref{fre2} that
 \for
 \lefteqn{ \noo \summ_{k=1}^n
 (1-P_k)a_kP_k \rrm^2 \lel  \noo \summ_{k,l}
 (1-P_k)a_kP_kP_{l}a_{l}^*(1-P_{i_l}) \rrm }\\
 & & \lel  \noo \summ_{k=1}^n (1-P_k)a_kP_ka_k^*(1-P_k) \rrm
 \kl  \noo \summ_{k=1}^n (1-P_k)a_ka_k^*(1-P_k) \rrm\\
 & & \lel  \noo \summ_{k=1}^n (1-P_k)E(a_ka_k^*)(1-P_k) \rrm
   \\
% \mel\for
 & & \lel \noo \summ_{k=1}^n E(a_ka_k^*)^{\frac12}(1-P_k)(1-P_k)E(a_ka_k^*)^{\frac12}
 \rrm\\
 & & \kl  \noo \summ_{k=1}^n
 E(a_ka_k^*)^{\frac12}E(a_ka_k^*)^{\frac12}  \rrm
 \lel  \noo \summ_{k=1}^n E(a_ka_k^*)  \rrm \pl .
 \mel
The calculation for the third term is very similar.\qd

We will now prove a converse of Proposition \ref{Voi}.\lz

\begin{lemma} \label{conv} Let $i_1,...,i_n$ be mutullaly
different and {\rm $a_1\in {\rm \AA}_{i_1}$,... $a_n\in {\rm
\AA}_{i_n}$}. Let $D_{\phi}^\M$, $D_{\varphi_{i_k}}^{A_{i_k}}$ and
$D_{\varphi}^B$ be the density of $\phi$, $\varphi$  in $L_1(\M)$,
$L_1(A_{i_k})$  and $L_1(B)$, respectively. Then
 \begin{eqnarray}
  \noo \summ_{k=1}^n a_kD_{\phi}^\M \rrm_{L_1(\M)} &\le& \summ_{k=1}^n \noo
  a_kD_{\varphi_{i_k}}^{A_{i_k}}\rrm_{L_1(A_{i_k},\phi)}\pl , \label{s1}
   \\
  \noo \summ_{k=1}^n a_kD_{\phi}^\M \rrm_{L_1(\M)} &\le& \noo
  (\summ_{k=1}^n  D_{\varphi}^B E(a_k^*a_k)D_{\varphi}^B)^{\frac12}
    \rrm_{L_1(B)}\pl , \label{s2}
   \\
     \noo \summ_{k=1}^n D_{\phi}^\M a_k\rrm_{L_1(\M)} &\le& \noo
  (\summ_{k=1}^n  D_{\phi}^B E(a_ka_k^*)D_{\phi}^B)^{\frac12}
    \rrm_{L_1(B)}\pl . \label{s3}
  \end{eqnarray}
\end{lemma}

\begin{proof}[\bf Proof:] Since $\phi\circ E=\phi$, we have a natural
family of embeddings $\iota_{p}:L_p(B)\to L_p(\M)$ satisfying
 \begin{equation}\label{theta} \iota_p([D_{\phi}^{B}]^{\frac{1-\theta}{p}}b[D_{\phi}^{B}]^{\frac{\theta}{p}})
 \lel
 [D_{\phi}^{\M}]^{\frac{1-\theta}{p}}b[D_{\phi}^{\M}]^{\frac{\theta}{p}}
 \end{equation}
for all $b\in B$ (see \cite{JX}). Moreover, we  have an
$E$-preserving  conditional expectation $E_j:\M\to A_j$ (see
\eqref{condex}). Therefore, we have a family of isometric
embeddings $\iota_p:L_p(A_j)\to L_p(\M)$.
$(\iota_p)_{0<p<\infty}$ satisfying
 \begin{equation}\label{theta2} \iota_p([D_{\varphi_j}^{A_j}]^{\frac{1-\theta}{p}}a[D_{\varphi_j}^{A_j}]^{\frac{\theta}{p}})
 \lel
 [D_{\phi}^{\M}]^{\frac{1-\theta}{p}}a[D_{\phi}^{\M}]^{\frac{\theta}{p}}
 \end{equation}
for all $a\in A_j$. The triangle inequality implies \eqref{s1}.
For the proof of \eqref{s2}, we consider $z=\sum_{k=1}^n
a_kD_{\phi}^\M$ with $a_k\in {\rm \AA}_{i_k}$. Then, we have
 \[ E(z^*z) \lel D_{\phi}^\M \summ_{j,k} E(a_j^*a_k) D_{\phi}^\M
 \lel D_{\phi}^\M \summ_{k} E(a_k^*a_k) D_{\phi}^\M \pl .\]
Thus by Lemma \ref{bpE} and \eqref{theta}, we deduce
 \begin{align*}
 \noo \summ_{k=1}^n a_k D_{\phi}^\M\rrm_1 &\le  \noo
 E(z^*z)^{\frac12}
  \rrm_1 \lel \noo \summ_{k=1}^n D_{\phi}^\M
 E(a_k^*a_k)D_{\phi}^\M \rrm_{\frac12} \\
 &= \noo \summ_{k=1}^n D_{\phi}^B
 E(a_k^*a_k)D_{\phi}^B \rrm_{\frac12}
  \pl .
  \end{align*}
Assertion \eqref{s3}  follows by taking adjoints, an isometry on
$L_1(\M)$.\qd

Given a von Neumann algebra $\N$ and a conditional expectation
$\E:\N\to B$, we will now construct an embedding of
$\kz_n(\N,\E)$ in the predual of an appropriate amalgamated free
product. Let us assume that $\varphi_\N$ is  a normal faithful
state on $\N$  satisfying $\varphi_\N\circ \E=\varphi_\N$. In
order to achieve mean $0$ elements, we will work with the algebras
$\B=\ell_\8^2(\N)$ and the state $\varphi_\B(x,y)=\frac 12
(\varphi_\N(x)+\varphi_\N(y))$. Let us denote by $E_2$ the
conditional expectation
$E_2(x,y)=(\frac12(\E(x)+\E(y)),\frac12(\E(x)+\E(y)))$. Clearly,
$\varphi_\B \circ E_2\lel \varphi_\B$ and $E_2$ is a conditional
expectation onto (a copy of) $B$. Given $x\in \N$, we can now
define the symmetrization $\eps x\lel (x,-x)\in \B$ which
obviously satisfies $E_2(\eps x)=0$. We will drop the indices
$_2$,  $_\N$,  and  $_\B$ in the following, because there are
immediate from their context. Given these data, we may now define
$A_i\lel \B$ for $i=1,..,n$ and $\M_n=\ast_BA_i$ the amalgamated
free product over $B$. \lz

\begin{prop} \label{comp1} Let $\M_n$, $A_1,...,A_n$, $\B$, $B$ and $E$ be as
above. The space $\kz_{n}(\N,\E)$ is $3$-complemented in
$L_1(\M)$.
\end{prop}

\begin{proof}[\bf Proof:] Let $D_{\phi}^{\M}$, $D_{\varphi}^B$,
$D_{\varphi}^\N$, $D_{\varphi}^\B=
(\frac12D_{\varphi}^\N,\frac12D_{\varphi}^\N)$ be the densities of
$\phi$, $\varphi$  in $L_1(\M)$, $L_1(B)$, $L_1(\N)$,
$L_1(\B)=L_1(\N)\oplus L_1(\N)$, respectively. We define the map
$T:L_1(\N)\to L_1(\M)$ by
 \[ T(xD_{\varphi}^\N) \lel \summ_{k=1}^n \pi_k(\eps x)D_{\phi}^{\M}
 \pl .\]
According to Lemma \ref{conv} (see \eqref{s1}), we have
 \[ \noo T(xD_{\varphi}^\N)\rrm_1 \kl  n \noo (\eps x) D_{\varphi}^\B\rrm_{L_1(\B)}
 \lel n \noo xD_{\varphi}^\N\rrm_{L_1(\N)}
  \pl .\]
Similarly, we deduce from the fact that $\eps x$ is mean $0$ and
Lemma \ref{conv} (see \ref{s2})  that
 \for
 \noo T(xD_{\varphi}^\N)\rrm &\le  \noo D_{\varphi}^B\summ_{k=1}^n
 E(x^*x) D_{\varphi}^B \rrm_{L_\frac12(B)}^{\frac12}
 \lel      \sqrt{n}  \noo (D_{\varphi}^B E(x^*x)D_{\varphi}^B)^{\frac12}
  \rrm_{L_1(B)}
  \pl .
   \mel
For an analytic element $x\in \N_a$, we deduce from the action
$\si_t^{\phi}$ on the free product (see \eqref{modg}) that
 \begin{eqnarray}\label{rtrans}
 T(D_{\varphi}^\N x) &=& T(\si_{-i}^{\varphi}(x)D_{\varphi}^\N)
 \lel \summ_{k=1}^n
 \si_{-i}^{\phi}(\pi_k(\eps x))D_{\phi}^\M \lel
  \summ_{k=1}^n D_{\phi}^\M \pi_k(\eps x) \pl .
 \end{eqnarray}
By continuity this extends to all elements $x\in \N$. Hence,
\eqref{s3} implies the missing inequality and we deduce
 \[ \noo T:\kz_{n}(\N,\E) \to L_1(\M)\rrm \kl 1 \pl .\]
We define the map $S:\bar{\N}\to \bar{\M}$ by
 \[ S(y) \lel  \summ_{k=1}^n \pi_k(\eps y) \pl .\]
We refer to Lemma \ref{dualk} for calculating  the dual space
$\kz_{n}(\N,\E)^*$ (formed with the anitilinear duality).
According to Proposition \ref{Voi}, we have $\noo
S:(\bar{\N},\kz_{n}(\N,\E)^*)\to \bar{\M}\rrm\le 1$. Moreover,
  \for
  \lefteqn{  \lb S(y),T(xD_{\varphi}^{\N})\rb_{(\bar{\M},L_1(\M))}
  \lel    tr_\M(S(y)^*T(xD_{\varphi}^\N)) }\\
   & &  \lel \summ_{k=1}^n tr_\M(\pi_k(\eps   y)^*\pi_k(\eps x)D_{\phi}^\M)
   \lel \summ_{k=1}^n \varphi(y^*x) \lel
 \lb y, x\rb_n \pl .
   \mel
This shows $S^*T$ coincides with the natural embedding of
$\kz_n(\N,\E)$ into its bidual and hence $T$ and $S$ are
isomorphisms. (Note that $S^*$ is linear!) Now, we show that
\[ S^*(L_1(\M))\subset \kz_n(\N,d) \pl .\]
Indeed, given an element $z=a_1\cdots a_nD_{\phi}^\M$ with
$a_1\cdots a_n\in {\rm \AA}_{i_1}\ten_B\cdots \ten_B {\rm
\AA}_{i_n}$.  We observe that for $S^*(z)\neq 0$ to hold we must
have $n=1$. If $z=\pi_k((a_1,a_2))D_{\phi}^\M$ is given by
$a_1,a_2 \in \N$, then
 \for
  \lb S(y),z\rb   &=&tr((\pi_k(\eps
 y))^* \pi_k((a_1,a_2))D_{\phi}^\M) \lel
  \frac{1}{2}
 (\varphi(y^*a_1)-\varphi(y^*a_2))\\
 &=&  \frac12 tr_\N(y^*(a_1-a_2)D_{\varphi}^\N)
   \pl .
 \mel
This implies $S^*((a_1,a_2))D_{\phi}^\M=\frac12
(a_1-a_2)D_{\varphi}^\N$. By continuity and density, we deduce
the inclusion  $S^*(L_1(\M))\subset L_1(\N)$. Hence,
$S^*T=id_{\kz_n(\N,\E)}$ and the assertion is proved.\qd

\begin{theorem}\label{comp2} Let $N$ be a von Neumann algebra with normal faithful trace $\tau$,
$d\in L_1(N,\tau)$ a positive element with full support. Let
$\varphi_d(x)=\tau(xd)$. On  $\ell_\infty^2(N)$ the normal
faithful state is given by
$\varphi_2(x_1,x_2)=\frac{1}{2}\varphi(x_1)+\frac12\varphi(x_2)$.
Then $\kz_{n}(N,d^{\frac12})$ is $3$-completely complemented in
$L_1(\ast_{i=1}^n (\ell_\infty^2(N),\varphi_2))$.
\end{theorem} \lz

\begin{proof}[\bf Proof:]
Let us refer to the proof of  Example \ref{mainex} on how we may
identify  $\kz_n(N,d)$ with the image $i_1(\kz_n(\N,\sqrt{d}))$ in
Haagerup space  $L_1(N)$.  For fixed $m\in \nz$, we consider
$\N=M_m(N)$ and the conditional expectation $\E(x\ten
y)=\phi_d(y)x$ onto $B=M_m$. On $B$ we use the normal faithful
state $\varphi_B=\frac{tr}{m}$. We recall that
$\B=\ell_\infty^2(\N)=M_m(\ell_\infty^2(N))$. For $i=1,..,n$ we
have $A_i\lel \B$. We deduce from Lemma \ref{matmod} that the
amalgamated product $\M_n$ satisfies
 \[ \M_n \lel \ast_{M_m} A_i\lel M_m(\ast_{i=1}^n
  (\ell_\infty^2(N),\varphi_2))
 \pl .\]
Let us define  $\M_n^1 \lel \ast_{i=1}^n
(\ell_\infty^2(N),\varphi_2)$ and $\phi^1_n$ for the faithful
normal state $\phi^1_n\lel \ast_{i=1}^n \varphi_2$. Using the fact
that $\frac{tr}{m}$ is a normalized trace it is clear that
 \[ L_1(\M_n)\lel L_1(M_m,\frac{tr}{m})\wet L_1(\M_n^1) \]
and the density of $D_{\phi}^{\M_n}$ is given by
$id_{\ell_2^m}\ten D_{\phi^1_n}^{\M_n^1}$. Then, we may apply
Proposition \ref{comp1} for $m=1$ and deduce that
  \[ T(xD_{\phi_d}^N) \lel \summ_{k=1}^n \pi_k(\eps
 x)D_{\phi_n}^{\M_n^1} \]
is a contraction. However, using Proposition \ref{comp1} for $m$,
we deduce that also
 \begin{align*}
  (id_{L_1(M_m)}\ten T)\big([x_{ij}](id_{\ell_2^m} \ten D_{\phi_d}^N)\big) &=
  \summ_{k=1}^n (id_{M_m}\ten \pi_k)([\eps
  x_{ij}])(id_{\ell_2^m}\ten D_{\phi^1_n}^{\M_n^1}) \\
 &= \summ_{k=1}^n (id_{M_m}\ten \pi_k)([\eps
  x_{ij}])D_{\phi}^{\M_n}
  \lel  T_m([x_{ij}]D_{\varphi}^{\N})
  \end{align*}
is a contraction. Hence $T$ is a complete contraction. Similarly,
using the concrete form
 \[ S_m((a\ten b_1,a\ten b_2)D_{\phi}^\M)
 \lel (a \ten \frac{b_1-b_2}{2}) D_{\varphi}^\N
 \lel a\ten (\frac{b_1-b_2}{2})D_{\phi_d}^N \]
we obtain  $S_m=id \ten S$. Therefore $S$ has cb-norm less than
$3$ and still $ST=id_{K(N,d^{\frac12})}$.\qd

Before, we proof our main results we need some results on the
QWEP property.\lz

\begin{lemma}\label{ultraqwep}  Let $(M^s)$ be  a family of QWEP von
Neumann algebras, then the von Neumann algebra $\M=(\prod_{s,\U}
M^s_*)^*$ also has QWEP.
\end{lemma} \lz

\begin{proof}[\bf Proof:] According to Kirchberg \cite{Ki}, we know that
$\prod M^s$ is  QWEP and thus \cite{Ki}, $(\prod M^s)^{**}$ is
QWEP. Following Groh \cite{Gr}, we observe that  the space of
functionals $\prod_{s,\U} M^s_*$ on $\prod M^s$ is left and right
invariant under the action of $\prod M^s$ and hence there is a
central projection $z_\U$ such that $\M\cong z_\U(\prod
M^s)^{**}$. Equivalently, $\M=(\prod M_n)^{**}/I_\U$, where
$I_\U=z_\U(\prod M^s)^{**}$. Thus $\M$ is QWEP. \qd

The next lemma follows an argument of Pisier/Shlyaktenko
\cite[Lemma 2.5]{PS} and is based results of Dykema \cite{Dyk}.
\lz

\begin{lemma}\label{dyk}  Let $A$ be a finite dimensional
algebra and $\varphi$ be a state on $A$, then the free product
$\ast_{i=1}^n (A,\psi)$ is QWEP.
\end{lemma}

\begin{proof}[\bf Proof:]
In \cite{Ki}   Kirchberg proved that $L(\ff_n)=\ast_{i=1}^n
L^\infty$ has QWEP. In particular, if $\varphi$ is a trace, then
we can find a trace preserving embedding of $A$ into $L(\ff_2)$.
Therefore $\ast_{i=1}^n (A,\varphi)$ is the range of a conditional
expectation $\E:L(\ff_{2n})\to \ast_{i=1}^n (A,\varphi)$. Let us
assume that $\varphi(x)=tr(xd)$ such that $\ln d$ is a duaydic
rational. Then we can apply Dykema's results \cite[Theorem 1 and
Theorem 4]{Dyk} to deduce that $\ast_{i=1}^n
(A_i,\varphi)=M_0\oplus \M_1$ where $\M_0$ is finite dimensional
or $\{0\}$  and $\M_1$ is a type $III_{\la}$ factor with
centralizer $L(\ff_\8)$ and the induced state is almost periodic.
In particular, $\M_1$ is a factor. According to Connes result
\cite[Theorem 4.4.1]{C2}, $\M_1$ is a isomorphic to a crossed
product of $L(F_{\infty})\rtimes G$ with $G$ a discrete amenable
group. Hence there exists a non-normal conditional expectation
from the discrete core $L(F_{\infty})\ten B(\ell_2(G))$ onto
$\M_1$. Since $L(F_{\infty})\ten B(\ell_2(G))$ is QWEP by
Kirchberg's \cite{Ki} results, we deduce that $\M_1$ is QWEP. For
general $\phi$, we use an ultraproduct argument in approximating
$\ln d$ by duadic rationals. \qd

\begin{theorem}\label{Ko} Let
$N$ be a semifinite von Neumann algebra with trace $\tau$ and $d$
a positive density in $L_2^{ap}(N,\tau)$. Then $K(N,d)$ is
$3$-completely complemented in the predual of a von Neumann $\M$
algebra. If $N$ is hyperfinite, then there is such an $\M$ with
QWEP.
\end{theorem}\lz

\begin{proof}[\bf Proof:] Since $d\in L_2^{ap}$ and in view of Lemma
\ref{dlimit1}, it suffices to show this for a positive density
$d\in L_2(N,\tau)$. By normalization we can assume $\tau(d^2)=1$
and defined the state $\varphi_d(x)=\tau(d^2x)$. According to
Lemma \ref{limit1}, $K(N,d)$ is completely contractively
complemented in $\prod_{n,\U} \kz_{n}(N,d)$. According to Theorem
\ref{comp2}, $\prod_{n,\U}\kz_n(N,d)$ is $3$  completely
contractively complemented in
 \begin{eqnarray}
  \M_* &=& \prodd_{n,\U} L_1(\M_n)
  \end{eqnarray}
for every free ultrafilter on the integers. The first assertion
is proved. Now, we assume in addition that $N$ is hyperfinite.
According to Lemma \ref{dykprep}, $K(N,d)$ is completely
contractively complemented in
 \[ \prodd_{\al,\U} K(A_\al,d_\al)  \pl ,\]
where the $A_\al$'s are  finite dimensional. Thus $\M_n(\al)
=\ast_{i=1}^n \ell_\infty^2(A_\al)$ is a free product of finite
dimensional von Neumann algebras. According to Lemma \ref{dyk},
all the $\M_n(\al)$'s  (and $\M_n^{op}(\al)$) are QWEP. Using
Lemma \ref{ultraqwep}, we deduce that $\prod_{n,\U}
L_1(M_n(\al))$ is a predual of a QWEP von Neumann algebra. Using
a further ultraproduct
 \[ \prodd_{\al,\U'} \prodd_{n,\U} L_1(\M_n(\al))\]
we stay in the class of preduals with QWEP (after another
application of Lemma \ref{ultraqwep}). \qd

\begin{cor} \label{kick1} Let $N$  be semifinite hyperfinte von
Neumann algebras with trace $\tau$ and $d$ be a density in
$L_2^{ap}(N,\tau)$. Then
 \[ K(N,d)\wet K(N,d) \]
$9$ completely complemented in
 \[ \M_*\wet \M_* \]
for some von Neumann algebra $\M$ with QWEP.
\end{cor}\lz

{\bf Proof of Proposition \ref{kick00} and Theorem 2:} By Lemma
\ref{2matrix} and Lemma \ref{G}, we see that $G$ is completely
contractively complemented in $K(M_2(\ell_\8(L^\8(\mu))),d)$ where
$d(t)=\kla
\begin{array}{cc} \sqrt{t} &0 \\ 0 &\sqrt{1-t} \end{array}\mer$.
According to Lemma \ref{dykprep}, $d\in L_2^{ap}$. Hence, the
assumptions of  Theorem \ref{Ko} and Corollary \ref{kick1} are
satisfied and Theorem 2 follows immediately. The lower estimate in
Proposition \ref{kick00} follows immediately by complementation.
\qed

\section{Appendix}

We will explain how the central limit procedure leads to  a
concrete realization of the von Neumann algebra $\M$ having $OH$
in its predual. This procedure relies on  Speicher's central
limit theorem. Our basic data are given by a $\si$-finite measure
space $(\Om,\Si,\mu)$, the commutative von Neumann algebra
$L_{\infty}(\mu)=L_\infty(\Om,\Si,\mu)$ and two densities
$d_1,d_2$. On $N=M_2(L_\infty(\mu))$ we define the density
 \[ d\lel \kla\begin{array}{cc}d_1& 0\\ 0& d_2 \end{array}\mer
\quad \mbox{and the weight} \quad
 \phi_d(x)\lel \intt tr(d(\om)x(\om)) d\mu(\om) \pl .\]

\begin{lemma}\label{grost} Assume in addition
$T=\phi_d(1)<\infty$ and let
$\phi_2(x_1,x_2)=\frac{1}{2T}(\phi_d(x_1)+\phi_d(x_2))$. Let
$\nen$ and $\M_n\lel \ast_{i=1}^n (\ell_\infty^2(N),\phi_2)$ and
$\phi_n=\ast_{i=1}^n \phi_2$ the product state with density $D\in
L_1(\M_n)$. Then $\M_n$ is QWEP and there exists a mappings
$u_n:L_\infty(\mu)\to \M_n$, $v_n:L_\infty(\mu)\to \M_n$  such
that
 \begin{align*}
 \lefteqn{\max\left\{\noo \summ_k b_k\ten u_n(f_k)\rrm_{\B\ten_{min}\M_n},
 \noo \summ_k b_k\ten v_n(f_k)\rrm_{\B\ten_{min}\M_n^{op}}
 \right\}\kl 3 }\\
 &\pl
 \max\left\{\sqrt{\frac{T}{n}} \noo \summ_k b_k \ten
 f_k\rrm_{\B\ten_{\min}L_\infty(\mu)},
 \noo \summ_{k,j} b_k^*b_j \intt
 \bar{f_k}f_jd_2d\mu\rrm_{\B}^{\frac12},
 \noo \summ_{k,j} b_kb_j^* \intt
 \bar{f_j}f_kd_1 d\mu\rrm_{\B}^{\frac12}\right\}
 \end{align*}
for all $C^*$-algebras $\B$ and $b_k\in \B$, $f_k\in
L_\infty(\mu)$. Moreover,
 \[ tr_{\M_n}(u_n(f)D_{\phi}^{\frac12}v_n(g)D_{\phi}^{\frac12})
 \lel \intt_{\Om} fg\sqrt{d_1d_2} d\mu \]
for all $f,g\in L_\infty(\mu)$.
\end{lemma}

\begin{proof}[\bf Proof:] The QWEP property follows from Lemma \ref{dyk}. We define $U_n:N\to \M_n$ by
 \[ U_n(x)\lel \sqrt{\frac{T}{n}}\summ_{k=1}^n \pi_k(\eps x) \pl
 .\]
Since $\pi_k(\eps x)$ and $\pi_j(\eps x)$ are freely independent
mean $0$ elements, we obtain by approximation with analytic
elements that
 \begin{align*}
 tr(U_n(x)D_{\phi}^{\frac12}U_n(y)D_{\phi}^{\frac12})
 &= \frac{T}{n} \summ_{k=1}^n \frac{1}{T}tr_{\ell_\infty^2(N)}(\eps xd^{\frac12}y \eps d^{\frac12}) \lel
 tr_N(xd^{\frac12}yd^{\frac12}) \pl .
 \end{align*}
We define
 \[ u_n(f)\lel U_n\kla \begin{array}{cc}0&f\\0&0\end{array}\mer \quad \mbox{and} \quad
  v_n(f)\lel U_n\kla \begin{array}{cc}0&0\\f&0\end{array}\mer \]
Then, we deduce
 \[ tr(u_n(f)D_{\phi}^{\frac12}v_n(g)D_{\phi}^{\frac12}) \lel
 \intt f\sqrt{d_2}g\sqrt{d_1} d\mu \pl . \]
The second assertion is  proved. In order to prove the first
assertion, we consider a von Neumann algebra   $\B$ and
 \[ x\lel \summ_{k=1}^m b_k\ten f_k \in \B\ten N \pl .\]
We define $E:\B\ten \ell_{\infty}^2(N)\to \B$ by $E(b\ten n)\lel
\phi_2(n)b$. According to Lemma \ref{matmod}, we know that the
algebras
 \[ A_k \lel id\ten \pi_k(\B\ten_{min}\ell_{\infty}^2(N)) \]
are freely  independent over $E$. Therefore, we map apply
Proposition \ref{Voi} and deduce
  \begin{align*}
 \lefteqn{   \noo id\ten u_n(x) \rrm_{\B\ten_{min}\M_n}
 \kl  3 \sqrt{\frac{T}{n}} \max\left\{
   \noo \kla \begin{array}{cc} 0&x\\ 0&0 \end{array}\mer \rrm_{\B\ten_{min}
   \ell_\infty^2(N)}, \right.}\\
   &\quad \left.  \sqrt{n} \noo E\bigg((id\ten u_n)(x)^* (id\ten
   u_n)(x) \bigg)\rrm_{\B}^{\frac12}, \sqrt{n} \noo E\bigg((id\ten u_n)(x)(id\ten
   u_n)(x)^*\bigg)\rrm_{\B}^{\frac12} \right\} \pl .
  \end{align*}
Using the multiplication with  a matrix unit, we get
 \[ \noo \kla \begin{array}{cc} 0&x\\ 0&0 \end{array}\mer \rrm_{\B\ten_{min}
   \ell_\infty^2(N)} \lel \noo x\rrm_{\B\ten_{min}L_\infty(\mu)}
   \pl .\]
It is easily checked that
 \begin{align*}
 E\bigg((id\ten u_n)(x)^* (id\ten
   u_n)(x) \bigg)
   &= \summ_{k,j} b_k^*b_j  \intt \bar{f_k}f_j \frac{d_2}{T} d\mu
 \end{align*}
and
 \begin{align*}
 E\bigg((id\ten u_n)(x) (id\ten
   u_n)(x)^* \bigg)
   &= \summ_{k,j} b_kb_j^*  \intt f_k\bar{f_j} \frac{d_1}{T} d\mu
 \end{align*}
This yields the first inequality for $\M_n$. Note that
$\M_n^{op}=\ast_{i=1}^n (\ell_\infty^2(N^{op}),\phi_2)$. We
denote by $m_1\circ m_2=m_2m_1$ the reversed multiplication. We
may apply Proposition \ref{Voi} to the same map $U_n$ defined on
$N^{op}$. Since $L_\infty(\mu)^{op}=L_{\infty}(\mu)$ the first
term is unchanged. For the second term, we note that
 \begin{align*}
 E\bigg((id\ten v_n)(x)^* \circ (id\ten
   v_n)(x) \bigg)
  &= \summ_{k,j} b_k^*b_j E\kla \kla \begin{array}{cc} 0&\bar{f}_k\\
  0 &0 \end{array}\mer \circ
  \kla \begin{array}{cc} 0&0\\
  f_j &0 \end{array}\mer\mer \\
  &= \summ_{k,j} b_k^*b_j \intt f_j\bar{f_k} \frac{d_2}{T} d\mu \pl .
 \end{align*}
By the commutativity of $L_{\infty}(\mu)$, we obtain the same
right hand side for $u_n$ and $v_n$. \qd

\begin{rem} {\rm Let us consider the space case $\Om=\{1,...,m\}$,
$\mu$ the counting measure and $d_1(j)=\la_j$,
$d_2(j)=\la_j^{-1}$. This yields $T=\sum_{j=1}^m
\la_j+\la_j^{-1}$ and we find map $u_n:\cz^m\to \M_n$,
$v_n:\cz^m\to \M_n^{op}$ such that
 \begin{align*}
 &\max\left\{\noo \summ_k b_k\ten u_n(e_k)\rrm_{\B\ten_{min}\M_n},
 \noo \summ_k b_k\ten v_n(e_k)\rrm_{\B\ten_{min}\M_n^{op}}
 \right\}\\
 &\quad \le 3 \max\left\{\sqrt{\frac{T}{n}} \sup_k \noo b_k\rrm
 \noo \summ_{k} b_k^*b_k \la_k^{-1}\rrm_{\B}^{\frac12},
 \noo \summ_{k} b_kb_k^* \la_k\rrm_{\B}^{\frac12}\right\}
 \end{align*}
holds for all $C^*$-algebras $\B$ and $b_k\in \B$. If $n\ge T\max
\la_k^{-1}$, the first term can be estimated by the second term.
Let $\pi_l:\M_n\to B(L_2(\M_n))$ the natural representation,
$J(x)=x^*$ the conjugation operator and $\pi_r(x)\lel
J\pi_l(x^*)J$ the representation of $\M_n^{op}$ in
$\pi_l(\M_n)'$. Then, we  deduce from $\sqrt{d_1d_2}=1$ that
 \[ (D_{\phi}^{\frac12},\pi_l(u_n(e_k))\pi_r(v_n(e_j))D_{\phi}^{\frac12})
 \lel tr(D_{\phi}^{\frac12}u_n(e_k)D_{\phi}^{\frac12}v_n(e_j))\lel
 \delta_{kj} \pl .\]
This is the key combinatorial estimate from \cite[Lemma 2.4]{PS}.
However, we also have to know that $\M_n$ is QWEP, see Lemma
\ref{dyk}. This provides an alternative (though similar) approach
to the key probabilistic tools based on the module version of
Voiculescu's inequality.}
\end{rem}

In order to find a concrete realization of the von Neumann algebra
supporting $OH$, we will use Speicher's central limit theorem
\cite[Theorem 5]{Sp} (adapted to our setting).\lz

\begin{theorem}\label{moments}[Speicher] Let $N$ be a von Neumann algebra,
$\varphi:N\to \cz$ a positive faithful state and $T=\varphi(1)$.
Let $\varphi_2(x_1,x_2)\lel \frac{1}{2T}\varphi(x_1,x_2)$ and
$\M_n=\ast_{i=1}^n (\ell_\infty^2(N),\varphi_2)$
 \[ U_n(x)\lel \sqrt{\frac{T}{n}} \summ_{k=1}^n \pi_k(\eps x) \pl. \]
Let $\phi_n$ be the corresponding state associated to the vacuum
vector. Then  for all $x_1,...x_m\in N$ and
 \[ \lim_n \phi_n(U_n(x_1)\cdots U_n(x_m))\lel \begin{cases} 0 & m \mbox{ odd} \\
 \summ_{(i_1,j_1),...(i_{m/2},j_{m/2})\in NCP_m} \prod\limits_{k=1}^{\frac m2}
 \varphi(x_{i_k}x_{j_k}) & m \mbox{ even} \end{cases}\pl .\]
Here $NCP_m$ stands for the non-crossing pair partitions of
$\{1,...,m\}$.
\end{theorem}\lz

Momentarily, we still assume $\phi_d(1)=T$. Then, we may pass to
the limit for $n\to \infty$ and deduce for a single self-adjoint
element $x\in N$ that
 \[ \lim_n \phi_n(U_n(x)^m)\lel \card(NCP_m) \pl  \phi_d(x^2)   \pl
 .\]
This means the limit object has the distribution of a
semi-circular random variable. In the previous section, we have
have already seen how to define this limit object as an operator.
We consider a free ultra-filter $\U$ and
$\tilde{\M}_\U=(\prod_{n,\U} (\M_n)_*)^*$. Let us denote
$\phi_\U=(\phi_n)$ and $e_{\U}$ the support of $\phi_{\U}$. Then
$\phi_{\U}$ defines a normal faithful state on
$\M_{\U}=e_{\U}\tilde{M}e_{\U}$. For a bounded element $x\in N$,
we know that $U_{\infty}(x)=(U_n(x))_n$ is a bounded sequence and
hence an element in $\tilde{M}_{\U}$. If $x$ is in addition
$\varphi_d$ analytic, we deduce from Raynaud's isomorphism
\cite{Ray} that
 \[ U_{\infty}(x)D_{\phi_\U}^{\frac12}\lel (U_n(x)D_{\phi_n}^{\frac12})_n \lel
 (D_{\phi_n}^{\frac12}U_n(\si_{\frac{i}{2}}^{\varphi_d}(x)))_n \lel
 D_{\phi_\U}^{\frac12}U_{\infty}(\si_{\frac{i}{2}}^{\varphi_d}(x)) \pl .\]
This implies $\phi_{\U}( U(x)^*(1-e)U(x))$ and hence
$(1-e)U(x)e=0$. Hence, for analytic elements $x\in N$, we deduce
that $U_{\infty}(x)=eU(x)e$ still satisfies
 \begin{equation}
 \phi_{\U}(U_{\infty}(x_1)\cdots U_{\infty}(x_m))\lel
  \summ_{(i_1,j_1),...(i_{m/2},j_{m/2})\in NCP_m} \prod\limits_{k=1}^{\frac m2}
 \varphi_d(x_{i_k}x_{j_k}) \label{mm} \end{equation}
for even $m$ and $0$ else. In particular, for selfadjoint $x$,
the element  $U_{\infty}(x)$ is semicircular and hence satisfies
(see e.g. \cite[Corollary 2.6.5]{VDN})
 \[ \noo U_{\infty}(x)\rrm\kl 2 \varphi_d(x^2)^{\frac12} \pl .\]
By continuity, we may now extend $U_{\infty}$ to a  continuous map
$U_{\infty}:L_2(N_{sa},\phi_d)\to \M_{\U}$ which still satisfies
the moment formula \eqref{mm}. In order to remove the assumption
$\phi_d(1)=T<\infty$, we follow the approach in Theorem \ref{Ko}
and first pick an increasing sequence $(f_{\al})$ in the
centralizer of $N$ converging to $1$ such that
$\phi(f_{\al})=T_{\al}<\infty$. Using the construction above, we
obtain $\M_{\al}$ and a map $U_{\al}:L_2(f_{\al}N_{sa},\phi_d)\to
\M_{\al}$. Then, we consider $\tilde{M}=(\prod_{\al,\U'}
(\M_{\al})_*)^*$ and the support $e$ of the ultraproduct state
$\phi=(\phi_{\al})_{\al}$. Using the same trick as above, we find
deduce that for analytic elements $(U_{\al}(x))_{\al}$ commutes
with $e$. We use $\M=e\tilde{M}e$ and obtain a map
$U_{\infty}:L_2(N_{sa},\phi_d)\to \M$ still satisfying
\eqref{mm}. We extend $U_{\infty}$ by linearity to complex linear
map. Since
 \[ \kla \begin{array}{cc} 0&f\\0&0 \end{array}  \mer \lel
 \frac{1}{2} \kla \begin{array}{cc} 0&f\\\bar{f}&0 \end{array}\mer  +
 \frac{-i}{2} \kla \begin{array}{cc} 0&if\\\overline{if}&0
 \end{array}\mer
 \]
we see that $U_{\infty}$ is well defined for matrices which have
only entries in the right upper or left lower corner. According
Theorem \ref{Ko} we obtain a completely contractive embedding
$u_1:K(N,\sqrt{d})\to L_1(\M)$ given by
 \begin{equation}
  u_1\bigg((D_{f_{\al_0}\phi_{d}}^{\frac12}x,yD_{f_{\al_0}\phi_{d}}^{\frac12})+kerQ\bigg)
  \lel D_{\phi}U_{\infty}(x)+U_{\infty}(y)D_{\phi}
   \label{emb}
 \end{equation}
such $\noo u\rrm_{cb}\le 1$. Since we are interested in
$K(L_\infty(\mu),\sqrt{d_1},\sqrt{d_2})$, we have to restrict
$u_1$ to matrices which only have an entry in the right upper
corner. It is also more appropriate to consider the slightly
smaller algebra $\M(\mu,d_1,d_2)$ generated by the image under
$U_{\infty}$ of
 \[ L^2\lel \left\{ \kla \begin{array}{cc}0&f\\\bar{f}&0 \end{array} \mer \mitt
 \intt |f|^2(d_1+d_2)d\mu<\infty\right\}\pl \subset L_2(N_{sa},\varphi_d) \pl .\]
Note that $\M(\mu,d_1,d_2)$ is complemented in $\M$ because  we
have
 \begin{equation}\label{moddd}  \si_t^{\phi}(U_{\infty}(x))\lel
 U_{\infty}(\si_t^{\varphi_d}(x)) \end{equation}
for all $x\in N$. Using complex linear combinations of
selfadjoint elements, we see that $u_1$ maps into
$L_1(\M(\mu,d_1,d_2))$. We may apply Lemma \ref{grost} and deduce
that the map $u_{\infty}(f)=U_{\infty}\kla
\begin{array}{cc}0&f\\0&0\end{array}\mer$ provides a completely
bounded map on $K(L_{\infty}(\mu),\sqrt{d_1},\sqrt{d_2})^*$ with
$cb$-norm $\le 2$. (Here we also used the density of
$L_{\infty}(\mu)$ in $K(L_{\infty}(\mu),\sqrt{d_1},\sqrt{d_2})^*$
see Lemma \ref{wstar}). Thus
$K(L_{\infty}(\mu),\sqrt{d_1},\sqrt{d_2})^*$ is $2$-cb
complemented in $\M(\mu,d_1,d_2)$.

The advantage of this smaller algebra is also given by the fact
that unitaries (contraction) on $L_2(\mu)$ extend to automorphism
on $\M(\mu,d_1,d_2)$. More precisely, let us assume that a unitary
$w:L_2(\mu)\to L_2(\mu)$ commutes with the canonical
multiplication operators $M_{d_1}$, $M_{d_2}$ associated with
$d_1$, $d_2$, respectively. Then we define $W:L^2_{sa}\to
L^2_{sa}$ as follows
 \begin{equation} \label{wform}
 W\kla\begin{array}{cc} 0& f\\\bar{f}&0 \end{array}\mer
 \lel \kla\begin{array}{cc} 0& w(f)\\\overline{w(f)}&0
 \end{array}\mer \pl .\end{equation}
Since $w$ is complex linear, we find that $W(x)$ is selfadjoint if
$x$ is.  We also have
 \[ \phi_d(W(x)W(y))\lel \phi_d(xy) \pl \]
for all $x,y\in L^2_{sa}$ because $w$ commutes with $d_1$ and
$d_2$.  In view of \eqref{mm}, we see that
 \[ \pi_w(U_{\infty}(x))\lel U_{\infty}(W(x)) \quad , \quad x\in L^2_{sa} \]
extends to a $\phi$-preserving automorphism of $\M(\mu,d_1,d_2)$.
Note also that in view of \eqref{moddd} and
 \[ W\circ \si_t^{\varphi_d}\lel \si_t^{\varphi_d} \circ W \pl ,\]
$\pi_w$ commutes with $\si_t^{\phi}$. For a general contraction
$v$, we work in $L_2(\mu)\oplus L_2(\mu)=L_2(\mu;\ell_2^2)$ and
define the unitary
 \[ w(v)\lel \kla\begin{array}{cc} v & \sqrt{1-vv^*}\\
                                -\sqrt{1-v^*v} & v^*
                                \end{array}\mer \pl .\]
Using the densities $d_1\ten 1$, $d_2\ten 1$, it suffices to note
that the canonical image of $L^2\subset L^2\oplus L^2$ is
$\si_t^{\varphi_d}$ invariant. Hence, $\pi(u(x))\lel u(x,0)$
extends to a $*$-preserving homomorphism onto a $\si_t^{\phi}$
invariant subalgebra of $\M_{d\ten 1}$. Let $E$ be the canonical
conditional expectation, we see that
 \[ \pi_v\lel E\circ \pi_{w(v)}\pi\]
is a completely positive map such that $\phi\circ \pi_v=\phi$ and
 \[ \pi_v(u_{\infty}(f))\lel u_{\infty}(vf) \pl .\]
This implies $\pi_{v_1v_2}\lel \pi_{v_1}\pi_{v_2}$. Using the
moment formula it is also clear that the map $\pi(v)=\pi_v$ is
continuous with  respect to the strong topology.

We may apply this in particular to a projection $p\in L_2(\mu)$
which commutes with $d_1$ and $d_2$ and denote by $P:L^2_{sa}\to
L^2_{sa}$ the induced projection. Using Speicher's moment formula
\ref{moments} and the definitions of free cumulants, we deduce
that
 \[ k_m[U(x_1),....,U(x_m)]  \lel \begin{cases} 0 & \mbox{ if } m>2 \\
 \phi_d(x_1x_2) & \mbox{ if } m=2 \end{cases} \pl.
 \]
In particular,  it follows that $U_{\infty}(P(L^2))$ and
$U_{\infty}((1-P)(L^2))$ are $^*$-free (see \cite[section
1.4]{Sp2} and the references therein) with respect to $\phi$.
Therefore $\M(\mu,d_1,d_2)$ is the free product of
$\pi_p(\M(\mu,d_1,d_2))$ and $\pi_{1-p}(\M(\mu,d_1,d_2))$.

Resuming our results, we obtain
 \begin{enumerate}
 \item[i)]  There exists a (weak$^*$-continuous) map
  \[ u_{\infty}:K(L_{\infty}(\mu),\sqrt{d_1},\sqrt{d_2})^*\to  \M(\mu,d_1,d_2)
  \]
such that $u_1^*u_{\infty}=id$ and
  \[ \noo u_{\infty}\rrm_{cb}\noo u_1^*:\M(\mu,d_1,d_2) \to K(L_{\infty}(\mu),\sqrt{d_1},\sqrt{d_2})^*\rrm_{cb}
  \le 2 \]
and
  \[ tr(u_{\infty}(f)^*D_{\phi}^{\frac12}u_{\infty}(g)D^{\frac12}_{\phi})\lel
  \intt \bar{f}g \sqrt{d_1d_2} d\mu \pl. \]
\item[ii)] If $w:L_2(\mu)\to L_2(\mu)$ is contraction
which commutes with the multiplication map $M_{d_2}$ and
$M_{d_1}$, then there exits a $\phi$-preserving normal completely
positive map $\al_w:\M\to \M$ such that
 \[ u_{\infty}\circ w \lel \al_{w}\circ u_{\infty} \quad \mbox{and}
 \quad u_{\infty}\circ w \lel \al_{w}\circ u_{\infty}
 \pl. \]
If $u$ is unitary $\al_{w}$ is $^*$-homomorphism.
 \item[iv)] If $p$ is a projection in $B(L_2(\mu))$ which commutes with $d_1$ and $d_2$, then there exists freely independent
subalgebras $\M_p$, $\M_{1-p}\subset \M(\mu,d_1,d_2)$ such that
 \[ (\M(\mu,d_1,d_2),\phi)\lel (\M_p,\phi)\ast (\M_{1-p},\phi) \pl
 \]
and $u(pL_2((d_1+d_2)\mu)\subset \M_p$.
\end{enumerate}

It is now very easy to provide a concrete construction of
$\M(\mu,d_1,d_2)$. Indeed, we may use the Fock space
 \[ \F(L_2(\rz;L^2_{sa}))\lel \summ_{n\ge 0} \oplus (L_2(\rz;L^2_{sa}))^{\ten_n} \]
and then for a selfadjoint  $x=\kla
\begin{array}{cc}0&f\\ \bar{f}&0
\end{array}\mer$
 \[ u_{\F}(1_{[0,t]}x)(h_1\ten \cdots h_n)\lel 1_{[0,t]}x \ten h_1\ten \cdots \ten h_n
 +(1_{[0,t]}x,h_1)\pl  h_2\ten \cdots \ten h_n \pl .\]
Of course the scalar product $(1_{[0,t]}x,h)$ depends on
$d_1,d_2$ and $\mu$, namely for $h=\kla
\begin{array}{cc} 0& g\\ \bar{g} &0
\end{array}\mer$ we have
 \[ (1_{[0,t]}x,1_{[0,t]}h)\lel t\intt f\bar{g}d_1d\mu + t\intt
 \bar{f}gd_2d\mu \pl .\]
Let us denote by $\Omega$ the vacuum vector corresponding to
$n=0$. According to Speicher's central limit theorem
\cite[Theorem5]{Sp}, we have
 \begin{align*}
  (\Om,u_{\F}(x_1)\cdots u_{\F}(x_m)\Om)
  &= \summ_{r=0}^m \summ_{\{V_1,...,V_r\}\in NCP_m} Q(V_1)\cdots
  Q(V_r)
  \end{align*}
and
 \[ Q(\{i_1,...,i_l\}) \lel \lim_{N\to \infty}
 N(\Om,u(1_{[0,\frac{1}{N}]}x_{i_1})\cdots
 u(1_{[0,\frac{1}{N}]}x_{i_l})\Om) \]
Following \cite[Theorem6]{Sp}, we see that only pair partitions
provide a non-zero term and in this case
 \[ Q((i,j))\lel (\Om,u_{\F}(1_{[0,1]}x_i)u_{\F}(1_{[0,1]}x_j)\Om)\lel (x_i,x_j) \pl .\]
Now, we don't need $L_2(\rz)$ anymore and may work on
$\F(L^2_{sa})$. If $\rho(x)=(\Om,x\Om)$ denotes the state
corresponding to the vacuum vector, we deduce
 \begin{align*}
 (\Om,u_{\F}(x_1)\cdots u_{\F}(x_m)\Om)&=
 \summ_{(i_1,j_1),...(i_{m/2},j_{m/2})\in NCP_m} \prod\limits_{k=1}^{\frac m2}
 (x_{i_k},x_{j_k})\\
  &= \phi(U_{\infty}(x_1)\cdots U_{\infty}(x_m)) \pl .
  \end{align*}
We deduce that $\M$ is spatially isomorphic to $\{u_\F(x)|x\in
L^2_{sa}\}''\subset B(\F(L^2_{sa}))$. This formula also provides
the crucial link to the free quasi-free states in free
probability discussed by Shlyakhtenko \cite[Remark 2.6]{S1}.
Indeed, the algebra $\Gamma(\Hh_{\rz},\U_t)$ is generated by
semicircular elements $\{s(h)|h\in \Hh_{\rz}\}$ satisfying
 \[ \phi_U(s(h_1)\cdots s(h_m))\lel 2^{-m}
 \summ_{(i_1,j_1),...(i_{m/2},j_{m/2})\in NCP_m} \prod\limits_{k=1}^{\frac m2}
 (x_{i_k},x_{j_k})_U \pl .\]
Here $(\pl ,\pl)_U$ is a complex sesquilinear form given by
 \[ (x,x)_U \lel (\frac{2}{1+A^{-1}}x,x)\pl ,\]
where $U_t=A^{it}$ is a one parameter family of unitaries acting
on the complex Hilbert space $\Hh_\cz=\Hh_\rz+i\Hh_{\rz}$ and
leaves $\Hh_{\rz}$ in variant. In our case, we have
 \[ \kla \kla \begin{array}{cc}0&f\\\bar{f}&0\end{array}\mer,
 \kla\begin{array}{cc}0&g\\\bar{g}&0\end{array}\mer\mer_U
 \lel 2\intt f\bar{g}d_1d\mu + 2\intt
  \bar{f}gd_2d\mu \pl .\]
Calculating the real part, we find the real scalar product is
given by
 \[ (x,x)_{\Hh_{\rz}}\lel 2 \intt (\bar{f}g+\bar{g}f) (d_1+d_2)
 \pl d\mu \pl .\]
However, by \cite[(2)]{S1}, we have
 \[ (x,y)_U\lel (y,A^{-1}x)_U \pl .\]
This implies
 \[ A^{-1}\kla \begin{array}{cc}0&f\\\bar{f}&0\end{array}\mer
 \lel  \kla \begin{array}{cc}0& \frac{d_1}{d_2}f\\ \frac{d_2}{d_1} \bar{f}&0\end{array}\mer   \pl .\]
We deduce
 \[ U_t\kla \begin{array}{cc}0&f_1\\ f_2&0\end{array}\mer
  \lel \kla \begin{array}{cc}0& d_2^{it}d_1^{-it}f_1 \\
  d_1^{it}d_2^{-it}f_2&0\end{array}\mer   \pl .\]
Using Shlayhtenko's result   \cite[Theorem 6.10]{S1}, we obtain
the following application.\lz

\begin{cor}\label{typeIII} $\M(\mu,d_1,d_2)$ is isomorphisc to
$\Gamma(L^2_{sa},U_t)$. If the multiplication operators associated
to $d_1/d_2$ and $d_2/d_1$ have no point spectrum, then
$\M(\mu,d_1,d_2)$ is type $III_1$ factor.
\end{cor}\lz

\noindent Let us now describe how this applies to the operator
space $OH$. We  recall
 \[ d\mu(t)\lel \frac{1}{\pi\sqrt{t(1-t)}} \quad \mbox{and}\quad
 d_1(t)\lel t\pll , \pll  d_2(t)=1-t \pl .\]
We consider $\ell_\infty(L_\infty([0,1],\mu))$ with the product
measure $\mu_{\infty}=m\ten \mu$, $m$ the counting measure. The
space $G$ is the quotient of $L_2^c((1\ten d_1^{-1})
\mu_\infty)\oplus_1 L_2^r((1\ten d_2^{-1})\mu_\infty)$ with
respect to subspace $E$ spanned by elements $(f,-f)$. In view of
Lemma \ref{G} and \eqref{emb}, we  see that  the mapping
 \[ \hat{u}_1\circ (I_1,I_2)(f,g)\lel
  \hat{u}_1(fd_1^{-\frac12},gd_2^{-\frac12})\lel
 D_{\phi}u_{\infty}(fd_1^{-1})+
 u_{\infty}(gd_2^{-1}) D_{\phi}   \]
produces a $2$-cb embedding $u_1$ of  $G$. We may apply this to
the unit vectors $f_k\in G$ and find
 \begin{align*}
  b_k &=  u_1(f_k) \lel
 D_{\phi}u_{\infty}(e_k1_{[\frac12,1]}d_1^{-1})+
 u_{\infty}(e_k1_{[0,\frac12]}d_2^{-1})D_{\phi} \pl .
 \end{align*}
Note again that $u_{\infty}(e_k1_{[0,\frac12]}d_2^{-1})$ and
$u_{\infty}(1_{[\frac12,1]}d_1^{-1})$ are linear combinations of
semicircular random variables and $b_k$ is the restriction of the
sum of vector states
 \[ \psi_k(T) \lel ((u_{\infty}(e_k1_{[\frac12,1]}d_1^{-1}))^*\Om,
 T\Om) + (\Om,T
 u_{\infty}(e_k1_{[0,\frac12]}d_2^{-1})\Om)  \pl .\]
According to condition iii) above,
$\M_{\infty}=\M(\mu_{\infty},1\ten d_1,1\ten d_2)$ is  an
infinite free product of $\M(\mu,d_1,d_2)$
 \[ (\M_\infty,\phi) \lel \ast_{i=1}^\infty (\M(\mu,d_1,d_2), \phi) \pl
 .
 \]
Indeed, the moment formula \ref{mm} tells us that the algebra
$A_k=\pi_{e_k}\M_{\infty}$ are all isomorphic to
$\M(\mu,d_1,d_2)$, free from each other and their union generates
$\M_{\infty}$. This means that the $b_k$'s  are `free copies' (in
the sense of $L_1$) of $b_1$. (Note also that the construction of
$\M(\mu,d_1,d_2)$ only requires one limit procedure because
$\phi_d$ is a state in that case. In view of Corollary
\ref{typeIII}, we know that $\M(\mu,d_1,d_2)$ and $\M_{\infty}$
are type III$_1$. Using condition iv), we see that there is a
$^*$-preserving homomorphism $\rho:B(\ell_2)\to
CP_{\phi}^{\si}(\M_{1\ten d})$, the normal completely positive
$\phi$-preserving maps on $\M_{1\ten d}$, given by
 \[ \rho(v)\lel \pi_{v\ten id_{L_2(\mu)}} \pl .\]
Using the definition of $\pi_v$ it is easily checked that $\rho$
is continuous with respect to the strong topology and weak
topology on bounded sets. In particular, if $v_n$ converges to
$0$ weakly, $\rho(v_n)$ converges to $\rho(0)(x)=\phi(x)$. This
provides an alternative way  of showing that $\M_\infty$ is a
factor and that the only space invariant element is the identity.
It is easily checked that
 \[ \pi_v^*(u_{\infty}(f)D_{\phi})\lel u_{\infty}(v^*f)D_{\phi}
 \pl .\]
In particular, the space $F_{oh}={\rm span}\{b_k|k\in \nz\}$ is
homogeneous.\lz

\begin{prob} Is $F_{oh}$ competely isometrically isomorphic to
$OH$? \end{prob}\lz

\begin{prob} Are $\M_{\infty}$ and $M(\mu,d_1,d_2)$ isomorphic?
\end{prob}
\lz

The same arguments apply to $C_p=[R,C]_{\frac1p}$ using
$\mu_{\al}$ and the same densities $d_1$ and $d_2$. A posteriori,
one can also show directly that the image of
$u_{\infty}K(N,\sqrt{d_1},\sqrt{d_2})^*$ is completely
complemented and therefore obtain an alternative prove for the
embedding of $OH$ and in general $[C,R]_{\frac1p}$.

\newcommand{ \pp } {\prod _\U}
\newcommand{ \flip } {\rm flip}
\newcommand{ \ee } {\eta}

\end{document}